\documentclass{article}
\usepackage{a4wide}
\usepackage{authblk}
\usepackage[letterpaper,top=2cm,bottom=2cm,left=1cm,right=1cm,marginparwidth=1.75cm]{geometry}

\usepackage{hyperref} 

\usepackage{lineno}
\usepackage{microtype}
\usepackage{subcaption}
\usepackage[onehalfspacing]{setspace}
\usepackage{epsfig}

\usepackage{tikz}
\usetikzlibrary{patterns}
\usepackage{pgfplots}
\pgfplotsset{compat=newest}
\usepackage{comment}
\usepackage{amssymb}
\usepackage{amsmath}
\usepackage[noabbrev,capitalise,english]{cleveref} 

\usepackage{color}
\definecolor{ForestGreen}{RGB}{34,139,34}

\usepackage{bm}

\newcommand{\revision}[1]{\textcolor{black}{#1}}
\newcommand{\franz}[1]{\textcolor{black}{#1}}

\newcommand{\twick}[1]{\textcolor{black}{#1}}

\newcommand{\tom}[1]{\textcolor{black}{#1}}

\newcommand{\FL}{{\mathcal F}}
\newcommand{\SO}{{\mathcal S}}
\newcommand{\IN}{{\mathcal I}}
\newcommand{\bSI}{{\boldsymbol{\Sigma}}}
\newcommand{\E}{\mathbf{E}}
\renewcommand{\L}{\mathbf{L}}
\newcommand{\tr}{\operatorname{tr}}

\usepackage{algorithm, algorithmicx, algpseudocode}

\usetikzlibrary{shapes}
\usetikzlibrary{plotmarks}
\usepackage{float}

\newcommand{\eps}{\varepsilon}

\newcommand{\ds}{\displaystyle}

\let\vv =\v

\newcommand{\be}{\begin{equation}}
\newcommand{\ee}{\end{equation}}
\newcommand{\beq}{\begin{eqnarray}}
\newcommand{\eeq}{\end{eqnarray}}
\newcommand{\beqs}{\begin{eqnarray*}}
\newcommand{\eeqs}{\end{eqnarray*}}
\newcommand{\et}{\end{theorem}}
\newcommand{\bex}{\begin{example}}
\newcommand{\eex}{\end{example}}
\newcommand{\br}{\begin{remark}}
\newcommand{\er}{\end{remark}}
\newcommand{\bc}{\begin{corollary}}
\newcommand{\ec}{\end{corollary}}
\newcommand{\bl}{\begin{lemma}}
\newcommand{\el}{\end{lemma}}
\newcommand{\bp}{\begin{proposition}}
\newcommand{\ep}{\end{proposition}}
\newcommand{\bd}{\begin{definition}}
\newcommand{\ed}{\end{definition}}
\newcommand{\bas}{\begin{assumption}}
\newcommand{\eas}{\end{assumption}}

\newtheorem{remark}{Remark}

\newcommand{\R}{\mathbb{R}}

\newcommand{\PP}{\mathbb{P}}

\newcommand{\node}  {{\cal N}_h}        

\newcommand{\caE}{{\cal E}}

\newcommand{\caK}{{\cal K}}

\usepackage{stmaryrd}
\def\dd{{\rm d}}
\newcommand{\jumpEK}[1]  {\llbracket#1\rrbracket_{E,K} }  

\def\eb{{\bf e}}
\def\u{{\bf u}}
\def\U{{\bf U}}

\def\n{{\bf n}}

\def\f{{\bf f}}

\def\vv{{\bf v}}
\def\w{{\bf w}}

\def\z{{\bf z}}

\def\x{{\bf x}}
\def\y{{\bf y}}
\def\V{{\bf V}}
\def\W{{\bf W}}

\def\F{{\bf F}}
\def\0{{\bf 0}}

\newcommand{\bsigma}{\bm{\sigma}}

\newcommand{\bPi}{\bm{\Pi}}

\newcommand{\bB}{\bm{B}}
\newcommand{\bT}{\bm{T}}
\newcommand{\bzero}{\bm{0}}
\newcommand{\bu}{\bm{u}}
\newcommand{\bpp}{\bm{p}}

\newcommand{\bE}{\bm{E}}
\newcommand{\bS}{\bm{S}}
\newcommand{\bC}{\bm{C}}
\newcommand{\bI}{\bm{I}}
\newcommand{\bF}{\bm{F}}
\newcommand{\bV}{\bm{V}}
\newcommand{\bQ}{\bm{Q}}
\newcommand{\bq}{\bm{q}}
\newcommand{\bv}{\bm{v}}
\newcommand{\bw}{\bm{w}}
\newcommand{\bz}{\bm{z}}
\newcommand{\bR}{\bm{R}}
\newcommand{\bJ}{\bm{J}}

\newcommand{\nEK}{\n_{E,K}}

\newcommand{\bnabla}{\mbox{\boldmath{$\nabla$}}}

\newcommand{\bs}{\mbox{\boldmath{$\sigma$}}}
\newcommand{\e}{\mbox{\boldmath{$\varepsilon$}}}

\renewcommand{\div} {{\rm{div} \,}}

\newcommand{\restrictiona}[1]

\newcommand{\wzh}{\widehat{\mathbf{z}}_{h}}

\renewcommand{\vec}{\bm}

\newcommand{\eg}{e.g.} 

\usetikzlibrary{shapes.geometric, arrows}

\tikzstyle{startstop} = [rectangle, rounded corners, minimum width=3cm, minimum height=1cm,text centered, draw=black, fill=red!30]
\tikzstyle{process} = [rectangle, minimum width=3cm, minimum height=1cm, text centered, draw=black, fill=orange!30]
\tikzstyle{decision} = [diamond, minimum width=3cm, minimum height=1cm, text centered, draw=black, fill=yellow!30]
\tikzstyle{arrow} = [thick,->,>=stealth]

\usepackage[english]{babel}

\begin{document}

\author[1]{Roland Becker}
\author[2]{Franz Chouly\footnote{\it Corresponding author: 
\url{fchouly@cmat.edu.uy}.}}
\author[3]{Michel Duprez}
\author[4]{Thomas Richter}
\author[5]{Pierre-Yves Rohan}
\author[6]{Thomas Wick}

\affil[1]{Department of Mathematics, Université de Pau et de l’Adour (UPPA), Avenue de l’Université, BP 1155, 64013 Pau CEDEX, France}
\affil[2]{Universidad de la República, Facultad de Ciencias, Centro de Matem\'atica, 11400 Montevideo, Uruguay}
\affil[3]{Team-project MIMESIS, Inria de l'Université de Lorraine, MLMS team, Icube, Universit\'e de Strasbourg, CNRS UMR 7357, 2 Rue Marie Hamm, 67000 Strasbourg, France}
\affil[4]{Otto von Guericke University Magdeburg, Institute for Analysis and Numerics, Magdeburg, Germany}
\affil[5]{Arts et Métier Institute of Technology, Institut de Biomécanique Humaine Georges Charpak, Université Sorbonne Paris Nord, Paris, France}
\affil[6]{Leibniz University Hannover, Institute of Applied Mathematics, Welfengarten 1, 30167 Hannover, Germany}


\title{Dual Weighted Residual-driven adaptive mesh refinement to enhance biomechanical simulations}


\maketitle

 \abstract{
\franz{%
This chapter describes how \emph{a posteriori} error estimates targeting a user-defined quantity of interest, using the Dual Weighted Residual (DWR) technique, can be easily applied for biomechanical simulations in current engineering practice. The proposed method considers a very general setting that encompasses complex geometries, model non-linearities (hyperelasticity, fluid-structure interaction) and multi-goal oriented techniques. The developments are substantiated with some numerical tests.} 
 }



\section{Introduction}
\label{sec:Introduction}

\franz{The importance of finite element analyses (FEA) for biomechanical investigations has increased considerably worldwide in recent years. 
A survey of applications using simulation modeling for the healthcare sector can be found for instance in \cite{Mielczarek2012}.}

In this context{,} one major issue is meshing, since the reliability of the predicted {mechanical response} arising from computer simulation {heavily relies} on the quality of the underlying finite element mesh
{\cite{gratsch_posteriori_2005}}.
The patient-specific mesh has to be built from segmented medical images (CT, MRI, ultra-sound){,} and {has} to conform to anatomical details with potentially complex topologies and geometries {\cite{bijar2016atlas}}. 
In general, the quality of a given mesh is assessed through purely geometrical criteria, that allow in some way to quantify the distortion of the geometry of the elements 
\cite{bucki-2011}.
Beyond mesh quality, mesh density is another, related, parameter that must be controlled during biomechanics simulations. 
Moreover, solutions must be obtained on commodity hardware within clinical time scales: milliseconds (for surgical training);  minutes (for surgical assistance); hours (for surgical planning).
\franz{To limit the impact of the discretization error, a common practice is to carry out a "convergence study", where the initial mesh is refined uniformly until the quantity of interest to the practitionner stops fluctuating.}

\franz{In this chapter, 
we investigate the capability of
\emph{a posteriori} error estimates \cite{ainsworth-oden-2000,verfurth-2013} to provide useful information about the {discretization error}, \emph{i.e.},
the difference between the finite element solution and the exact solution of the same boundary value problem on the same geometry.
\emph{A posteriori} error estimates are quantities computed 
from the numerical solution, that indicate the magnitude of the local error. These estimates are at the core of mesh adaptive techniques \cite{nochetto-2009}. 
Many {\it a posteriori} error estimation methods have been developed in the numerical analysis community. These methods have different theoretical and practical properties. However, despite their great potential, 
error estimates have rarely been considered for patient-specific finite element simulations {in the biomechanical community}.} 

To the best of our knowledge, the existing works that address this issue are rare: 
see for instance \cite{Od18} in the context of tumor growth
or
\cite{Bui2017_error_brain_ijnmbe,Bui2016_TBME}
in the context of neurosurgery.
In the two aforementiond works, the authors
study the discretization error (based on the energy norm) of real time simulations using the recovery-based 
technique of Zienkiewicz and Zhu \cite{zienkiewicz1992superconvergent}. 
This approach is inexpensive and allows for the dealing of real time simulations. However, the error in the energy norm might not provide useful information for applications where one is interested in the error of a real physical quantity of interest. To overcome this difficulty, estimates based on duality arguments are common for {\it a posteriori} error estimation, see \eg{} \cite{becker-rannacher-2001,endtmayer2018two,endtmayer2021,Heuveline2003,paraschivoiu-peraire-patera-1997,prudhomme-oden-1999,maday-patera-2000,giles-suli-2002,suttmeier2004,RiWi15_dwr,wick2016goal,wick2017pum,wick2024}.
\franz{This chapter builds up particularly on  two studies that have been carried out using the Dual Weighted Residuals (DWR) method, as presented in \cite{becker-rannacher-2001}: a study in small strain elasticity, see
\cite{duprez2020}, and another one in large strain elasticity, see \cite{bui:hal-04208610}.}
The reader can also refer to  \cite{larsson-2002,whiteley-tavener-2014,gonzalezestrada2014,granzow2018} for previous applications of goal oriented error estimation in non-linear elasticity. 
\revision{It also addresses new aspects such as 
fluid-structure interaction~\cite{Richter2017} 
and multi-goal error control and adaptivity~\cite{EnLaWi18}.}

This chapter is organized as follows. In Section~\ref{sec:method1} and Section~\ref{sec:Methods}, we describe the linear elastic and hyperelastic settings for soft tissue, the corresponding finite element discretization, the DWR \textit{a posteriori} error estimation as well as the algorithms for mesh refinement. 
Section~\ref{sect:fluidstructure} presents techniques for fluid-structure interaction.
In Section~\ref{sec:Results}, we illustrate  and validate the methodology for different test cases. 
The results are discussed in  Section~\ref{sec:discussion} and perspectives are provided in 
Section~\ref{sec:Conclusions}.

\section{\franz{Small strain elasticity}}\label{sec:method1}

\franz{To emphasize the main ideas of the method and to avoid technical difficulties, we start by describing it in the context of small strain elasticity and a linear goal functional. Next sections will be devoted to more general cases.}

\subsection{Setting}\label{sec:toy problem}

\franz{
Let us denote by $\Omega$ in $\R^d$, $d=2,3$, the domain for an elastic body in reference configuration. 
{We assume small deformations, combined with the plain strain assumption when $d=2$}. 
The boundary $\partial\Omega$ of the elastic body
consists {of} two {disjoint} 
parts $\Gamma_D$ and $\Gamma_N$, with meas$(\Gamma_D) > 0$.
The unit outward normal vector on the boundary
$\partial\Omega$ is denoted {by} $\n$.
A displacement $\u_D=\0$ is applied on $\Gamma_D${, and the} body is subjected {to volume forces $\bB \in
L^2(\Omega;\R^d)$ and} surface loads $\bB_N \in L^2(\Gamma_N;\R^d)$.
The 
virtual work of external loads in the body and on its surface is
\begin{equation*} 
\ell_E(\w) := \int_{\Omega} \bB \cdot \w ~\dd\x + \int_{\Gamma_{N}} \bB_N \cdot \w ~\dd\mathrm{s}.
\end{equation*}
For two displacement fields $\vv,\w :  \Omega \rightarrow \R^d$, we introduce 
the (internal) virtual work associated to passive elastic properties:
\begin{equation*} 
a(\vv,\w) := \int_{\Omega} \bs (\vv) :\e (\w) ~\dd\x. 
\end{equation*}
{The notation $\e(\vv) = \frac12 (\bnabla \vv + \bnabla \vv^{^T})$ represents the linearized strain tensor field, and
$\bs = (\sigma_{ij})$, $1 \le i,j \le d$, stands for the Cauchy stress tensor field. 
} 
}
%
\franz{Finally we take into account genuinely the active properties of soft tissue as a linear anisotropic pre-stress which reads 
\begin{equation*}\label{eq:defla}
\ell_A(\w) := - \beta T \int_{\omega_A} \left ( \e (\w) \eb_A \right ) \cdot \eb_A  ~\dd\x,
\end{equation*}
where $\omega_A$ is the part of the body where muscle fibers are active, $T \geq 0$ is a scalar which stands for the tension of the fibers, $\eb_A$ is a field of unitary vectors that stands for muscle fibers orientation,
{and} {$\beta \in [0,1]$} is the activation parameter. When $\beta=0$ there is no activation of the muscle fibers, and the value $\beta=1$ corresponds to the maximum of activation. This modelling can be viewed as a linearization of some more sophisticated active stress models of contractile tissues in large strain (see, \emph{e.g.}, {\cite{cowin2001cardiovascular,payan-ohayon-2017}}).}

\noindent The small strain elasticity problem reads:
\begin{equation}\left\{\begin{array}{l}
\textrm{Find a displacement }\u \in \V\textrm{ such that}\\\noalign{\smallskip}
a(\u, \vv) = \ell(\vv), \quad \forall \, \vv \in \V,
\end{array}\right.\label{eq:primal_weak}
\end{equation}
where $\ell(\cdot) = \ell_E(\cdot) + \ell_A(\cdot)$, and 
where 
$\u$ and 
$\vv$ {lie in} the space of admissible displacements
\begin{align*}
\V := \big \{ \vv\in H^{1}(\Omega;\R^d) \,|\, ~\vv = \0 \text{~on~} \Gamma_{D} \big \}.
%
\end{align*}
%
From the displacement field, we are interested in computing a linear quantity
\begin{equation}\label{def:J}
J: \V \ni \u \mapsto J(\u) \in \R, 
\end{equation}
which can be defined according to a specific application and  the interest of each practitioner. {Thereby, the quantity $J$ will be called {\emph{quantity of interest (QoI)}}.} 
\twick{Typical QoIs are point evaluations (if solution is sufficiently regular), line integrals 
over (parts of) the boundary or (sub-)domain integrals. Physically, they can represent deflections, stresses or mean values of the solutions.}
We will provide \franz{various expressions in } Section \ref{sec:Results}.

\subsection{\franz{Finite element approximation}}\label{sec:finite element method}

Consider a family of meshes $(\caK_h)_{h>0}$ 
constituted of triangles and assumed to be subordinate to the decomposition of the boundary 
$\partial\Omega$ into $\Gamma_D$ {and} $\Gamma_N$. 
%
For a mesh $\caK_h$, {we} denote by $\caE_h$ the set of \franz{edges/faces}, 
{by} $\caE^{int}_h:=\{E\in \caE_h:E\subset \Omega\}$ the set of interior \franz{edges/faces},
and {by} $\caE^N_h:=\{E\in \caE_h:E\subset \Gamma_N\}$ the set of boundary \franz{edges/faces} that correspond to Neumann conditions {(we assume that any boundary \franz{edge/face} is either inside $\Gamma_N$ or inside $\Gamma_D$)}.
For an element $K$ of $\mathcal{K}_h$, we set $\caE_K$ the set of \franz{edges/faces} of $K$, 
$\caE_K^{int}:=\caE_K\cap \caE_h^{int}$ and  $\caE_K^{N}:=\caE_K\cap \caE_h^{N}$. {We also assume that each element $K$ is either completely inside 
$\omega_A$ or completely outside it.}
Let $\bs$ be a second-order tensorial field in $\Omega$, {which is assumed to be} piecewise continuous.
We define the jump of $\bs$ across an interior edge $E$ 
of an element $K$, 
at a point $\y\in E$, as follow{s}
\begin{equation*}
\jumpEK{\bs} (\y):=\lim\limits_{\alpha\rightarrow 0^+}
\left (\bs(\y+\alpha\nEK) - \bs(\y-\alpha\nEK) \right)  \nEK,
\end{equation*}
where $\nEK$ is the unit normal vector to $E$, pointing out of $K$.

The finite element space $\V_{h}\subset\V$
is built upon
continuous Lagrange finite elements of degree $k=1,2$ 
(see, \emph{e.g.}, \cite{ern_guermond2004}), 
{\it i.e.}
\begin{equation*}
\franz{\V_{h}:=\left\{\vv_{h}\in\mathcal{C}^0(\overline{\Omega};\R^d):\vv_{h|K}\in\PP_k(K;\R^d),
\forall K\in\mathcal{K}_h,\vv_{h}=\0\mbox{ on }\Gamma_D\right\}.}
\end{equation*}
\franz{The finite element approximation of the small strain elasticity Problem \eqref{eq:primal_weak} is, as usual:}
\begin{equation}
\left\{\begin{array}{l}
\textrm{Find }\u_{h} \in \V_{h}\mbox{ such that }\\\noalign{\smallskip}
a(\u_{h}, \vv_{h}) = \ell(\vv_{h}), \qquad \forall \vv_{h} \in \V_{h}.
\end{array}\right.
\label{eq:primal_weak_disc}
\end{equation}

\subsection{\franz{A first taste of Dual Weighted Residuals in a fully linear setting}}
        \label{sec:goal_oriented_error_estimates}
\franz{
We provide a first insight into the
Dual Weighted Residual (DWR) technique
introduced in 
\cite{becker-rannacher-1996,becker-rannacher-2001}
in a very general (nonlinear) setting.
For the linear setting described above, we 
can simplify many calculations, to emphasize the main ideas, see
\cite{rognes-logg-2013}. 
}

\subsubsection{\franz{Estimator based on an exact dual problem}}

\franz{To $\u_h$, the solution to Problem \eqref{eq:primal_weak_disc},
we associate the
{weak residual} 
defined for all $\vv\in\V$ by}
%
\begin{equation*}\label{eq:def residual}
r(\vv) := \ell(\vv) - a(\u_h, \vv).
\end{equation*}
%
%
Moreover, we recall that we have the goal functional $J(\u)$ at hand. In the following, we are interested in minimizing the finite element discretization error measured in $J(\cdot)$ such that
\begin{equation*}
    J(\u) - J(\u_h)
\end{equation*}
is sufficiently small. This problem statement can be framed as a constrained optimization problem
on the continuous (non-discretized) level as follows:
\begin{equation*}
    \min J(\u) \quad\text{ s.t. } \quad a(\u,\vv) = \ell(\vv) \quad\forall \vv\in\V.
\end{equation*}
Employing the Lagrange formalism yields an unconstrained optimization problem with the Lagrangian $L:\V\times \V\to\mathbb{R}$ with
\begin{equation*}
    {L}(\u,\z) := J(\u) - a(\u,\z) + \ell(\z)
\end{equation*}
in which $\z\in\V$ is a Lagrange multiplier, i.e., the so-called adjoint solution. The solution of this 
saddle-point problem is obtained by its first-order necessary condition, by which we differentiate 
in both solution variables, namely $\u$ (primal solution) and $\z$ (adjoint solution). The optimality 
system reads
\begin{align*}
    J(\u)'(\delta \u) - a(\u,\z)'(\delta \u) &= 0 \quad\forall\delta \u\in \V,\\
     - a(\u,\z)'(\delta \z) + \ell(\z)'(\delta \z) &= 0 \quad\forall\delta \z\in \V.
\end{align*}
We notice that for linear PDEs and linear goal functionals the optimality system reads
\begin{align*}
    J(\delta\u) - a(\delta\u,\z) &= 0 \quad\forall\delta u\in \V,\\
     - a(\u,\delta\z) + \ell(\delta\z) &= 0 \quad\forall\delta \z\in \V.
\end{align*}
To summarize, let $\z$ denote the solution to the dual problem:
\begin{equation}
\left\{\begin{array}{l}
\mbox{Find }\z \in \V\mbox{ such that}\\\noalign{\smallskip}
a(\vv, \z) = J(\vv), \qquad \forall\, \vv \in \V.
\end{array}\right.
\label{eq:DefDual}
\end{equation}
%
%
\twick{
From 
\begin{equation*}
    J(\delta\u) - a(\delta\u,\z) = 0
\end{equation*}
by using $\delta\u = \u$, we obtain
\begin{equation*}
    J(\u) = a(\u,\z)
\end{equation*}
The same procedure can be done for the finite element solution, i.e.,
\begin{equation*}
    J(\u_h) = a(\u_h,\z_h).
\end{equation*}
Then, 
}
the DWR method, in a linear setting, relies on the fundamental observation that
\begin{align}
J(\u) - J(\u_h) = a(\u, \z) - a(\u_h, \z) = \ell(\z) - a(\u_h, \z) = r(\z).
\label{eq:J_resid}
\end{align}
From this, we design an error estimator of $J(\u) - J(\u_h)$ 
as an approximation of the residual $r(\z)$. We detail the different steps below.

\subsubsection{Numerical approximation of the dual problem and global estimator}\label{sec:dual}
The exact solution $\z$ to the dual system \eqref{eq:DefDual} is unknown in most practical situations, and thus needs to be approximated.
Let us consider a finite element space $\widehat{\V}_h \subset \V$. {This space is assumed to be} finer than $\V_h$. \franz{For instance, it can be made of continuous piecewise polynomials of higher order ($k+1$) or with a finer (nested) mesh.}
\revision{The approximation $\widehat{\z}_{h}$ of the solution $\z$ to the dual problem can be obtained:
\begin{enumerate}
    \item Either by solving the following discrete dual problem directly on the finer space $\widehat{\V}_h$:
\begin{equation}
\left\{\begin{array}{l}
\mbox{Find }\widehat{\mathbf{z}}_{h} \in \widehat{\V}_{h}\mbox{ such that}\\\noalign{\smallskip}
a(\widehat{\vv}_{h}, \widehat{\mathbf{z}}_{h}) = J(\widehat{\vv}_{h}), \qquad \forall\, \widehat{\vv}_{h} \in \widehat{\V}_h.
\end{array}\right.
\label{eq:DefDual discr1}
\end{equation}
\item Or by, first, solving an approximate dual problem on the original finite element space $\V_h$:
\begin{equation}
\left\{\begin{array}{l}
\mbox{Find }{\mathbf{z}}_{h} \in {\V}_{h}\mbox{ such that}\\\noalign{\smallskip}
a({\vv}_{h}, {\mathbf{z}}_{h}) = J({\vv}_{h}), \qquad \forall\, {\vv}_{h} \in {\V}_h.
\end{array}\right.
\label{eq:DefDual discr1 BIS}
\end{equation}
And then, by extrapolating the dual solution $\mathbf{z}_{h}$ on the finer space:
\[
\widehat{\mathbf{z}}_{h} = E_h {\mathbf{z}}_{h}
\]
where 
\[
E_h : {\V}_h \rightarrow \widehat{\V}_h
\]
is an extrapolation operator, 
described for instance in \cite{rognes-logg-2013}. This option can be a little less accurate howerver it is obviously much cheaper than solving a finite element problem on a finer space. 
\end{enumerate}}
\noindent We define
\begin{equation} \label{eq:globalestimatorh}
\eta_h := |r(\wzh)| 
\end{equation}
as the \emph{global estimator} that approximates the residual $|r(\z)|$. 

\subsubsection{Derivation of local estimators}\label{sec:local est}

Following \cite{becker-rannacher-2001,rognes-logg-2013}, we 
\franz{can split the global estimator into a sum of contributions on each mesh cell, which allows to refine cells that contribute the most to the error on $J$.
The \textit{local estimator} of the error $|J(\u)-J(\u_h)|$ is of the form}
\begin{equation}\label{estim j(u)-j(u_h)}
\sum\limits_{K\in\caK_h} \eta_K,
\quad
\eta_K := \left | 
\ds \int_K R_K \cdot (\wzh - i_h \wzh) \dd\x
\, + \sum_{E \in \caE_K} 
\int_E R_{E,K} \cdot (\wzh-i_h\wzh ) \mathrm{ds} \, 
\right |, \quad \forall K \in \caK_h,
\end{equation}
where the notation $i_h$ stands for the Lagrange interpolant onto $\V_h$. 
%
%
The local element-wise and edge-wise residuals are given explicitly by 
$$R_K:=\bB_K +{\bf div}\, \bs_A (\u_h)$$ and %
\begin{equation*} 
R_{E,K} :=\begin{cases}
-\dfrac{1}{2} 
\jumpEK{\bs_A(\u_h)}
& \text{if}~ E \in \caE^{int}_K, \\
\bB_E-\bs_A( \u_h ) \nEK  & \text{if}~ E \in \caE^N_K, 
\end{cases}
\end{equation*}
where $$\bs_A(\u^h) := \bs(\u^h) + \beta T ( \eb_A \otimes \eb_A ) {\chi_A}.$$
The notation $\chi_A$ stands for the indicator function of $\omega_A$, \emph{i.e.} $\chi_A=1$ in $\omega_A$ and $\chi_A=0$ elsewhere. The quantity $\bs_A(\u^h)$ represents the sum of passive and active contributions within the stress field.
The quantity $\bB_K$ (resp. $\bB_E$) is a computable approximation of $\bB$ (resp. $\bB_N$). 

\subsubsection{\franz{Some additionnal comments}}

\franz{Remark first that the following bound} always hold{s}
\[
\eta_h \leq \sum\limits_{K\in\caK_h} \eta_K,
\]
since compensation effects (balance between positive and negative local contributions) can occur for $\eta_h$, see, \emph{e.g.}, \cite{nochetto-veeser-verani-2009}. Thus $\eta_h$ is expected to be sharper than $\sum\limits_{K\in\caK_h} \eta_K$. In practice,$\sum\limits_{K\in\caK_h} \eta_K$ aims at quantifying the local errors for mesh refinement.

Each local estimator $\eta_K$ is made up of two contributions. On the one hand, the residuals $R_K$ and $R_{E,K}$ represent the local error in the natural norm. On the other hand, the contribution $(\wzh^i-i_h\wzh)$ coming from the dual problem can be interpreted as a \emph{weight} (or a sensitivity factor) that measures the local impact on the quantity of interest $J(\cdot)$, see, \emph{e.g.}, \cite[Remark 3.1]{becker-rannacher-2001}.
\label{rmq:weight}




\subsection{Algorithm for goal-oriented mesh refinement}\label{sec:algo}

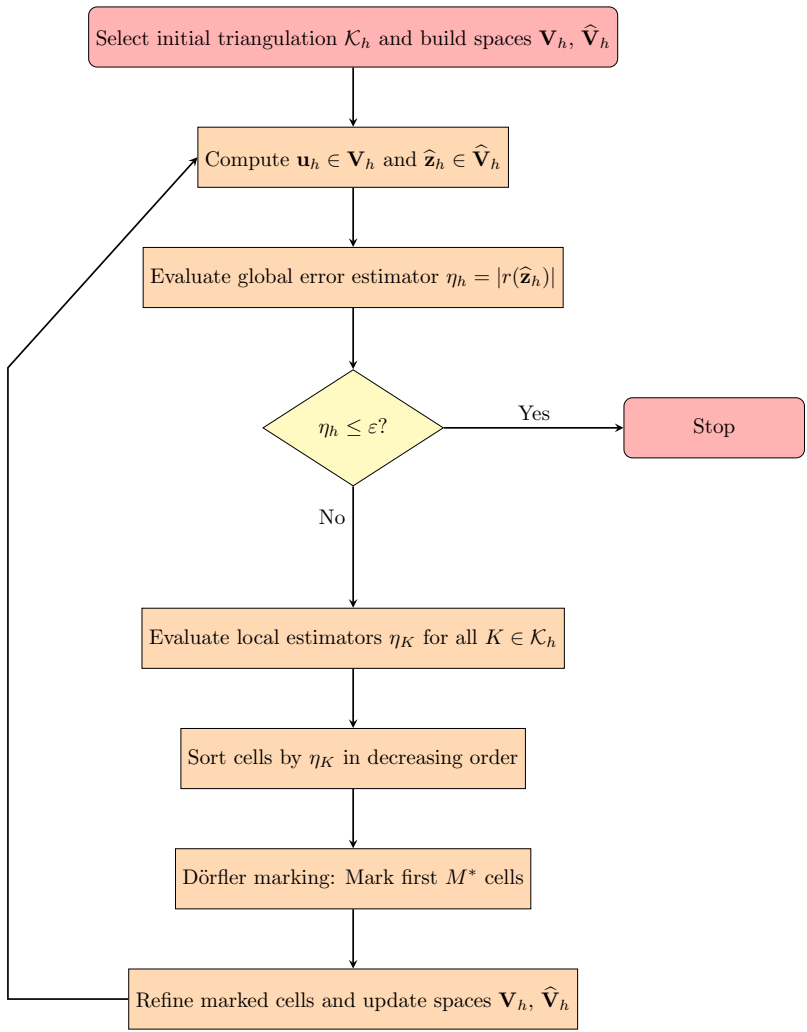
\begin{figure}[!ht]
    \centering
\begin{center}
\begin{minipage}{0.9\linewidth} 
\centering
\scalebox{0.8}{ 
\begin{tikzpicture}[node distance=2cm]

\node (start) [startstop] {Select initial triangulation $\mathcal{K}_h$ and build spaces $\V_h$, $\widehat{\V}_h$};
\node (compute) [process, below of=start] {Compute $\u_h \in \V_h$ and $\widehat{\z}_h \in \widehat{\V}_h$};
\node (error) [process, below of=compute] {Evaluate global error estimator $\eta_h = |r(\wzh)|$};
\node (decision) [decision, below of=error, yshift=-0.5cm] {$\eta_h \leq \varepsilon$?};
\node (stop) [startstop, right of=decision, xshift=4cm] {Stop};
\node (local) [process, below of=decision, yshift=-1.5cm] {Evaluate local estimators $\eta_K$ for all $K \in \mathcal{K}_h$};
\node (sort) [process, below of=local] {Sort cells by $\eta_K$ in decreasing order};
\node (mark) [process, below of=sort] {Dörfler marking: Mark first $M^*$ cells};
\node (refine) [process, below of=mark] {Refine marked cells and update spaces $\V_h$, $\widehat{\V}_h$};

\draw [arrow] (start) -- (compute);
\draw [arrow] (compute) -- (error);
\draw [arrow] (error) -- (decision);
\draw [arrow] (decision.east) -- ++(1.5,0) node[anchor=south] {Yes} -- (stop.west);
\draw [arrow] (decision.south) -- ++(0,-0.5) node[anchor=east] {No} -- (local.north);
\draw [arrow] (local) -- (sort);
\draw [arrow] (sort) -- (mark);
\draw [arrow] (mark) -- (refine);
\draw [arrow] (refine.west) -- ++(-2,0) -- ++(0,10.5) -- (compute.west);

\end{tikzpicture}
} 
\end{minipage}
\end{center}
    \caption{Flowchart for adaptive mesh refinement.}
    \label{fig:flowchart}
\end{figure}

\franz{Using the D\"orfler marking strategy \cite{dorfler1996}, we describe, in Algorithm \ref{algo:refinement}, a simple algorithm to refine the mesh by taking into account these quantities. 
The corresponding flowchart is given as well Figure~\ref{fig:flowchart}.
In this algorithm, there are two independent numerical parameters: first a parameter $0 < \alpha \leq 1$ 
that controls the level of refinement in D\"orfler marking, and then a tolerance threshold $\varepsilon > 0$ for the global estimator, that serves as a stopping criterion. A complete Python/FEniCS script is available, see \cite{Duprez2019}.}

\begin{algorithm}
\caption{ \franz{Refinement Algorithm}}
\begin{algorithmic}[1]
\State Select an initial triangulation $\mathcal{K}_h$ of the domain $\Omega$.
\State Build the finite element spaces $\V_{h}$ and $\widehat{\V}_h$. 
\While{$\eta_h > \epsilon$}
    \State Compute $\u_{h} \in \V_{h}$ such that:
    \State \quad $a(\u_{h}, \vv_{h}) = \ell(\vv_h), \quad \forall \vv_{h} \in \V_{h}$.
    \State Compute $\widehat{\z}_{h} \in \widehat{\V}_{h}$. 
    \State Evaluate the global error estimator $\eta_h = |r(\wzh)|$.
    \If{$\eta_h \leq \varepsilon$}
        \State \textbf{stop}.
    \EndIf
    \State Evaluate the local estimators:
    \begin{equation*}
        \eta_K := \left| 
        \int_K R_K \cdot (\wzh - i_h \wzh) \, \mathrm{d}\x
        + \sum_{E \in \caE_K} 
        \int_E R_{E,K} \cdot (\wzh^i - i_h \wzh) \, \mathrm{d}s\, 
        \right|, \quad \forall K \in \caK_h.
    \end{equation*}
    \State Sort the cells $\{K_1, \dots, K_N\}$ by decreasing order of $\eta_K$.
    \State \textbf{Dörfler marking}:
    \State Mark the first $M^*$ cells for refinement, where:
    \begin{equation*}
        M^* := \min \left\{ M \in \mathbb{N} \:\middle|\: 
        \sum_{i=1}^{M} \eta_{K_i} \geq \alpha \sum_{K \in \mathcal{K}_h} \eta_K 
        \right\}.
    \end{equation*}
    \State Refine all marked cells (and propagate refinement to avoid hanging nodes).
    \State Update the finite element spaces $\V_{h}$ and $\widehat{\V}_h$ accordingly.
\EndWhile
\end{algorithmic}
\label{algo:refinement}
\end{algorithm}

\clearpage

\section{\franz{Hyperelastic soft tissue}}
\label{sec:Methods}

\franz{This section now details how the previous approach can be extended to a fully nonlinear setting, with large strain elasticity and a possibly nonlinear goal functional.}


\subsection{\franz{Setting}} 


\franz{In this section we use the same notations as previously in Section~\ref{sec:method1}, but with the following changes.
The domain $\Omega$ is occupied by 
a 
hyperelastic body in its reference configuration 
and subjected to a given body force $\bB$. The unknown displacement field 
is still denoted by $\u$. 
The deformation gradient is denoted $\bF$, with $\bF := \bI + \nabla_X \u$, where $\bI$ stands for the identity matrix, and $\nabla_X$ denotes the gradient with respect to coordinates in the reference configuration. 
The first Piola-Kirchhoff stress tensor is denoted $\bPi$ and the second Piola-Kirchhoff stress tensor is denoted $\bS$.} 
\franz{Contractile structures like muscles in soft tissues can carry out an internal strain/stress for the tissue itself, see \cite{Costa2001}. The behaviour of contractile soft tissues can be modelled by using two different strategies: pre-stress or pre-strain activations, see, \eg{} \cite{OHAYON2017}. 
In this paper, we adopt the pre-stress approach in order to model the active behaviour of soft biological tissues. This method consists in adding an active stress tensor to the classical (passive) stress tensor \cite{Bovendeerd1992,Smith2004}. The total stress reads
\begin{equation}
\bPi = \bF \cdot \bS + \bPi_a,
\end{equation}
where $\bF \cdot \bS$ is the passive contribution
and where
$\bPi_a$ is an active stress tensor. The expression of this latter is chosen as
\begin{equation}
\bPi_a = \beta T_a (\bF f_0) \otimes f_0,
\end{equation}
in which $0 \leq \beta \leq 1$ is the activation parameter, $T_a$ stands for the active tension of the muscle fibre, and $f_0$ is a unitary vector which describes the fibre orientation. The coefficient $\beta$ stands for the rate of activation: when $\beta=0$ there is no activation (purely passive behavior) whilst $\beta=1$ corresponds to the maximum activation.
In this paper we employ a constant active tension $T_a$ for the fibre. However, the active tension can be modelled as a function of fibre stretch ratio, see, \eg{}, \cite{Mulvany1979,Rachev1997}.
}
\franz{The virtual works associated with the internal forces is: 
\begin{equation*} 
A(\u; \vv) := \int_{\Omega} \bPi(\u): \nabla_X \vv ~\dd\x,
\end{equation*}
where $\u$ and $\vv$ are admissible displacements. 
The virtual work relative to external forces reads:
\begin{equation*} 
L(\vv) := \int_{\Omega} \bB \cdot \vv \, 
~\dd\x
+
\int_{\Gamma_N} \bB_\Gamma \cdot \vv  ~\dd\mathrm{s},
\end{equation*}
and, finally, 
the hyperelastic problem in weak form reads
\begin{equation}
\left.
\begin{array}{l}
\text{Find a displacement }\bu 
\text{ such that } \\\noalign{\smallskip}
A(\u; \vv) = L(\vv), \, \forall \vv. 
\end{array}
\right.
\label{eq:primal_weak2}
\end{equation}}
%

\begin{remark}
\franz{The model presented in the previous section can be viewed as a linearization of this hyperelastic model. Indeed, if we decompose the (large) strain Green-Lagrange deformation tensor 
$\bE (:= \frac{1}{2}(\bF^T \bF - \bI ))$
into the sum of a linear term and a nonlinear term as
\[
\bE(\u) = \boldsymbol{\epsilon}(\u) + \boldsymbol{\rho}(\u), \quad
\boldsymbol{\epsilon}(\u) = \frac12 \left ( \nabla_X \u + (\nabla_X \u)^T \right ), \quad
\boldsymbol{\rho}(\u) = \frac12 (\nabla_X \u)^T \nabla_X \u,
\]
and if we suppose the constitutive law of the form
\[
\bS = \bS_s (\boldsymbol{\epsilon}) + \bS_l ( \boldsymbol{\rho} ),
\]
then 
\[
\bPi (\u) = \bF(\u) \bS(\u) 
= (\bI + \nabla_X \u)  
(\bS_s (\boldsymbol{\epsilon}(\u)) + \bS_l ( \boldsymbol{\rho}(\u) )) 
= \mathcal{C} : \boldsymbol{\epsilon}(\u) + \bPi_l (\boldsymbol{\rho}(\u))
\]
where $\mathcal{C} : \boldsymbol{\epsilon}(\u)$ is a small strain (linear) elasticity term, see Section~\ref{sec:toy problem}
($\mathcal{C}$ is a fourth-order elasticity tensor), and $\bPi_l$ contains the nonlinear terms, of high order in $\nabla_X \u$.}
\end{remark}

\franz{Let us still denote $\mathcal{K}_h$  a mesh of the domain $\Omega$.
Let us denote by $\bV_h $ a conforming finite element space made for instance of Lagrange finite elements. The finite element method to solve our hyperelastic problem reads 
%
\begin{equation}
\left\{
\begin{array}{l}
\text{Find a displacement }\bu_h \in \bV_h 
\text{ such that } \\\noalign{\smallskip}
A(\u_h ; \vv_h ) = L(\vv_h), \, \forall \vv_h \in \bV_h. 
\end{array}
\right.
\label{eq:FEM_primal_weak}
\end{equation}
}

\subsection{Dual problem for computing the weights}

One of the main ingredients of the DWR method is to solve an adjoint problem to extract information from the quantity of interest and get feedback about the regions where it is more, or less, influenced, by the approximation error. 
\twick{The motivation is the same as in Section~\ref{sec:method1}, i.e., 
the (discretization) error \franz{for a nonlinear quantity of interest $J$} 
\begin{equation*}
   |J(\bu)-J(\bu_h)|
\end{equation*}
shall be minimized with the PDE as a constraint:
\begin{equation*}
    \min |J(\bu)-J(\bu_h)|
\quad\text{s.t. } A(\cdot,\cdot; \cdot,\cdot) = L(\cdot)
\end{equation*}
Thus, 
to satisfy this constrained minimization problem, Lagrange multipliers are introduced in order to formulate 
a Lagrangian functional.}
\revision{This 
implies computing the following dual problem:
\begin{equation}
\left\{
\begin{array}{l}
\text{Find } {\bz} \in {\bV} \text{ such that } \\\noalign{\smallskip}
(A')^*(\bu | {\bz} ;{\bv}) = J'(\bu|{\bv}) \quad \forall {\bv} \in {\bV},
\end{array}
\right.
\label{eq:nonlinear_dual_problem}
\end{equation}
where $A'$ and $J'$ denote the Fr\'echet derivative of $A$ and $J$, respectively,
and $(A')^*$ is the adjoint form of $A'$.} 

\subsection{\revision{Discrete dual solution}}

Let us still denote $\widehat{\V}_h \subset \V$
a finite element space finer than $\V_h$, {for instance made of continuous piecewise polynomials of higher order ($k+1$).}
\revision{As for small strain elasticity, the approximation $\widehat{\z}_{h}$ of the solution $\z$ to the dual problem can be obtained:
\begin{enumerate}
    \item {Either by solving the following discrete dual problem directly on the finer space $\widehat{\V}_h$:
\begin{equation}
\left\{
\begin{array}{l}
\text{Find } \widehat{\bz}_h \in \widehat{\V}_h \text{ such that } \\\noalign{\smallskip}
(A')^*(\bu_h| \widehat{\bz}_h ;\widehat{\bv}_h) = J'(\bu_h|\widehat{\bv}_h) \quad \forall \widehat{\bv}_h \in \widehat{\V}_h.
\end{array}
\right.
\label{eq:nonlindiscrete_dual_problem}
\end{equation}}
\item {Or by, first, solving an approximate dual problem on the original finite element space $\V_h$:
\begin{equation}
\left\{
\begin{array}{l}
\text{Find } {\bz}_h \in {\bV}_h \text{ such that } \\\noalign{\smallskip}
(A')^*(\bu_h| {\bz}_h ;{\bv}_h) = J'(\bu_h|{\bv}_h) \quad \forall {\bv}_h \in {\bV}_h,
\end{array}
\right.
\label{eq:discrete_dual_problemnonlinbis}
\end{equation}}
And then, by extrapolating the dual solution $\mathbf{z}_{h}$ on the finer space:
\[
\widehat{\mathbf{z}}_{h} = E_h {\mathbf{z}}_{h}
\]
where 
\[
E_h : {\V}_h \rightarrow \widehat{\V}_h
\]
is an extrapolation operator, 
described for instance in \cite{rognes-logg-2013}. 
\end{enumerate}}

%

Remark that the dual problem \eqref{eq:nonlinear_dual_problem} is linear, so solving it is not expensive in comparison to \eqref{eq:primal_weak2}.
{Moreover, the bilinear form \eqref{eq:primal_weak2} (left-hand side) has been already assembled in the last Newton iteration in the resolution of the primal problem. }
\twick{A mathematical difficulty results from Galerkin orthogonality as both the continuous-level and finite element dual solutions enter into the a posteriori error estimator. For a full practical evaluation, the continuous-level solution must be solved as well with the help of the finite element method.
If both come from the same discrete function spaces, the dual sensitivity weight will vanish identically. Therefore, a high-order solution must be constructed, either by solving a higher-order finite element problem or by a local higher order reconstruction.}

For model problems or more complex problems such as \eqref{eq:primal_weak2}, and some intricate expressions of $J$, the practical calculation of $A'$ and $J'$ can be fastidious. For this purpose, we take advantage of the capabilities of automatic symbolic differentiation embedded into modern finite element software such as FEniCS or GetFEM++. 
Furthermore, this feature makes possible some genericity in the implementation: virtually nothing has to be changed in the program if the hyperelastic constitutive law is modified.

\subsection{The representation formula of Becker and Rannacher}

\franz{We introduce $r(\bu_h;\bv)$ the residual of Problem \eqref{eq:FEM_primal_weak} as
\begin{equation}
\label{eq_primal_residual}
r(\bu_h;\bv) := L(\bv) - A(\bu_h;\bv) \qquad \forall \bv \in \bV.
\end{equation}
This, roughly speaking, quantifies how well the hyperelasticity equations are approximated (it should tend to zero if the mesh is uniformly refined).
Thanks to the dual system \eqref{eq:nonlinear_dual_problem}, we obtain expression of the error on $J$ 
as the best approximation term involving the residual and the (exact) dual solution
(see \cite[Proposition 2.3]{becker-rannacher-2001}): 
\begin{equation}
J(\bu) - J(\bu_h) = \min_{\bv_h \in \bV_h} r(\bu_h; \bz - \bv_h) + R_m \label{eq:representation_DWR}
\end{equation}
where 
$R_m$ is the high-order remainder related to the error caused by the linearization of the nonlinear problem (the precise expression of which can be found in \cite{becker-rannacher-2001}). 
In practice, this quantity is, hopefully, negligible.
Note at this stage that there are various possibilities to represent the error on $J$, which are detailed in \cite{becker-rannacher-2001}, and for instance, in \cite{rognes-logg-2013}, the authors make use of another representation formula (Proposition 2.4) which is then approximated. 
%
Proceeding as usual in \emph{a posteriori} error estimation, \emph{i.e.,}
after performing integration by parts on the residual $r$, we localize the different contributions to the goal-oriented error as follows:
\begin{equation}
\lvert J(\bu) - J(\bu_h) \rvert \leq \sum_{K \in \mathcal{K}_h} \eta_K{(\bu_K,\bz_K)} + H.O.T.
\end{equation}
In the above expression, $K$ denotes any cell of the mesh $\mathcal{K}_h$, and expressions such as $\bu_K$ denote the local restriction of the finite element variable $\bu_h$ on the cell $K$. Moreover, $H.O.T.$ denotes high order terms, that are not considered in the implementation. In section 2.6 below, the detailed expression of $\eta_K$ is given.} 
\twick{Alternative approaches for the error localization are 
the filtering approach~\cite{braack2003posteriori} or a partition-of-unity~\cite{RiWi15_dwr}.
The latter was extended to space-time versions in \cite{ThiWi24,ENDTMAYER2024286}; see also the overview article \cite{wick2024}.
}

\subsection{Expression of the estimator} 

\franz{We provide below for each cell-wise contribution:\revision{
\begin{equation}
\eta_K = \left| \int_K \bR_{\bu} \cdot (\widehat\bz_h-i_h(\widehat\bz_h)) ~\dd\x
+  \int_{\partial K} \bJ \cdot (\widehat\bz_h-i_h(\widehat\bz_h)) ~\dd\mathrm{s}\, 
\right|
\end{equation}}
with, $E_h$ and $I_h$ are resp. the extrapolation (see \cite{rognes-logg-2013}) and the Lagrange interpolation, the interior residual
\begin{equation*}
\bR_{\bu} = \bB + 
\mathrm{div} \, \bPi(\bu_h)
\end{equation*}
and the stress jump 
{
\begin{equation*}
\bJ = \begin{cases}
   - \frac{1}{2}   [[\bPi(\bu_h)]]
       & \text{ if } F \not \subset \Gamma, \\
      \bT -  \bPi(\u_h)  \cdot n_F                                 & \text{ if } F \subset \Gamma_N, \\
      \bzero                                                                     & \text{ if } F \subset \Gamma_D,
      \end{cases}
\end{equation*}
for each facet $F$, where $n_F$ is the exterior normal to the facet $F$ of $\Gamma_N$. }
The jump can be defined for a function $\bv_h$ on a facet $F$ between two cells $K$ and $K'$ by 
$[[\bv_h]]=\bv_{h|K}\cdot n_K + \bv_{h|K'}\cdot n_{K'}$, where $n_K$ and $n_{K'}$ are the normal of $K$ and $K'$ on $F$.}

\subsection{Adaptive mesh refinement}

Using the error estimate on $J$, we implement a standard procedure for mesh refinement. 
As described in \cref{algo:adaptivity}, we start with an initial mesh called $mesh_i$, and by providing a guessed solution $\bu_i^{(0)}$, the nonlinear primal problem can be solved using Newton's method. 
Once accepting $\bu_i$ as the solution of the discrete primal problem, solving the discrete dual problem 
provides the dual solution $\bz_i$ on the current mesh $mesh_i$. The estimator $\eta_K$ is then computed by using the primal and dual solutions $\bu_i$ and $\bz_i$, respectively. From the estimator, different strategies can be used to mark the elements whose \emph{error} is high. In this paper, we use the D\"orfler marking strategy \cite{dorfler1996}. 

\subsection{\twick{Multigoal-oriented error control and adaptivity}}
\label{sec:multigoal_section}
{\color{black} 
For complex biomechanics situations in which several physical phenomena interact, multigoal-oriented error estimates are of interest. Therein, several quantities of interest,
say $S$, such as $J_1,\ldots, J_S$,
might be considered simultaneously. 
They allow to focus simultaneously into the different parts of the geometry and/or different \franz{solution} 
components and/or different types of quantities of interest such as point evaluations, line integrals 
or domain integrals. 
For such cases, one possibility is to construct a combined functional of the form
\begin{equation}
\label{eq_multiple_goals}
  J_c(\bu):= \sum_{i=1}^S w_i J_i(\bu),
  \end{equation}
with some \franz{positive or negative} weights $w_i\in\mathbb{R}$ for $i=1,\ldots,S$. Single goal functionals
are modeled as before with $S=1$ and $w_1 = 1$. In the case $S\geq 2$, 
\franz{it can \revision{easily happen}  situations where $J_c(\bu) \simeq 0$ due to cancellations of different $J_i$ with same \revision{absolute} values, but opposite signs.}
\revision{This can be seen as follows. As in the previous sections, at the end, we are interested in the discretization error $J(\bu) - J(\bu_h)$. For simplification, let us assume the linear case, where $J(\bu) - J(\bu_h) = J(\bu - \bu_h)$. Then, let us employ~\eqref{eq_multiple_goals} for some arbitrary element $\bv$, yielding
\begin{equation}
\label{eq_Jc_bv}
  J_c(\bv):= \sum_{i=1}^S w_i J_i(\bv).
\end{equation}
In order to obtain a meaningful error representation for the choice $\bv := \bu - \bu_h$, we must know the signs of $J_i(\bu) - J_i(\bu_h)$, thus we choose
\begin{equation}
\label{eq_multi_goal_signs_original}
  w_i := \omega_i \frac{\text{sign}(J_i(\bu) - J_i(\bu_h))}{|J_i(\bu_h)|}, \quad \omega_i \geq 0.
\end{equation}
Employing this expression and \franz{$\bv := \bu - \bu_h$} into~\eqref{eq_Jc_bv}, we obtain
\begin{equation}
    J_c(\bu - \bu_h):= \sum_{i=1}^S \omega_i \frac{\text{sign}(J_i(\bu) - J_i(\bu_h))}{|J_i(\bu_h)|} J_i(\bu - \bu_h).
\end{equation}
Therefore, to have a meaningful overall error in the combined goal functional $J_c$, the single errors $J_i(\bu - \bu_h)$ should not cancel out due to their signs.
The task is to get the information about the signs
\begin{equation}
    \text{sign}(J_i(\bu) - J_i(\bu_h)).
\end{equation}
We 
\franz{mention} that a \revision{first} sign computation
using an adjoint-adjoint problem
was introduced in~\cite{Ha08,HaHou03}. \revision{A brief explanation of the adjoint-adjoint is as follows (see~\cite{HaHou03}[Section 2]). Solving $S$ adjoint problems, namely  for each goal functional $J_i$, becomes costly \franz{if} $S\gg 1$. As an alternative, an error equation can be introduced as: Find the error $e$ in some admissible function space such that $A(e,w) = r(\bu_h,w)$ for all admissible test functions $w$. Here, $r(\bu_h,\cdot)$ is the primal residual introduced in~\eqref{eq_primal_residual}. With the error $e$ at hand, the $S$ goal functionals $J_i(e)$ can be evaluated directly. The mathematical operation to deduce $A(e,w) = r(\bu_h,w)$ is to take the adjoint of the adjoint equation (see again~\cite{HaHou03}[Section 2]).
}
}

\revision{In order to avoid the adjoint-adjoint problem, which means to solve a third PDE besides the primal PDE (our original problem statement) and the adjoint PDE (for the dual weights), an alternative way was proposed in~\cite{wick2017pum} by 
approximating~\eqref{eq_multi_goal_signs_original} with}
\begin{equation}
\label{eq_multi_goal_signs}
  w_i := \omega_i \frac{\text{sign}(J_i(\bu_h^{(2)}) - J_i(\bu_h))}{|J_i(\bu_h)|}, \quad \omega_i \geq 0,
  \end{equation}
where $\bu_h^{(2)}$ is some higher order approximation~\cite{endtmayer2021,wick2017pum} (see also~\cite{wick2024}). 
\revision{We mention that this sign approximation is also mathematically justified; see~\cite{wick2017pum}[Prop. 3.1].}

Using the above definition of $J_c$ in~\eqref{eq_multiple_goals}, under the condition that the signs are computed as shown in~\eqref{eq_multi_goal_signs}, the further procedure is the same as discussed in Section~\ref{sec_dual_problem}, namely we have for the dual problem
\begin{equation}
\left\{
\begin{array}{l}
\text{Find } ({\bz}_h,\bw_h) \in {\bV}_h^0  \times \bJ_h \text{ such that } \\\noalign{\smallskip}
(A')^*(\bu_h,\bpp_h| {\bz}_h,  {\bw}_h;{\bv}_h, {\bq}_h) = J_c'(\bu_h,\bpp_h|{\bv}_h,\bq_h) \quad \forall ({\bv}_h,\bq_h) \in {\bV}_h^0  \times \bJ_h.
\end{array}
\right.
\label{eq:discrete_dual_problemIH_multi_goal}
\end{equation}
We notice that on the right hand side $J_c$ is employed, \franz{that \revision{takes} into account all the different functionals}.

The adaptive mesh refinement can be carried out with the same previous algorithms. The only additional step is the sign computation~\eqref{eq_multi_goal_signs} for $w_i$ after the primal problem has been computed. Then, the $J_c$ is designed, 
which is followed by the computation of the dual problem.
As mentioned before, in~\cite{wick2017pum}[Prop. 3.1], the sign computation was simplified in comparison to~\cite{Ha08,HaHou03}. \revision{As explained in}~\cite{wick2017pum} a higher order solution \revision{$\bu_h^{(2)}$} of the primal problem is required. \revision{As before in the previous sections, the same procedure can be applied: either} by a higher order finite element solution or by subsequent higher order local interpolation of a low-order solution. Several comparisons about the consequences of either way were conducted in~\cite{endtmayer2021}. Then, with $J_c$, the dual problem~\eqref{eq:discrete_dual_problemIH_multi_goal}
can be solved. The rest of the adaptive algorithm and mesh refinement is the same as in the single goal case. 
\revision{Therein, it should be mentioned that for nonlinear goal functionals that calculation of the derivative $J_c'(\cdot)(\cdot)$ can be cumbersome in practice. However, the mathematical procedure is the same as usual by computing the directional derivative of $J_c(\bu)$
in the linearization point $\bu$ into the direction $\delta\bu$ in some admissible function space, i.e., we obtain $J_c'(\bu)(\delta\bu)$. For details, we refer to~\cite{EnLaWi18}[Section 4].}
}

\section{Fluid-structure interaction}
\label{sect:fluidstructure}

\twick{In the previous two sections, we have focussed on solid mechanics only. In the following, we extend those to the interaction with fluid flows. In fact,}
fluid-structure interaction (FSI) plays a central role in many biomechanical systems, where biological tissues deform in response to fluid flow, and vice versa. Prominent examples include blood flow in arteries~\cite{Formaggia2009,Failer2021,BoGaNe14,RazzaqDamanikHronOuazziTurek:2012,Croelal10}, airflow in the lungs~\cite{heil2011,Wall2008} \franz{or in the upper airways~\cite{fc2006,heil2011}}, cerebrospinal fluid dynamics~\cite{Cardillo2021}, and 
the motion of heart valves~\cite{VANLOON2006806,Astorino2009,Wick2012,Wi14_fsi_eale_heart,Hsu2014,Hoffman2021,Fuma2023,FENG2024116724}. In these systems, the mechanical interaction between a fluid (e.g., blood, air) and compliant biological structures (e.g., vessel walls, heart tissue) critically affects both function and pathology.

In cardiovascular biomechanics, for instance, the pulsatile nature of blood flow induces complex pressure and shear forces on arterial walls, which in turn deform and influence flow patterns. This two-way coupling is essential to understanding atherogenesis, aneurysm development, and vascular remodeling. Neglecting structural compliance can lead to significant errors in predicting flow fields and stresses within tissues~\cite{Formaggia2009,YangJaegerNeussRaduRichter2015,YangRichterJaegerNeussRadu2015,Failer2021,BoGaNe14,RazzaqDamanikHronOuazziTurek:2012}.

Similarly, in respiratory biomechanics, the interaction between airflow and lung tissue is vital for accurately modeling lung mechanics and disease progression in conditions like emphysema or pulmonary fibrosis~\cite{Faizal2019}. 
\franz{In the same way, interaction between the respiratory flow and the surrounding soft tissue in the upper airway is fundamental for pathologies such as obstructive sleep apnea~\cite{ashraf2022fluid} and snoring~\cite{Faizal2019}. It is also strongly involved in speech production processes~\cite{heil2011,mittal2013}, notably at the level of the vocal folds~\cite{sundstrom2025}.}
FSI is also critical in the design of prosthetic heart valves, where device performance is governed by dynamic interactions with blood flow~\cite{Arminio2024}.
\franz{Let us finally mention biomechanics of the digestive system where fluid-structure interactions occur during the whole trajectory
of the food bolus~\cite{engmann2013}.}

From a computational standpoint, the accurate resolution of such coupled phenomena is challenging due to the multiphysics coupling between incompressible fluids and hyperelastic solids, the large structural deformations in soft tissues and the moving boundaries and interfaces, which demand flexible and robust numerical schemes.

A critical numerical challenge in fluid-structure interaction problems - especially in biomechanics - is the so-called added mass effect. This phenomenon arises when a light or compliant structure is immersed in a denser fluid, leading to a strong inertial coupling between the two domains~\cite{Causin2005}. The effect is particularly pronounced in applications like blood flow in arteries or cerebrospinal fluid–brain tissue interaction, where the density of the fluid is comparable to or exceeds that of the surrounding tissue. In such cases, even small structural displacements can induce significant changes in the surrounding fluid field, which then exerts substantial feedback forces on the structure. From a computational standpoint, the added mass effect can lead to instabilities and loss of accuracy in partitioned FSI solvers, especially when the coupling is \franz{treated explicitly~\cite{fernandez2011}}. 

Mathematical modeling and simulation of fluid-structure interactions is well developed~\cite{Richter2017} \twick{(as well as the research monographs \cite{BaTaTe13,FrHoRiWiYa17,Formaggia2009,BoGaNe14})}, the problem however remains a highly challenging one. For biomechanical problems often strongly coupled or monolithic formulations must be used~\cite{Heil2008}. Efficient solvers are rare~\cite{Richter2015a,Failer2020,Boffi2024,JoLaWi19_fsi,BaQuaQua08,BaQuaQua08b,gee2011truly,CrDeFouQua11} and the computational cost is immense. Hence, adaptivity as a measure to reduce the complexity is of great importance. 

There are several contributions to goal oriented adaptivity for fluid-structure interactions starting with the early works of Dunne~\cite{Dunne2006} considering a Eulerian formulation,  van der Zee~\cite{Zee2011} treating the stationary Stokes problem coupled to a lower-dimensional solid and Richter and Wick~\cite{Richter2012,Wick2012a} considering the stationary nonlinear problem.  Failer applied the Dual Weighted Residual method to the fully coupled nonstationary and nonlinear fluid-structure interactions problem~\cite{Failerdiss,FaiWi18}. 
\twick{
Multigoal error control applied to stationary fluid-structure interaction and academic benchmark examples was done in~\cite{AhEndtSteiWi22}.}

In the following, we start by introducing the most basic fluid-structure interaction problem coupling the incompressible Navier-Stokes equation in a $d$-dimensional domain to an hyperelastic solid of same dimensionality. To focus on the essentials, we present a fully monolithic formulation and derive its variational form as basis for the upcoming discussion. Next, we briefly describe the basics discretizing the problem in space and time. \franz{We then give a detailed derivation of the adjoint problem with a special consideration of the coupling conditions between fluid and solid.} 

\subsection{Governing equations}

Taking cardiovascular mechanics as the prototypical typical biomechanical application of FSI, models of great complexity are studied to describe blood~\cite{Bessonov2015} and tissue~\cite{Holzapfel2000}. In this \twick{book chapter}, we will therefore limit ourselves to a simple model, the incompressible Navier-Stokes equations for the fluids velocity $\vv$ and pressure $p$
\begin{equation}\label{fsi:fluid}
\div\,\vv=0,\quad
\rho_f\big(\partial_t \vv + (\vv\cdot\nabla)\vv\big) 
-\div\bs_f = \rho_f\f\text{ in }\FL
\end{equation}
and the St. Venant Kirchoff model for the structure's deformation $\u$ and velocity $\vv$
\begin{equation}\label{fsi:solid}
d_t \u = \vv,\quad
\rho_s d_t \vv - \div\big( \F\bSI_s \big) = \rho_s \f\text{ in }\SO.
\end{equation}
Here, $\rho_f$ and $\rho_s$ are the densities of fluid and solid, $\F = \nabla\u$ the deformation gradient, 
\begin{equation}\label{fsi:sigmaf}
\bs_f := \rho_f\nu_f (\nabla\vv+\nabla\vv^T)-pI
\end{equation}
is the Cauchy stress of the fluid with the pressure $p$ and the kinematic viscosity $\nu_f$ and
\begin{equation}\label{fsi:sigmas}
\bSI_s := 2\mu_s \E +\lambda_s \tr(\E) I
\end{equation}
is the St. Venant Kirchhoff material law, where $\mu_s$ and $\lambda_s$ are the Lam\'e parameters and 
\begin{equation}\label{fsi:strain}
\E = \frac{1}{2}\big(\F^T\F-I)
\end{equation}
is the Green-Lagrange strain tensor. 

\subsubsection{Coordinate systems}

The two equations for fluid~\eqref{fsi:fluid} and solid~\eqref{fsi:solid} are formulated in different coordinate systems. While the Navier-Stokes equations are noted in the classical Eulerian system, where the coordinate $\x\in\mathbb{R}^d$ denotes a fixed point in space, the solid problem is given in the Lagrangian system and $\hat \x\in\mathbb{R}^d$ here denotes a \emph{material point} with a position $x(\hat \x,t)=\x+\u(\hat \x,t)\in\mathbb{R}^d$ that moves in time. This means that the solid problem is formulated on a fixed Lagrangian domain $\SO\subset\mathbb{R}^d$ while the fluid problem is given on the moving domain $\FL(t)\subset\mathbb{R}^d$, where the domain motion originates from the solution of the coupled system itself. Bringing these two frameworks together is the main technical difficulty of fluid-structure interactions. 

Basically, two different approaches exist: first, we transform the fluid problem onto a coordinate system that matches the solid's one. This coordinate system is neither Eulerian nor Lagrangian but an artificial one. This approach is called \emph{Arbitrary Lagrangian Eulerian (ALE)} and it is the standard approach for most cases and also detailed in the following. The alternative way for coupling is to map the solid problem into Eulerian coordinates onto a moving framework $\SO(t)$ that is naturally given by its deformation as
\[
\SO(t) = \{\x = \hat \x + \u(\hat \x,t),\quad \hat \x\in\SO\}.
\]
This \emph{Fully Eulerian} approach goes back to pioneering works by Dunne~\cite{Dunne2006} and Cottet, Maitre and Milcent~\cite{Cottet2008}.  It has some theoretical advantages, all based in the fact that both the Eulerian and the Lagrangian system are physically sound, whereas the ALE formulation has no physical foundation. However, technically, the Fully Eulerian approach brings difficulties when it comes to discretization and solution, such that it is only considered in special cases, such as contact problems~\cite{Wahl2021}. We refer to~\cite{Richter2013,Wick2013,Richter2017,FrKnStWeWi25} for details. The DWR method has already been tested for the Fully Eulerian approach in early works of Dunne~\cite{Dunne2006}.

To construct the ALE formulation we must introduce a reference domain $\FL$, usually the domain at initial time $\FL=\FL(0)$, which matches to the solid. Then, we need a map 
\begin{equation}\label{fsi:alemap}
  T(t):\FL\mapsto \FL(t),
\end{equation}
that maps the reference domain onto the moving Eulerian domain $\FL(t)$ in such a way, that $\FL(t)$ always matches the Eulerian solid counterpart $\SO(t)$. This is realized by mimicking the physical Lagrangian-Eulerian map of the solid
\begin{equation}\label{fsi:mapsolid}
T(t):\SO\mapsto \SO(t),\quad T(\hat \x,t):= \hat \x + \u(\hat \x,t)
\end{equation}
namely, by defining the ALE map as
\begin{equation}\label{fsi:alemap:1}
T(t):\FL\mapsto \FL(t),\quad T(\hat \x,t):= \hat \x +\u(\hat \x,t),
\end{equation}
where $\u|_{\FL}$ is an extension of $\u|_{\SO}$ from the common interface $\IN = \partial\FL\cap\partial\SO$ to $\FL$. The most simple strategy is to harmonically extend the solid deformation
\begin{equation}\label{fsi:alemap:2}
-\Delta \u = 0\text{ in }\FL,\quad 
\u|_{\FL} = \u|_{\SO}\text{ on }\IN,\quad 
\u=0\text{ on } 	\partial\FL\setminus\IN. 
\end{equation}
While simple and efficient to realize, this construction fail for large deformations. Variants are discussed in~\cite{StTeBe03,Wick2011} or \cite[Sec. 5.3.5]{Richter2017}. 

Using some fundamental relations for transforming differential operators from $\FL(t)$ to $\FL$, see~\cite[Sec. 2.1.7]{Richter2017}, the Navier-Stokes equations in ALE reference system is given as
\begin{equation}\label{fsi:nsale}
	\rho_f J \big(\partial_t \vv + \hat\nabla \vv F^{-1} (\vv-\partial_t \u)\big)
	-\widehat{\div}\big(J\hat\bsigma_f F^{-T}\big) = \rho_f J \f,
	\quad 
	\widehat{\div}\big(JF^{-1}\vv\big) =0\text{ on }\FL,
\end{equation}
where we denote by an hat the derivatives with regard to the reference system and with the Cauchy stress in ALE coordinates as
\[
\hat\bsigma_f = \rho_f\nu_f(\hat\nabla\vv F^{-1} + F^{-T}\hat\nabla\hat \vv^T)-pI
\]
Further, $F:=\hat\nabla T$ is the deformation gradient of the ALE map and $J=\det(F)$ its gradient. The term $J\hat\bs_fF^{-T}$ is called the Piola transform~\cite[Lemma 2.12]{Richter2017}. 

The classical ALE approach is an iterative algorithm, where the fluid mesh is updated according to the ale map $T$ after every time step~\cite{DoHuePoRo04}. Then, if the geometry is approximated explicitly, the fluid problem simplifies to
\begin{equation}
\rho_f\big(\partial_t\vv + \big((\vv-\partial_t \u)\cdot\nabla\big) \vv\big)-\div\,\bsigma_f = \rho_f \f,\quad \div\,\vv =0\text{ in }\tilde \FL,
\end{equation}
where $\tilde \FL$ is the current domain. We will however base our approach strictly on~\eqref{fsi:nsale} without any remeshing as this simplifies the derivation of the adjoint problem and since this is the only way to achieve higher order accuracy in time. We refer to~\cite{Richter2015} for a detailed discussion of the temporal discretization  \tom{and to~\cite{Margenberg2021} for an analysis of the interplay of temporal and spatial discretization of fluid-structure interactions.}

\subsubsection{Coupling conditions}

The typical coupling conditions request continuity of velocities, the so called \emph{kinematic condition}
\begin{equation}\label{fsi:kinematic}
\vv|_{\FL} = \vv|_{\SO}\text{ on }\IN,	
\end{equation}
and continuity of normal stresses, the \emph{dynamic condition}, which, in ALE coordinates, is given as
\begin{equation}\label{fsi:dynamic}
F\bSI\vec n_f + J\hat\bsigma_f F^{-T}\vec n_s =0.
\end{equation}
By $\vec n_f=-\vec n_s$ we denote the outward facing normal vectors at the interface in reference state. 
The continuity of deformations which is part of the ALE map definition~\eqref{fsi:alemap:2} is often denoted the \emph{geometric condition}.

\subsubsection{Variational formulation}

We derive a monolithic variational formulation of the fluid-structure interaction problem that embeds the coupling conditions in the functions spaces by choosing one global space for the velocity on $\Omega = \FL\cup\IN\cup\SO$ such that
\[
\vv\in \V := H^1_0(\Omega;\Gamma_D^v)^d,
\]
where $\Gamma_D^v\subset\partial\Omega$ is the part of the outer boundary, either of the fluid- or the solid-domain, where the velocity has Dirichlet conditions and one global deformation
\[
\u\in \W :=H^1_0(\Omega;\Gamma_D^u)^d,
\]
where again $\Gamma_D^u$ is the part of the boundary, where the deformation has Dirichlet conditions. $\Gamma_D^v\neq \Gamma_D^u$ is often needed, e.g., if the fluid domain has a Neumann outflow boundary, hence $\vv\neq 0$, but the domain itself is fixed and shall not move, i.e. $\u=0$. Hereby, the kinematic and the geometric coupling conditions are strongly embedded in the function space. The velocity space $\V$ does not properly fit to the regularity of the solid's velocity, as the equation $\vv=d_t\u$, only gives $L^2$-control. We do not stress this difficult point here and refer to the sparse literature on the analysis of the fluid-structure interaction problem for details on the theory and the \twick{regularity~\cite{CouShk05,CouShk06}}.

As common test space for fluid's and solid's momentum equation use $\V$. Hereby the dynamic coupling condition results as natural condition by integration by parts, see~\cite[Sec. 3.4]{Richter2017}. Together with the space $Q_f=L^2(\FL)$ for the divergence, $\V_f=H^1_0(\FL)^d$ for the extension of the deformation to the fluid and $\L_s=L^2(\SO)^d$ for the solid's velocity-deformation coupling we obtain the variational problem:
\begin{equation}\label{fsi:var}
\begin{aligned}
\vv,\u,p&\in \V\times \W\times Q:\qquad
A(\vv,\u,p)(\phi,\psi_f,\psi_s,\xi) = F(\phi,\psi)\qquad
\forall \phi,\psi_f,\psi_s,\xi\in \V\times \V_f\times \L_s\times Q_f \\
A(\vv,\u,p)&:= 
\big(JF^{-T}:\hat\nabla\vv,\xi\big)_{\FL}
+ \big(\hat\nabla \u,\hat\nabla\psi_f\big)_{\FL}
+\big(d_t \u - \vv,\psi_s\big)_{\SO}\\
&\quad 
+\big(
\rho_fJ\big(\partial_t \vv + \hat\nabla\vv F^{-1}(\vv-\partial_t \u),\phi\big)_{\FL}
+ \big( \rho_s d_t \vv,\phi\big)_{\SO}
+ \big(J\hat\bsigma_f F^{-T},\hat\nabla\phi\big)_{\FL}
+ \big(F\bSI_s,\hat\nabla\phi\big)_{\SO}\\
F(\phi) &:= \big(J\rho_f \f,\phi\big)_{\FL}
+\big(\rho_s \f,\phi\big)_{\SO}
\end{aligned}
\end{equation}
Here we used the relation
\begin{equation}
\widehat{\div}(JF^{-1}\vv) = JF^{-T}:\hat\nabla\vv	
\end{equation}
with the product $A:B=\sum_{ij} A_{ij}B_{ij}$ to reformulate the divergence term, see~\cite[Lemma 2.61]{Richter2017} and avoid the appearance of the deformation's second derivatives. In the following we will skip all the hats in the ALE formulation. 

The setup of trial and test functions gives reason to concern, as the ``global space'' $\V$ balances the spaces $\V_f$ and $\L_s$. $\psi_s\in \L_s$ only gives $L^2$-control which is not sufficient to define a trace suitable as boundary condition for the kinematic coupling. Instead of going into detail we again refer to some details on the analysis~\cite{CouShk05,CouShk06}. However, this mismatch in function spaces also affects the adjoint problem discussed in Section~\ref{sec:fsi:adj}.  

\subsection{Discretization}

Discretization of the coupled fluid-structure problem~\eqref{fsi:var} is fairly standard by choosing discrete subspaces $\V^h\subset \V$, $\V_f^h\subset \V_f$, $\L_s^h\subset \L_s$ and $Q_f^h\subset Q_f$. The pair $\V^h\times Q_f^h$ must either be inf-sup stable or stabilization terms must be added to the system of equations. Usually the same space is considered for $\vv$ and $\u$. We refer to~\cite[Sec. 5.3]{Richter2017} for a detailed description of the finite element discretization. 

Discretization in time is more subtle, due to the non-standard appearance of nonlinearities coupling to the time derivatives. Once, the term $\rho_fJ\partial_t \vv$ includes the functional determinant of the ALE map $J=\det(I+\hat\nabla\u)$ and the domain convection term $-\rho_fJ\hat\nabla\vv F^{-1}\partial_t\u$ introduces even further couplings. It is hence not possible to derive a formulation where the temporal derivatives are isolated. We refer to~\cite{Richter2015} for a discussion. 

The resulting nonlinear problems can be efficiently solved in Newton iterations but the linear systems are large, ill-conditioned and lack any desirable structure such as positivity or symmetry, see~\cite{Richter2015a}. Generally speaking, efficient solvers must employ some degree of segregation of fluid and solid to counteract the overall condition number~\cite{Failer2020,Boffi2024}.   

\subsection{Dual Weighted Residual method for fluid-structure interaction}

Basically, the monolithic fluid-structure problem in ALE coordinates and their discretisation with finite elements is a quite common problem and the DWR method can be applied directly. However, the difficulty arises from the complexity of the equations with strong nonlinearities (e.g. in the time derivative) and the coupling conditions. These play a special role in the derivation of the dual equations, as the dual flow of information must be correctly reproduced.  
 
Pioneering work has been done by van der Zee~\cite{Zee2011}, who studied two different approaches to derive the correct conditions. First, leaving the problems in their own coordinate system and using the concept of shape derivatives to obtain gradient information on the interface and second, going via a transformation to a common reference system and taking the gradients with respect to this map. 
This second approach is similar to just considering the derivatives with respect to all quantities in the monolithic ALE formulation~\eqref{fsi:var} and this shows a peculiarity of the FSI system, since the ALE transformation is only an artificially introduced auxiliary variable, but different extensions of $u|_{\SO}$ to $u|_{\FL}$ should lead to the same physical solution. After discretisation, of course, this strict equality no longer applies. 

In the following, we discuss the adjoint of the FSI problem for reduced and simplified problems. The full adjoint of~\eqref{fsi:var} 
\franz{is described either in~\cite{Richter2012} or in~\cite[Sec. 8.3]{Richter2017}.  
}

\subsubsection{The adjoint of a stationary fluid-structure interaction problem}
\label{sec:fsi:adj}

We study a highly simplified problem coupling the vector-valued Poisson equation for the ``fluid-velocity'' $\vv$ to the vector-valued Poisson problem for the ``solid-deformation'' $\u$. We neglect the pressure and all nonlinearities apart from the linearized functional determinant $\tilde J := 1 + \div(\u)$,  to obtain the general structure. In the stationary setting $\vv=0$ in $\SO$:
\begin{multline}
\label{fsi:alelin}
\vv\in \V_f:=H^1_0(\FL)^d,\; \u\in \W:=H^1_0(\Omega)^d:\\
\underbrace{((1+\div\u)\nabla\vv,\nabla\phi)_\FL + (\nabla\u,\nabla\phi)_\SO + (\nabla\u,\nabla\psi_f)_{\FL}}_{=A(\vv,\u)} = (\f,\phi)\quad\forall \phi\in 
H^1_0(\Omega)^d,\; \psi_f\in H^1_0(\FL)^d. 
\end{multline}
This problem corresponds to
\begin{equation}\label{fsi:alelin:strong}
\begin{aligned}	
-\Delta \u = 0,\quad 
	-\div\big((1+\div\u)\nabla\vv\big) &=\f && \text{in }\FL,\\
	-\Delta \u &=\f && \text{in }\SO,\\
	\vv=0,\quad \u|_\FL=\u|_\SO,\quad 
	(1+\div\u)\nabla\vv\vec n_f + \nabla\u\vec n_s &=0&&\text{on }\IN,
\end{aligned}
\end{equation}
which has the same structure and coupling conditions as the stationary fluid-structure interaction problem. Note that $\psi_f=0$ on $\IN$ such that the extension problem is the Poisson problem with Dirichlet conditions $\u|_{\FL} = \u|_{\SO}$. We notice that the stationary problem does not suffer from uncertainties regarding the regularity at the interface, as $\vv|_{\SO}=0$. 
 
Taking the gradient of~\eqref{fsi:alelin} in direction $\delta\vv\in \V_f$ and $\delta u\in \W$ gives
\begin{equation}\label{fsi:alelin_b}
A'(\vv,\u)(\delta \vv,\delta \u;\phi,\psi_f)
=
\big(\div(\delta \u)\nabla\vv,\nabla\phi\big)_\FL +
\big((1+\div\u)\nabla\delta\vv,\nabla\phi\big)_\FL + \big(\nabla\delta\u,\nabla\phi\big)_\SO + \big(\nabla\delta\u,\nabla\psi_f\big)_{\FL}.
\end{equation}
Hereby, the adjoint solution $\phi\mapsto \z\in H^1_0(\Omega)^d$ and $\psi_f\mapsto \w_f\in H^1_0(\FL)^d$ is given as
\begin{equation}\label{fsi:alelineadj}
\underbrace{
\big(\div(\chi)\nabla\vv,\nabla\z\big)_\FL +
\big((1+\div\u)\nabla\zeta_f,\nabla\z\big)_\FL +
\big(\nabla\chi,\nabla\z\big)_\SO + \big(\nabla \chi,\nabla\w_f\big)_{\FL}}_{=:A'(\vv,\u)(\chi,\zeta_f;\z,\w)}
= j(\vv,\u)(\chi,\zeta_f)\quad\forall 
\chi\in \V_f,\;\zeta_f\in \W,
\end{equation}
with the test functions $\delta\u\mapsto \chi\in \W$ and $\delta\vv\mapsto \zeta_f\in \V_f$.

The 
\twick{strong} formulation corresponding to~\eqref{fsi:alelineadj} gets
\begin{equation}
\begin{aligned}
-\div\big((1+\div\u)\nabla\z\big) &= j_{\zeta_f}(\vv,\u)
&&\text{in }\FL\\
-\Delta \w_f-\nabla(\nabla\vv:\nabla\z) &= j_\chi(\vv,\u)
&&\text{in }\FL\\
-\Delta \z &= j_\chi(\vv,\u)
&&\text{in }\SO\\
\z|_\FL = \z|_{\SO},\quad \w_f=0,\quad 
\partial_{n_f}\w_f + \partial_{n_s}\z + (\nabla\vv:\nabla\z)\vec n_f &=0
&&\text{on }\IN.
\end{aligned}
\end{equation}
The structure of this simplified adjoint transfers to the full nonlinear fluid-structure interaction problem in the stationary setting. Here, $\z|_\FL$ can be computed in a ``nearly decoupled'' way as extension of $\z|_{\SO}$ to $\FL$. The only feedback to $\SO$ is by means of the boundary term $(\nabla\vv:\nabla\z)\vec n_f$, which arises from the ALE deformation. Further similar terms appear when the complete nonlinear system is considered. \tom{By $j_{\zeta_f}(\cdot,\cdot)$ and $j_\chi(\cdot,\cdot)$ we denote the functional in strong formulation, i.e. $\langle j_{\zeta_f}(\vv,\u),\chi_f\rangle+\langle j_{\chi}(\vv,\u),\chi\rangle = j(\vv,\u)(\chi,\zeta_f)$.}

\twick{
An open-source implementation of a stationary fluid-structure interaction problem with the primal residual of the DWR estimator can be found on github\footnote{\url{https://github.com/tommeswick/goal-oriented-fsi}}~\cite{Wi21_git_YIC}, which is based on the finite element library deal.II~\cite{dealII96,deal2020}. This code serves as basis for a FSI multigoal numerical example that is conducted in Section~\ref{sec:Results}.
}

\subsubsection{The adjoint of a nonstationary fluid-structure interaction problem}

Now, we study a simplified nonstationary problem. Again, we neglect most of the ALE coupling terms, but this time we include the domain convection to keep the overall nature of the couping:%
\begin{equation}
\label{fsi:nsalelin}
\big(\partial_t \vv - \u\cdot\nabla\vv,\phi\big)_\FL+
\big(\nabla\vv,\nabla\phi\big)_\FL 
+\big(\partial_t \vv,\phi\big)_\SO
+ \big(\nabla\u,\nabla\phi\big)_\SO + 
\big(\nabla\u,\nabla\psi_f\big)_{\FL}
+\big(\partial_t \u-\vv,\psi_s\big)_{\SO}
= \big(\f,\phi\big).
\end{equation}
We do not specify exact function spaces, but consider the problem in a space-time variational formulation. 
This problem corresponds to
\begin{equation}
\begin{aligned}	
-\Delta \u=0,\quad
\partial_t \vv - \u\cdot\nabla\vv - \Delta \vv &= \f&&\text{in }\FL\\
\partial_t \u=\vv,\quad 
\partial_t \vv -\Delta \u &=f&& \text{in }\SO,\\	
\u|_\FL = \u|_\SO,\quad
\vv|_\FL = \vv|_\SO,\quad
\partial_{n_f}\vv + \partial_{n_s}\u &=0
&&\text{on }\IN.
\end{aligned}
\end{equation}
Again, the overall coupling structure matches that of the full nonlinear coupled problem. The gradient is defined by the form
\begin{equation}
\begin{aligned}
A'(\vv,\u)(\delta\vv,\delta\u;\phi,\psi_f,\psi_s) &=
\big(\partial_t \delta\u - \delta\u\cdot\nabla\vv-\u\cdot\nabla\delta\vv,\phi\big)_\FL+
\big(\nabla\delta\vv,\nabla\phi\big)_\FL 
+\big(\partial_t \delta\vv,\phi\big)_\SO
+ \big(\nabla\delta \u,\nabla\phi\big)_\SO\\
&\qquad +
\big(\nabla\delta \u,\nabla\psi_f\big)_{\FL}
+\big(\partial_t \delta \u-\delta\vv,\psi_s\big)_{\SO}.
\end{aligned}
\end{equation}
We define the adjoint variables $\phi\mapsto \z\in H^1_0(\Omega)^d$, $\psi_f\mapsto \w_f\in H^1_0(\FL)^d$ and $\psi_s\mapsto \w_s\in L^2(\Omega)^d$ and state the adjoint problem, sorted by the test functions $\delta \u\to\chi\in H^1_0(\Omega)^d$, $\delta \vv\to \zeta\in H^1_0(\Omega)^d$, as
\begin{multline}
\big(\partial_t \chi - \chi\cdot\nabla\vv,\z\big)_\FL
-\big(\u\cdot\nabla\zeta,\z\big)_\FL+
\big(\nabla\zeta,\nabla\z\big)_\FL 
+\big(\partial_t \zeta,\z\big)_\SO
+ \big(\nabla\chi,\nabla\z\big)_\SO\\
+
\big(\nabla\chi,\nabla\w_f\big)_{\FL}
+\big(\partial_t \chi,\w_s\big)_{\SO}
-\big(\zeta,\w_s\big)_{\SO}  = j(\vv,\u)(\chi,\zeta).
\end{multline}
The corresponding classical form gets
\begin{equation}
	\begin{aligned}
 -\Delta \z + \div(\z\otimes \u)&= j_\zeta(\vv,\u),&\quad 
	-\partial_t \z - \nabla\vv^T\z 	- \Delta \w_f &=j_\chi(\vv,\u)&&\text{in }\FL,\\
-\partial_t \z -\w_s &=j_\zeta(\vv,\u),&\quad 
-\partial_t \w_s	-\Delta \z &= j_\chi(\vv,\u)&&\text{in }\SO,
\\
\w_f=0,\quad
-(\z\otimes\u)\vec n_f &+\partial_{n_f}\z =0,
&
\partial_{n_f}\w_f+\partial_{n_s}\z=0
,\quad
 \z|_{\FL}&=\z|_{\SO}
&&\text{on }\IN.
	\end{aligned}
\end{equation}

\tom{As usual, the system runs backward in time. And as in the stationary case, the adjoint variable $\z|_{\FL}$ is mostly an extension of $\z|_{\SO}$ into the solid domain and feedback to the solid is only by means of the ALE mapping, here shown in the term $-(\z\otimes\u)\vec n_f$.
  The role of the adjoint coupling conditions and the flow of information between the two subproblems in the context of optimization problems is analyzed in~\cite{RichterWick2013,FailerRichter2019opt}. 
}


\section{Numerical experiments} 
\label{sec:Results}

\franz{In this section, we present the performance of the DWR method in controlling the discretization error in simulations of soft-tissue biomechanics.
%
The simulations \twick{in Section~\ref{sec_example_artery} and Section~\ref{sec_example_silicone}} have been realized thanks to the \texttt{python} library \texttt{FEniCS} and the code is available online~\cite{Duprez2023}.}
\twick{The computations in Section~\ref{sec_example_FSI_multigoal}
are realized with an extension of the open-source code on github\footnote{\url{https://github.com/tommeswick/goal-oriented-fsi}}~\cite{Wi21_git_YIC}, which is based on the finite element library deal.II~\cite{dealII96,deal2020}. }



\subsection{First test case: human artery with fiber activation}\label{sec:artery}
\label{sec_example_artery}
\begin{figure}[!ht]
\begin{center}
\includegraphics[scale=0.80]{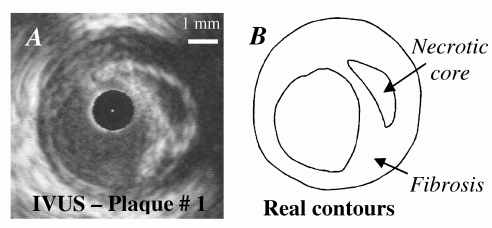} 
\includegraphics[scale=0.12]{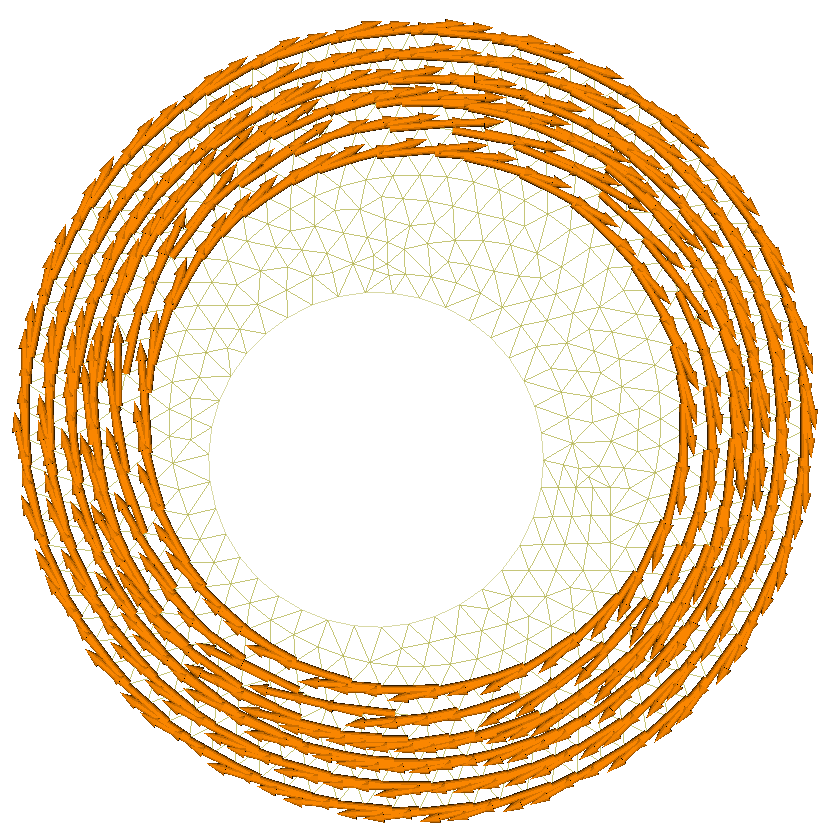} 
\includegraphics[scale=0.18]{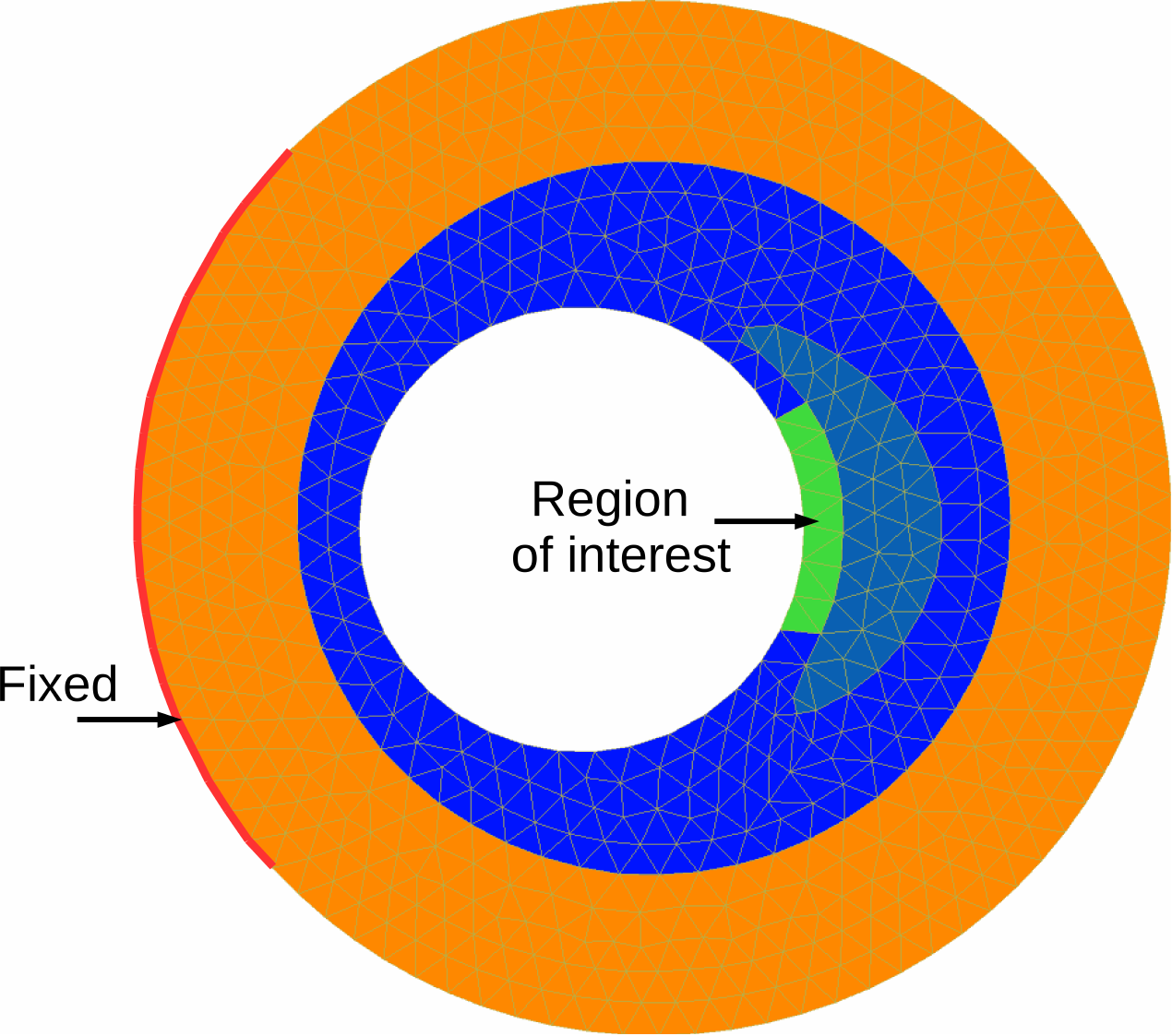}
\end{center}
\caption{Artery {model}. Geometry, from \cite{floch_vulnerable_2009} (left), fiber orientation
 (center) and region of interest (right).}
\label{fig:artery medical img}\end{figure}

\begin{figure}[!ht]
\begin{center}
\includegraphics[scale=0.08]{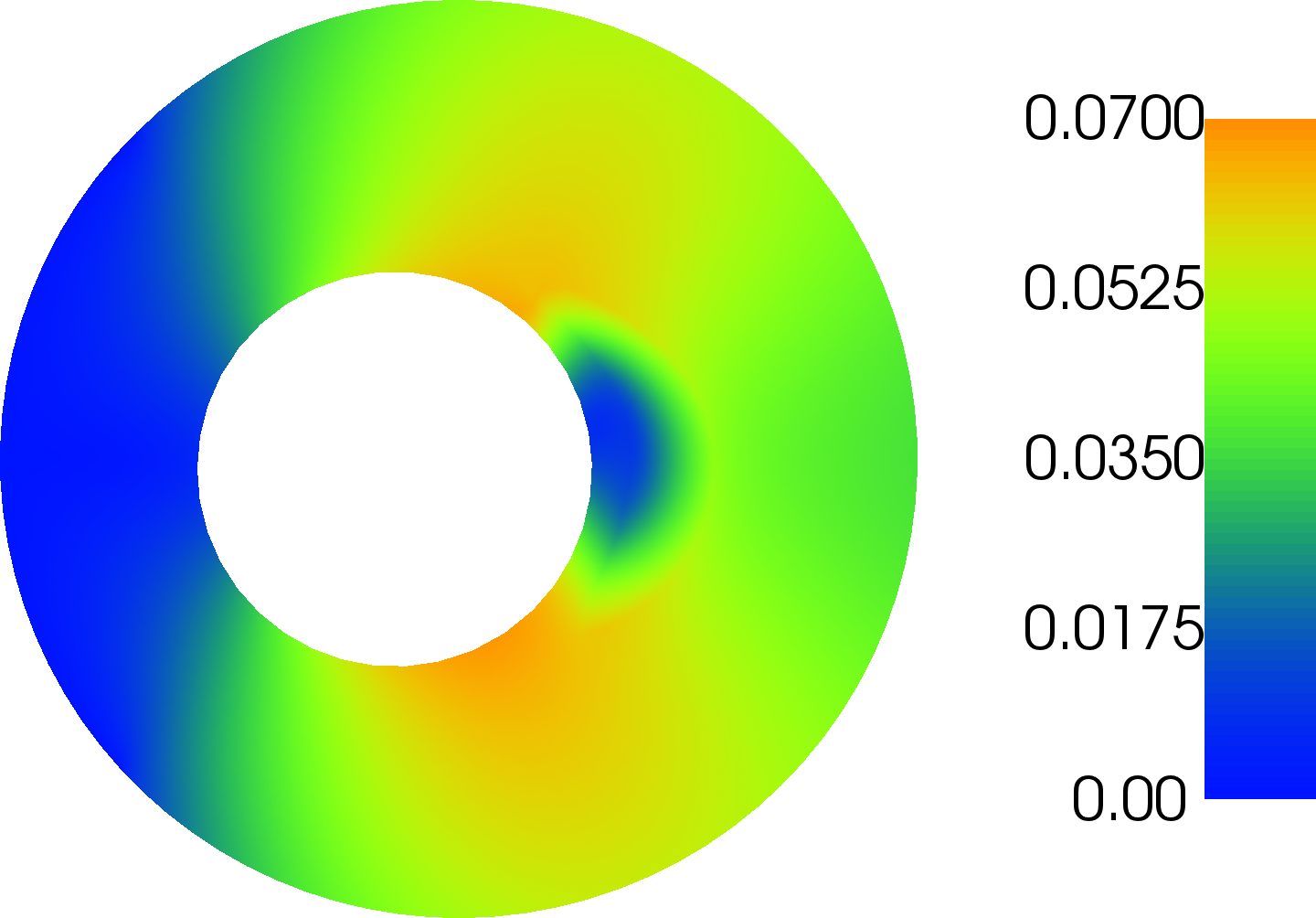}\hspace{1.0em}
\includegraphics[scale=0.08]{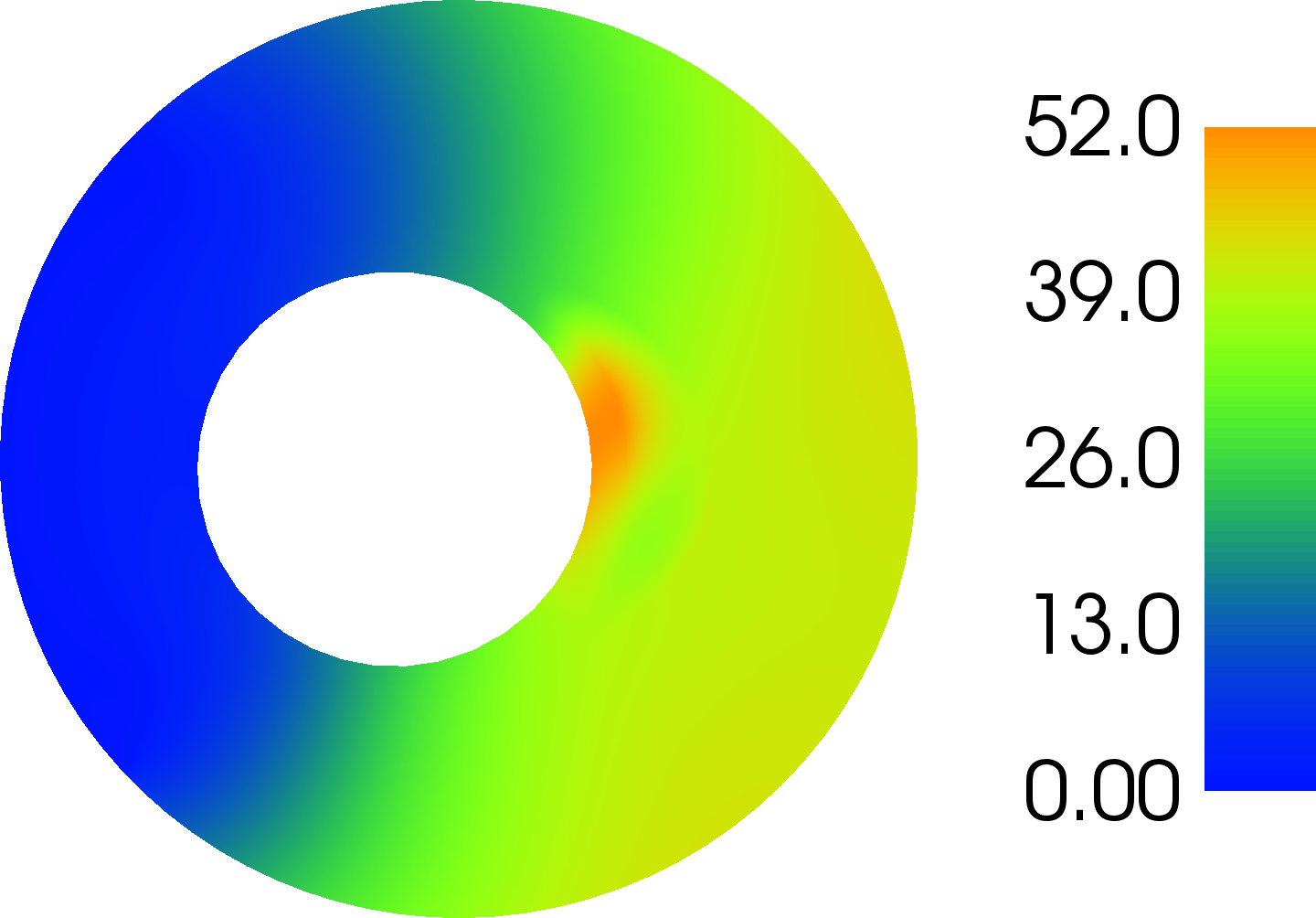}
\end{center}
\caption{Artery model. 
\franz{Displacement (left)
and dual solution for $J$ (right).}}
\label{fig:artery dual}\end{figure}

\franz{{As a first example, we showcase the performance of the method in the full linear setting (with small deformations) for the analysis of the mechanical response of an artery with vulnerable coronary plaque to internal loading.
This example comes from \cite{duprez2020}.
Rupture of the cap induces the formation of a thrombus which may obstruct the coronary artery, cause an acute syndrome and the patient death.} The geometry (see Figure \ref{fig:artery medical img} (left)) comes from \cite{floch_vulnerable_2009} where the authors develop a methodology to reconstruct the thickness of the necrotic core area and the calcium area as well as Young's moduli of the calcium, the necrotic core and the fibrosis. Their objective is the prediction of the vulnerable coronary plaque rupture.
{As represented in Figure \ref{fig:artery medical img} (left),
the diameter of the Fibrosis is equal to 5 mm.} 
{Following \cite{floch_vulnerable_2009}, we set different elastic parameters in each region: {$E_n=0.011$ MPa, $\nu=0.4$ in the necrotic core and $E_s=0.6$ MPa, $\nu=0.4$ in the surrounding tissue (contrast $E_s / E_n \simeq 55$)}.}
{No {volumetric} force field is applied: $\f=\0$.} 
We consider 
muscle fibers only in the media layer, where
smooth muscle cells are supposed to be perfectly oriented in the circumferential direction $\eb_A = \eb_\theta$, where $(\eb_r,\eb_\theta)$ is the basis for polar coordinates{{, see Figure \ref{fig:artery medical img} (center)}}. 
{Other parameters for fiber activation have been chosen as {$T = 0.01$ {MPa}} and $\beta=1$. 
}
As depicted in Figure \ref{fig:artery medical img} (right), the artery is fixed on the red portion {of} external boundary $\Gamma_D$. 
{Elsewhere, on the remaining part of the boundary, a homogeneous Neumann condition is applied: $\F=\0$. }
In the same figure, the green part represents the region of interest $\omega$, which has been defined in order to be relevant in the study of vulnerable coronary plaque rupture.
The quantity of interest for this example is:
\begin{equation}\label{def:quantities}
J(\u):=\displaystyle\int_{\omega} ( u_x+u_y ) \, d\x.
\end{equation}
Figure \ref{fig:artery dual} 
depicts the magnitude of the solution in terms of displacements (left)
and the dual solution associated to $J$ (right).
\revision{In this example, the dual solution is computed directly using the finer finite element space.}}

\begin{figure}[!ht]
\begin{center}
\includegraphics[scale=0.119]{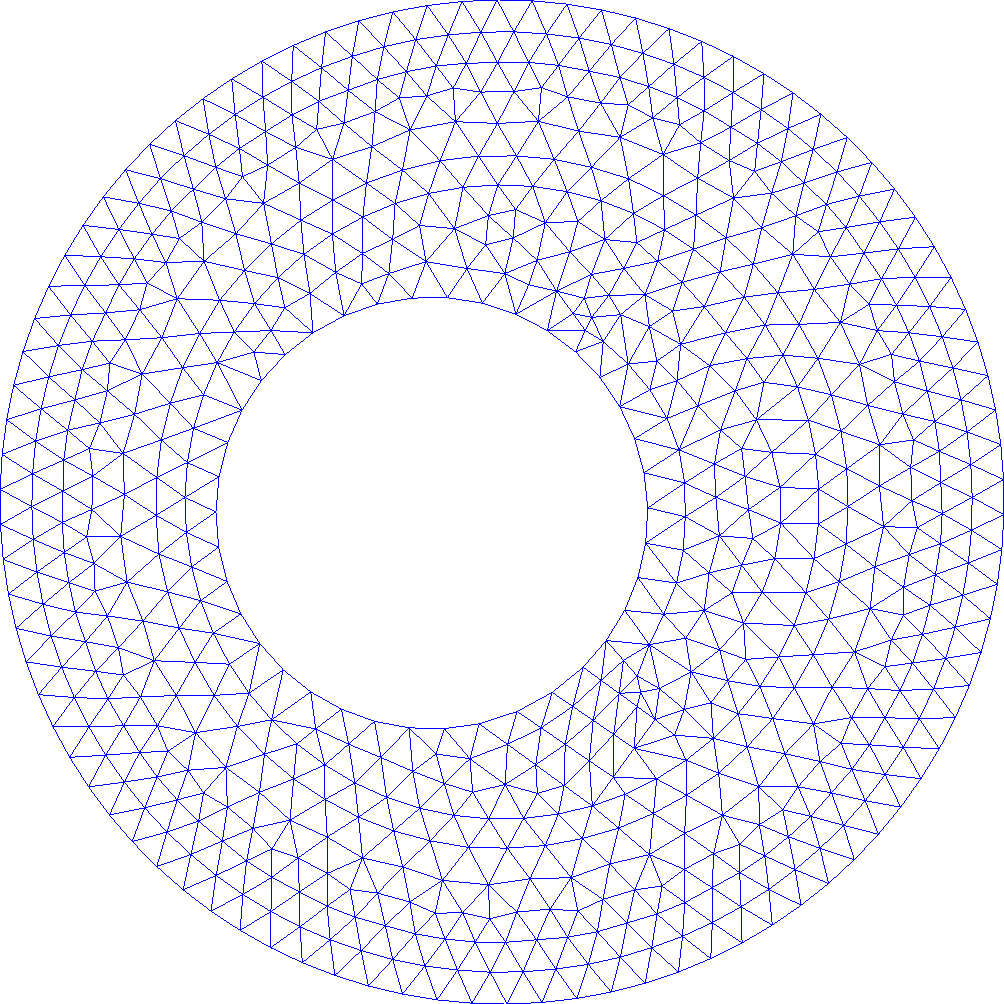}
\hspace{1em}
\includegraphics[scale=0.119]{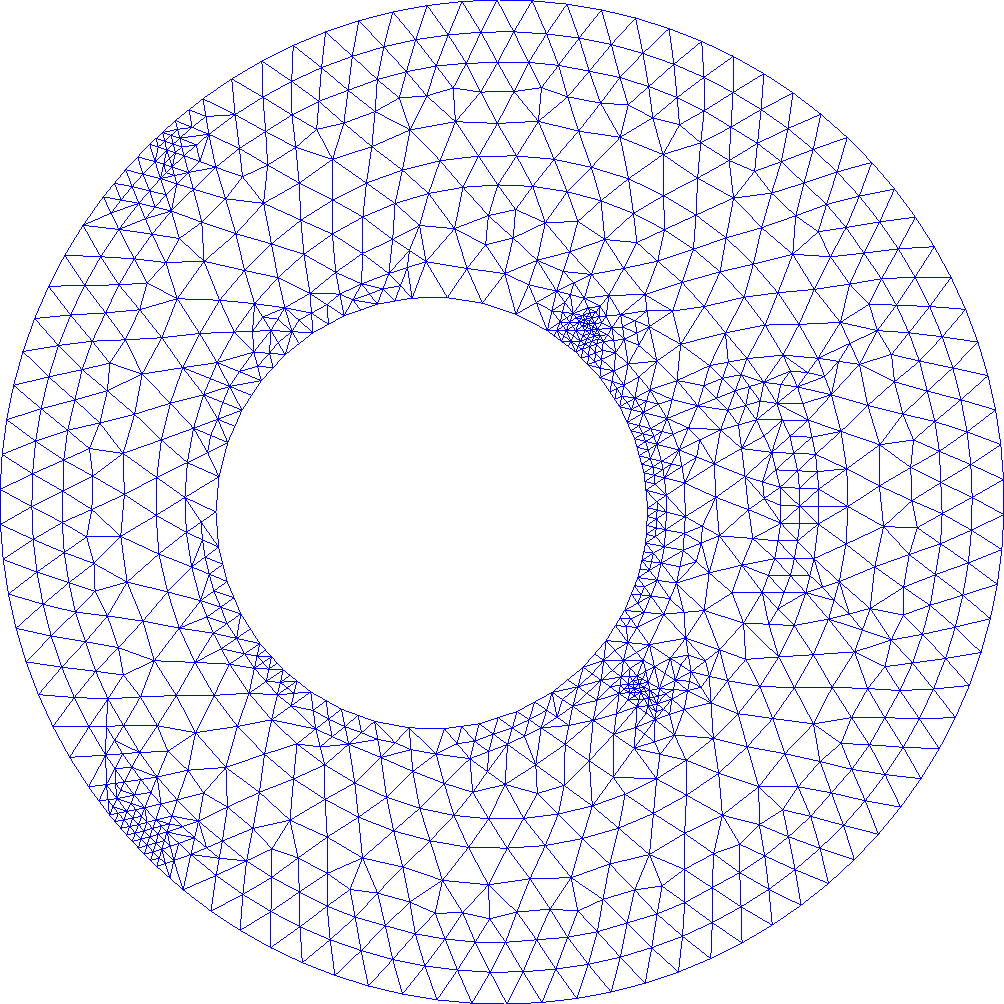}
\hspace{1em}
\includegraphics[scale=0.119]{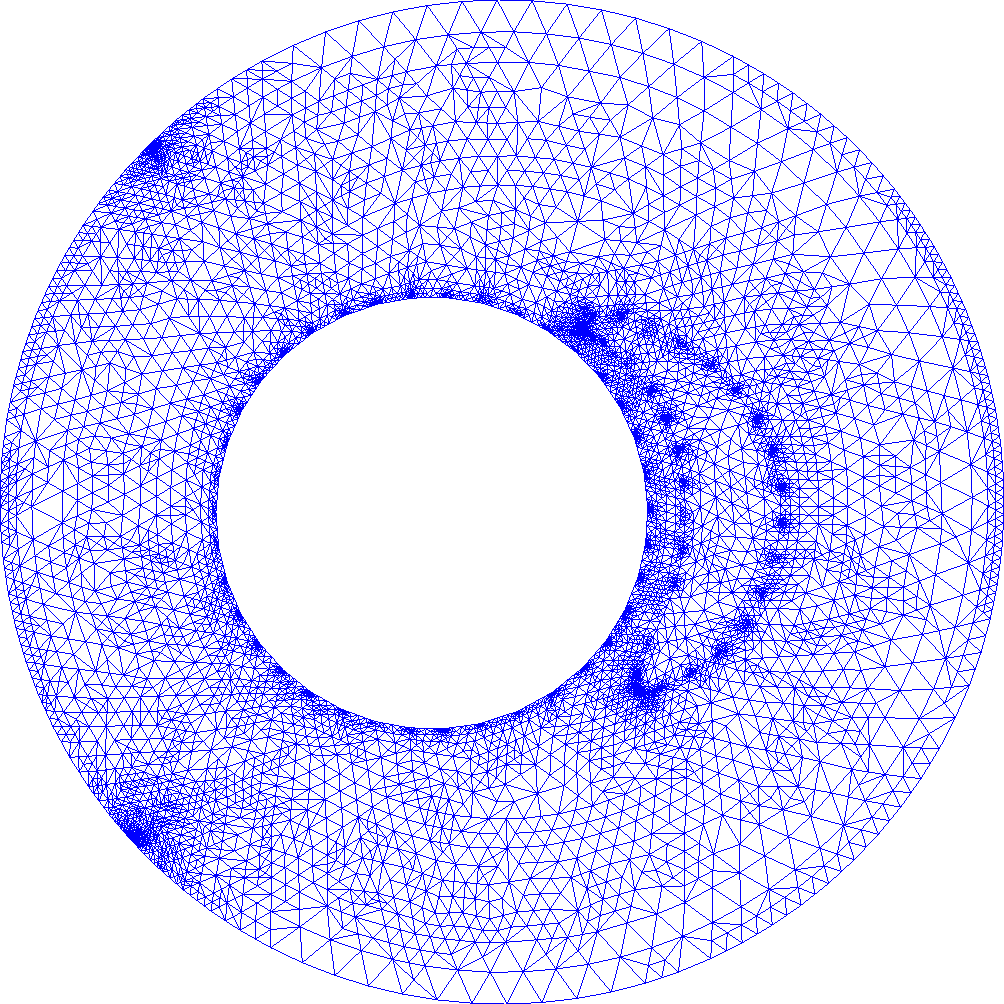} 


\end{center}
\caption{
Artery mesh. 
\franz{Refinement driven by the QoI $J$.
Initial mesh (left) with 1242 cells and a relative error of {$38.3$ \%}, adapted meshes after 2 iterations (center) with 2079 cells and a relative error of {$5.2$ \%} and after 6 iterations (right)
 with 15028 cells and a relative error of {3.4 \%}.}
%
 }
\label{fig:artery meshes}
\end{figure}

\franz{In Figure \ref{fig:artery meshes}, we present the mesh after 2 and 6 iterations 
for the quantity of interest $J$. 
The refinement occurs in some specific regions{,} such as {those} near Dirichlet-Neumann transitions and concavities on the boundary. {Our results also show that the proposed method leads to}  the strong refinement near the interface between the necrotic core and the fibrosis, where stresses are localized because of the material heterogeneity.} 



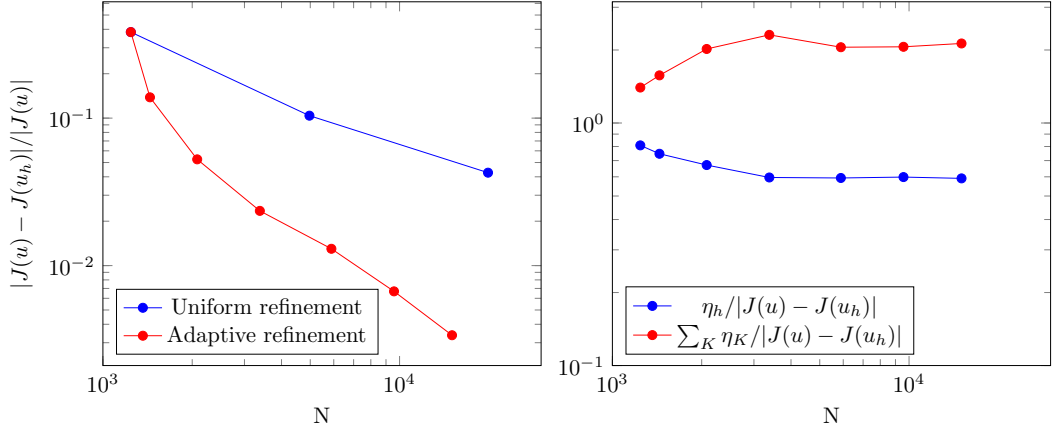
\begin{figure}[!ht] 
\begin{center}
\begin{tikzpicture}[thick,scale=0.85, every node/.style={scale=1.0}] \begin{loglogaxis}[xlabel=N,ylabel=$|J(u)-J(u_h)|/|J(u)|$,
xmin=1e3,xmax=3e4,xtickten={3,4}
,legend pos=south west, legend columns=1]
 \addplot[color=blue,mark=*] coordinates { 
(1242.0,0.383042850931)
(4968.0,0.103889179043)
(19872.0,0.0426608966597)
 };
  \addplot[color=red,mark=*]coordinates { 
(1242.0,0.383042850931)
(1441.0,0.138364267347)
(2079.0,0.0525334816088)
(3381.0,0.023503151137)
(5891.0,0.0130017151056)
(9578.0,0.00668885135909)
(15028.0,0.00337139785765)

 };
 \legend{ Uniform refinement, Adaptive refinement}
\end{loglogaxis} 
\end{tikzpicture} 
\begin{tikzpicture}[thick,scale=0.85, every node/.style={scale=1.0}] \begin{loglogaxis}[xlabel=N,
xmin=1e3,xmax=3e4,ymin=1e-1,xtickten={3,4}
,legend pos=south west, legend columns=1]
 \addplot[color=blue,mark=*] coordinates { 
(1242.0,0.808706493463)
(1441.0,0.746815968529)
(2079.0,0.67093292762)
(3381.0,0.596218349801)
(5891.0,0.593460833375)
(9578.0,0.598516436079)
(15028.0,0.591185139005)
 };
  \addplot[color=red,mark=*]coordinates { 
(1242.0,1.40028694858)
(1441.0,1.57197067518)
(2079.0,2.0195018399)
(3381.0,2.30936216451)
(5891.0,2.05327320072)
(9578.0,2.06154249942)
(15028.0,2.12780291125)
  };
 \legend{ $\eta_h/|J(u)-J(u_h)|$, $ \sum_K\eta_K/|J(u)-J(u_h)|$ }
\end{loglogaxis} 
\end{tikzpicture}
\end{center}
\caption{\franz{Artery {model}. Left: relative error for 
the quantity of interest $J$ \emph{vs.} the number $N$ of cells in the case of uniform (blue) and adaptive (red) refinement.
Right: effectivity indices of $\eta_h$ (blue) and $\sum_K\eta_K$ \emph{vs.} the number of cells $N$ for the quantity of interest $J$.} 
}
\label{fig:artery graph eff}
\end{figure}

\franz{Figure \ref{fig:artery graph eff} (left) depicts the relative goal-oriented error 
$|J(\u) - J(\u_h)| / |J(\u)|$
versus the number $N$ of cells in the mesh, both for uniform refinement {(blue)} and adaptive refinement (red). \revision{Here the quantity $J(u)$ has been approximated using an overrefined mesh.}
The stopping criterion $\varepsilon$ has been fixed at {$5.10^{-6}$}. 
{Remark that, for the intial mesh $N=1242$, the relative value of the discretization error is large (about $38$ \%), because the mesh does not resolve properly the discontinuity of material parameters $E_n$ and $E_s$ at the boundary of the necrotic core. 
This is typically a situation where the discretization error is not negligible relatively to the modelling errors, even for a mesh with more than thousand cells. Of course, the modelling errors are not quantified here, but there are generally considered acceptable when they are below 20\%.
The adaptive algorithm allows to recover this interface, as illustrated by Figure \ref{fig:artery meshes}, and to reduce the error, which is around 5 \% after only two iterations.}
In Figure \ref{fig:artery graph eff} (right), we depict the effectivity {indices} for the global estimator $\eta_h$ and the sum of local estimators $\sum_K \eta_K$. 
A recommandation for the practitionner in such a situation would be to carry out 2 iterations of refinements, which allows to decrease the discretization error below an acceptable threshold, while keeping a reasonable use of computational ressources and a reasonable increase of the degrees of freedom. Refining too much may not be necessary since it increases a lot the degrees of freedom, for a moderate gain in terms of accuracy.}

\subsection{Second test case: silicone samples}
\label{sec_example_silicone}
\begin{figure}[!ht]
     \centering
~\hfill  \begin{minipage}{0.25\linewidth}
  \includegraphics[width=1\textwidth]{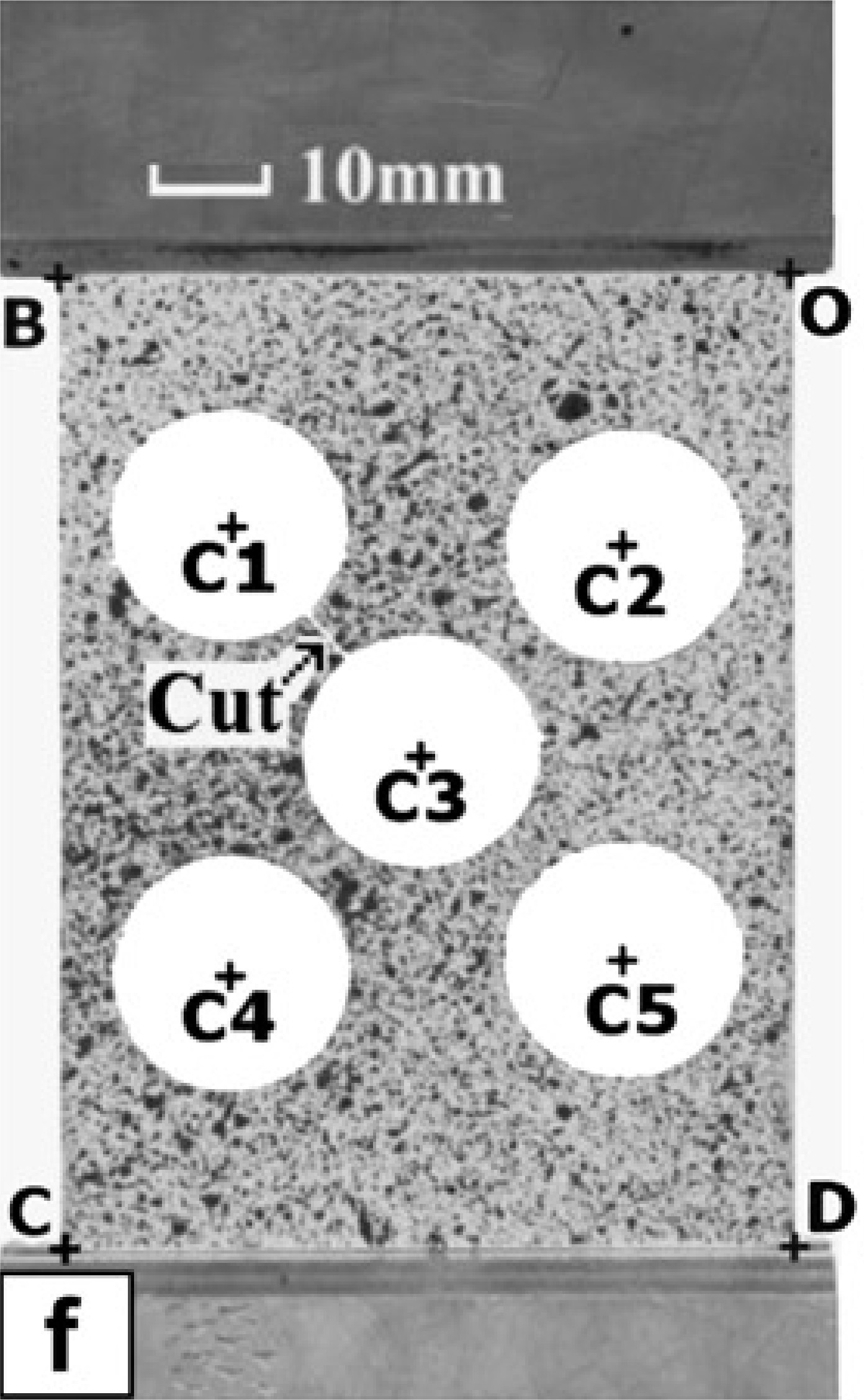}\vspace*{-0cm}
\end{minipage}\hfill
    \begin{minipage}{0.3\linewidth}
     \begin{tabular}{|c|c|c|}\hline
          &X (mm)&Y (mm)  \\\hline
          O& 0&0\\ \hline
        B& -62&0\\ \hline
      C& -62&-82.5\\ \hline
D& 0&-82.5\\ \hline
          $C_1$& -47.5&-21.4\\ \hline
    $C_2$& -14.0&-23.0\\ \hline
      $C_3$& -31.5&-41.0\\ \hline
        $C_4$& -47.7&-59.1\\ \hline
        $C_5$& -14.5&-58.0\\ \hline
     \end{tabular}\end{minipage}\hfill
     \begin{minipage}{0.25\linewidth}
     \includegraphics[width=1\linewidth]{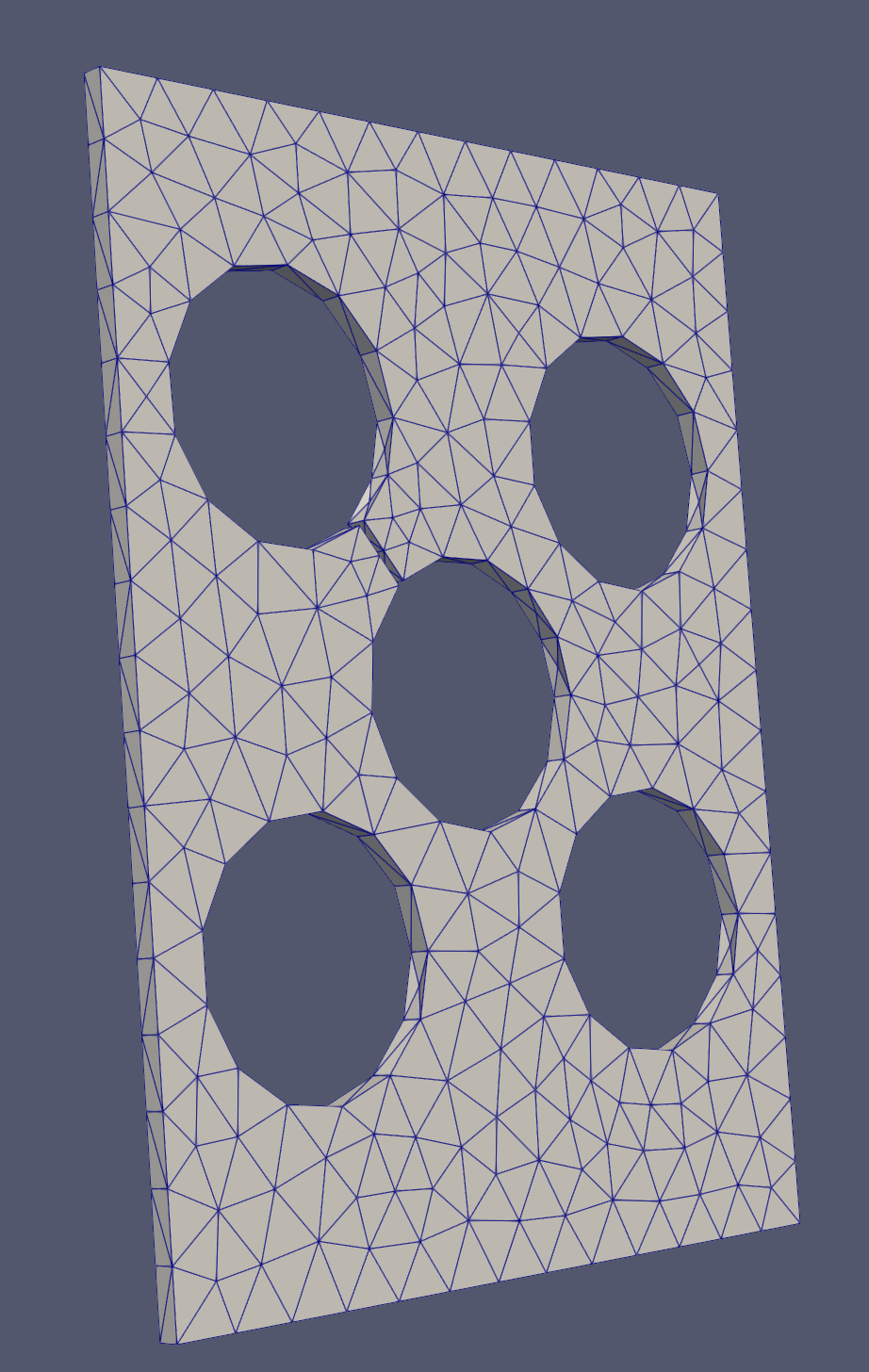}\vspace*{-0.3cm}
     \end{minipage}\hfill~
     \caption{Geometry of the silicone sample, position of the holes and initial mesh.}
     \label{fig:geometry}
    \end{figure}

\franz{Now we provide numerical results for the nonlinear situation of hyperelasticity, with an example based on a silicone sheet, which properties are close to those of soft tissue, while allowing to have reproductible and accurate experimental measurements. The results are a summary from \cite{bui:hal-04208610}.
The experimental procedure is briefly recalled here for the sake of clarity. The reader is referred to \cite{meunier2008mechanical} for more information.
Simple tensile tests are performed on dumbbell-shaped samples of silicone rubber (RTV 141) having an initial gauge length $l_{0}$ of 82.5 mm, a gauge width $b_{0}$ of 61.5 mm, and a gauge thickness $e_{0}$ of 1.75 mm.  The sample contains five holes of diameter 20 mm
and the position of the centers of the holes and the corners are given in Figure \ref{fig:geometry}. There is also a cut between the circles C1 et C3.}

Tested samples are deformed using a universal mechanical testing machine (MTS 4M) (see \cite{meunier2008mechanical}).
Dirichlet boundary conditions are imposed on the bottom edge of the dumbell silicone sample. On the left and right boundaries, we impose a homogeneous Neumann boundary condition ($\F = \0$). In the experiment, the following Neumann boundary condition is imposed on the top edge: 
$\F = \left (f_A / (b_{0} \times e_{0} )\right ) \n$ 
such that $\int_{\Gamma_N} \F \cdot \n \, \mathrm{ds} = 20$ N, with $\n$ the normal vector.
This force implies an observed vertical displacement of 57.3mm. In the simulations, we do it the other way round: we fix the bottom and we impose a displacement of 57.3mm on the top. We guess the corresponding traction force on the top boundary.
\franz{Thus, we will consider the following quantities of interest 
$$J(\bu) 
=\int_{\mbox{top}}(\bPi(\bu) \cdot \n)\cdot \n \,\mathrm{d}\mathrm{s},$$
where the integral is taken on the top of the silicone band.}
\revision{In this example, the dual solution is computed using the primal finite element space and then is extrapolated onto the finer finite element space.}

\begin{figure}[!ht]
    \centering
    \includegraphics[width=0.3\linewidth]{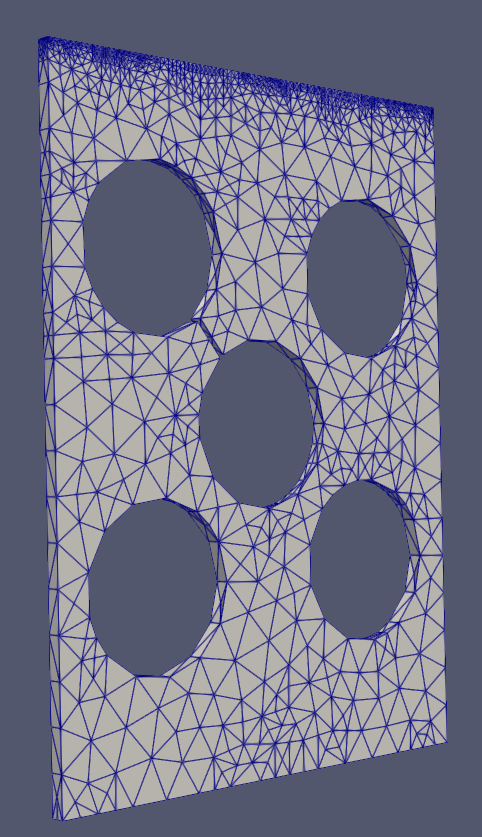}
        \includegraphics[width=0.3\linewidth]{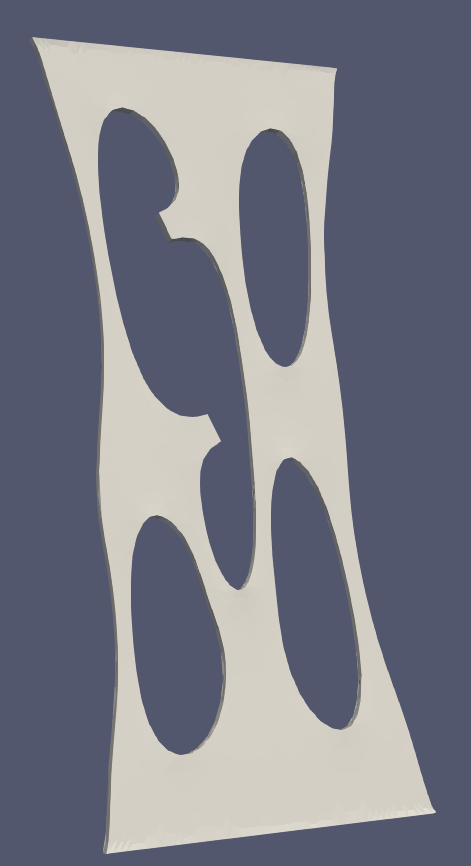}      
        \caption{\franz{Hyperelastic (nonlinear)} test case 
        (silicone sample).  Refined mesh for the Haine-Wilson model (left); Deformed geometry (right).}
    \label{fig:deformed mesh}
\end{figure}

\begin{figure}[!ht]
\centering
    
\begin{tikzpicture}\begin{loglogaxis}[
width = .49\textwidth, xlabel = \# degrees of freedom, ymin=1e-2,xmin=7000,xmax=800000,
            legend style = { font=\tiny,at={(0,1.1)},anchor=south west, legend columns =1,
			/tikz/column 2/.style={column sep = 10pt}}]

\addplot coordinates { 
(9597.0,0.0836046958956285)
(59637.0,0.03854)
(413510.0,0.01605)
};
\addplot coordinates{
(9597.0,0.0836046958956285)
(12568.0,0.08293439020081689)
(20734.0,0.06503395099262627)
(30854.0,0.03629264102734441)
(49369.0,0.025795263320612643)
(86952.0,0.015809590358542867)
(165252.0,0.014201305642279224)
(332323.0,0.01272454857468837)
(619881.0,0.011069776437639513)
};
\legend{\begin{tabular}{l}$|J(\bu_h{,}\bpp_h)-J(\bu_{\mbox{fine}}{,}\bpp_{\mbox{fine}})|/J(\bu_{\mbox{fine}}{,}\bpp_{\mbox{fine}})$ \\ Uniform\medskip\end{tabular},\begin{tabular}{l}$|J(\bu_h{,}\bpp_h)-J(\bu_{\mbox{fine}}{,}\bpp_{\mbox{fine}})|/J(\bu_{\mbox{fine}}{,}\bpp_{\mbox{fine}})$ \\Adaptative\medskip\end{tabular}}
\end{loglogaxis}
\end{tikzpicture}
\quad
\begin{tikzpicture}\begin{semilogxaxis}[width = .49\textwidth, xlabel = \# degrees of freedom, ymin=-1,ymax=2,xmin=7000,xmax=500000,
            legend style = { font=\tiny,at={(0,1.1)},anchor=south west, legend columns =1,
			/tikz/column 2/.style={column sep = 10pt}}]

\addplot coordinates{
(9597.0,0.28269696504426267)
(12568.0,0.5208906772106272)
(20734.0,1.472220761528068)
(30854.0,1.246899661177141)
(49269.0,0.9439882211722505)
(86952.0,0.6725932950904422)
(165252.0,0.7131123753928059)
(332323.0,0.7140406151783296)
};
\legend{
$[J(\bu_h{,}\bpp_h)-J(\bu_{\mbox{fine}}{,}\bpp_{\mbox{fine}})]/\sum_T\eta_T$}
\end{semilogxaxis}
\end{tikzpicture}
\caption{\franz{Hyperelastic (nonlinear)} 
test case (silicone sample). Relative error of discretisation (left) and efficiency of the estimator (right).}
\label{fig:errorH1_manufactured}
\end{figure}
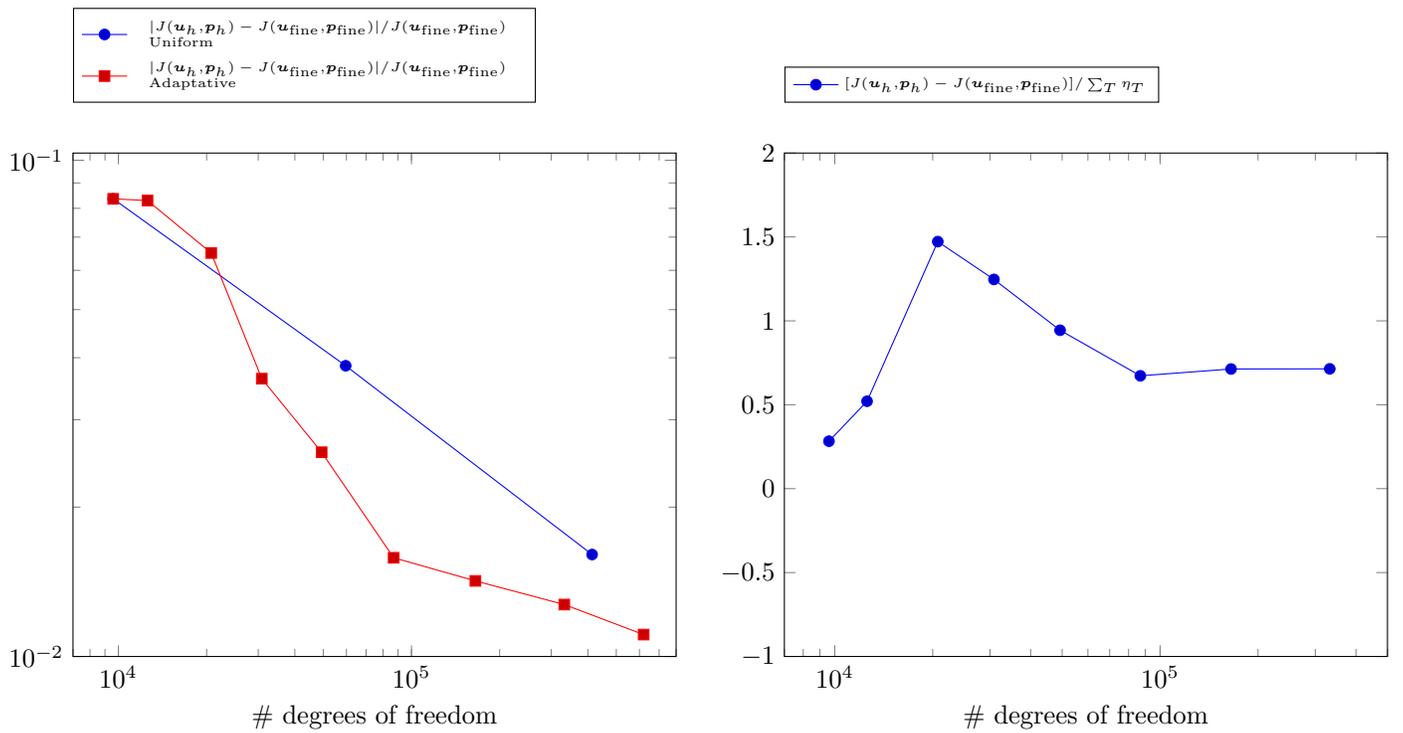

Table \ref{tab:parameter_law} recalls the value (estimated in \cite{meunier2008mechanical}) of the constitutive parameters used in the simulations.
\franz{In Table 
\ref{tab:law3}, we compare the model and discretisation error for the constitutive law of Haine-Wilson, 
for the quantity $J$. The model error corresponds to the relative error between a very fine FEM solution $\bu_{\mbox{fine}}$ {($4.5\times 10^6$ degrees of freedom)} and the experimental quantity of interest. The discretization error corresponds to the relative error between the computed solution on the current mesh and the computed solution on a very fine mesh, i.e.
$${
\mbox{model error}=\dfrac{|J(\bu_{\mbox{fine}})-20|}{20}
\mbox{ and discr. error}=\dfrac{|J(\bu_h)-J(\bu_{\mbox{fine}})|}{J(\bu_{\mbox{fine}})}.
}$$}
We give in Figure \ref{fig:deformed mesh} the refined mesh in the case of the Haine-Wilson law (left) and the deformed geometry when we apply the load (right). We remark that the refinement occurs mainly on the top of the silicone where the quantity of interest is localised but also near the holes.
\franz{The efficiency of the estimator is represented in Figure \ref{fig:errorH1_manufactured}.}

 \begin{table}[]
        \centering
\begin{tabular}{llll}\hline
Mooney&$C_{10}=0.14$&$C_{01}=0.023$&\\ \hline
Gent&$E=0.97$&$J_m=13$&\\ \hline
Haines-Wilson&$C_{10}=0.14$&$C_{20}=-0.0026$&$C_{30}=0.0038$\\
&$C_{01}=0.033$&$C_{02}=0.00095$&$C_{11}=-0.0049$\\\hline
        \end{tabular}       
        \caption{First test case (silicone sample). Values of the constitutive parameters of each hyperelastic law, following \cite{meunier2008mechanical}.}
        \label{tab:parameter_law}
    \end{table}

    \begin{table}[]
        \centering
        
\begin{tabular}{|c|c|c|c|c|}\hline
\multicolumn{5}{|c|}{Haine Wilson}\\
\hline
Number of &Number of & {Number of degrees} & Model&Discr. \\
     Iteration&Cells& {of freedom}& Error& Error\\\hline
1&1325&  9597  &&8.3 \% \\
2&1860&  12568 &  &8.2 \%\\
3&3384&  20734 &  &6.5 \% \\
4&5364&  30854 &&3.6  \%\\
5&9127&  49369 &0.5\%&2.6 \%\\
6&16924&   86952&&1.6 \%\\
7&33392&   165252 &&1.4 \%\\
8&68919&   332323 &&1.3 \%\\
9&130293&   619881 &&1.1 \%\\\hline
        \end{tabular}       
        \caption{\franz{Hyperelastic (nonlinear) test case (silicone sample) {with the Haine Wilson law}.  Relative model and discretisation error of the quantity of interest $Q$ with respect to the number of cells. The model error is computed thanks to the experimental data.} 
        }
        \label{tab:law3}
    \end{table}

\clearpage


\subsection{\twick{Multigoal error control in stationary fluid-structure interaction: stenosis in artery}}
\label{sec_example_FSI_multigoal}
\begin{figure}[!ht]
    \centering
    \includegraphics[width=15cm]{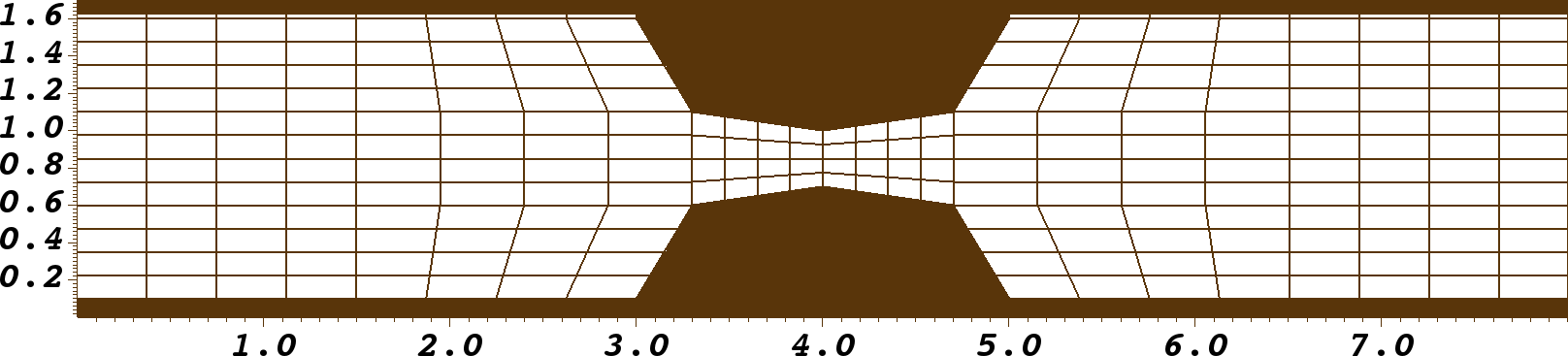}
    \caption{Example multigoal FSI: Geometry and initial mesh. The solid is prescribed in the brown-colored regions. The units are given in cm.}
    \label{fig:fsi_multigoal_a}
\end{figure}

In this section, we employ the multigoal technology from Section~\ref{sec:multigoal_section}. \twick{The following data on the overall domain (geometry) and material parameters are academic choices, but result from discussions with a medical doctor, Jeremi Mizerski, and were used prior in related studies; e.g.,~\cite{Wick2012,Wi14_fsi_eale_heart}.}
The geometry of the domain $\Omega$ and initial mesh are shown in Figure~\ref{fig:fsi_multigoal_a}.
The velocity inflow profile (non-homogeneous Dirichlet condition) is given by 
\begin{equation}
    \vv (t,0,y) = 0.1 (y-0.1)(y-1.6) \quad\text{on } \Gamma_{in} := \{ \x\in\mathbb{R}^2 |\; x=0\text{cm},\; 
    0.1\text{cm}\leq y\leq 1.6\text{cm} \}.
\end{equation}
On the top and bottom boundaries, we prescribe no-slip conditions, $\vv = \mathbf{0}$. At the outflow boundary, we have the do-nothing condition. For the displacements $\u$, we prescribe $\u = \mathbf{0}$ on $\partial\Omega$.

The material parameters are as follows: the kinematic viscosity is $\nu = 0.1$cm$^2$/s,
the fluid's and solid's densities are $\rho_f = 100$g/cm$^3$. The solid's Lam\'e parameter is $\mu = 1.0\times 10^6$dyne/cm$^2$ and Poisson's ratio is $0.4$.

\begin{figure}[!ht]
    \centering
    \includegraphics[width=14cm]{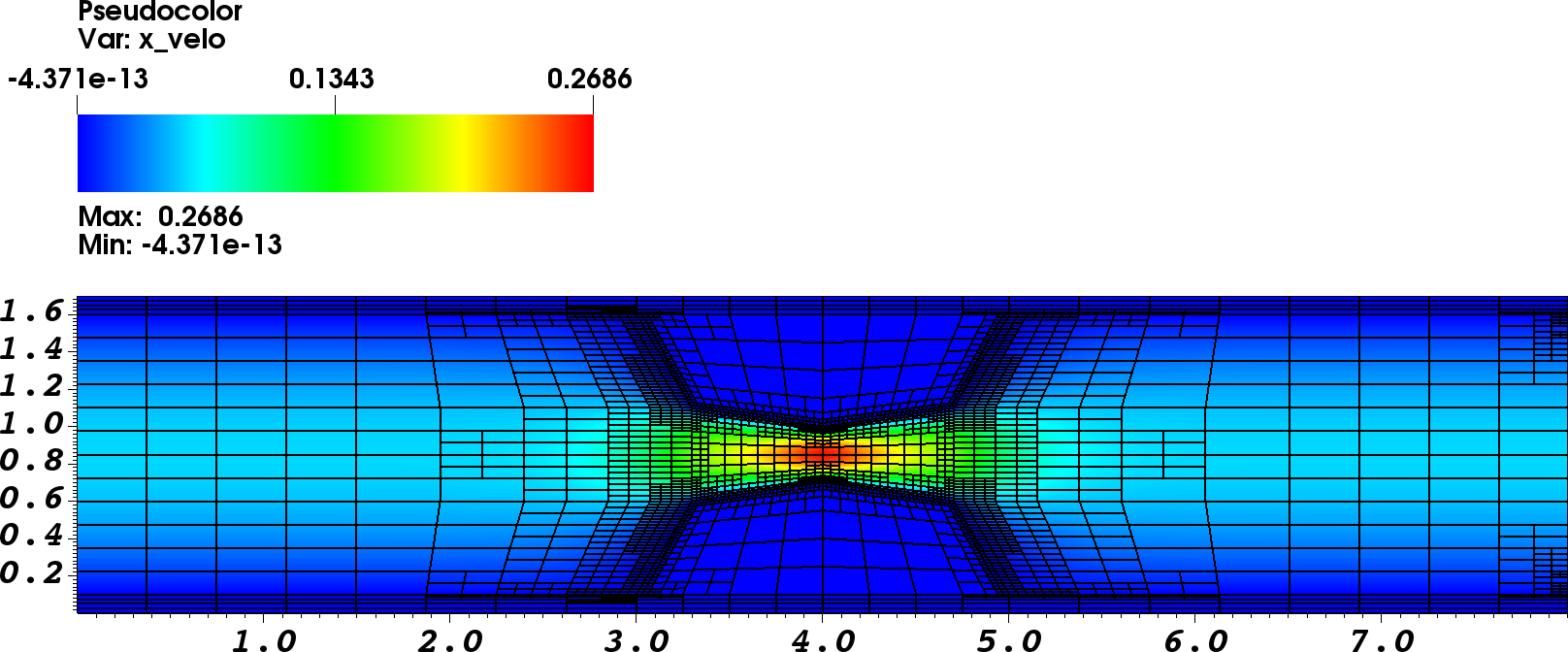}
    \includegraphics[width=14cm]{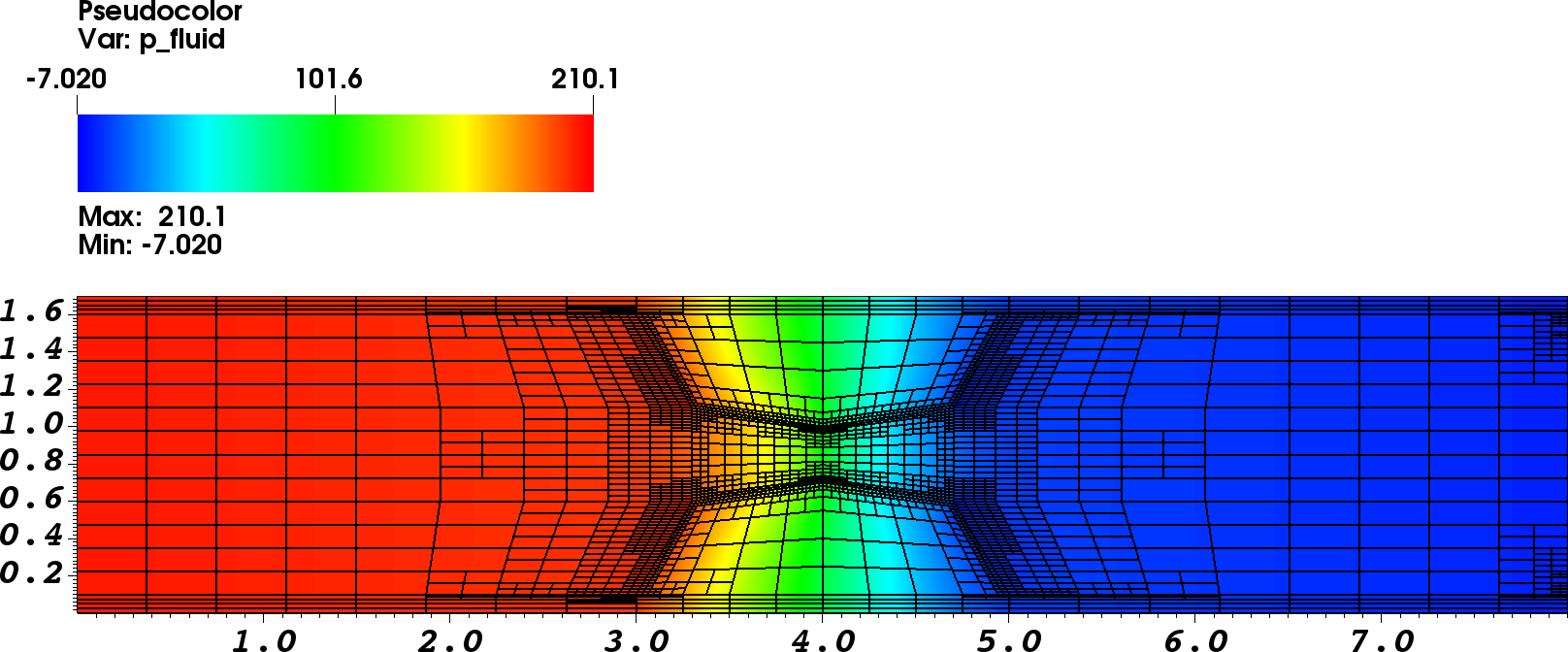}
    \includegraphics[width=14cm]{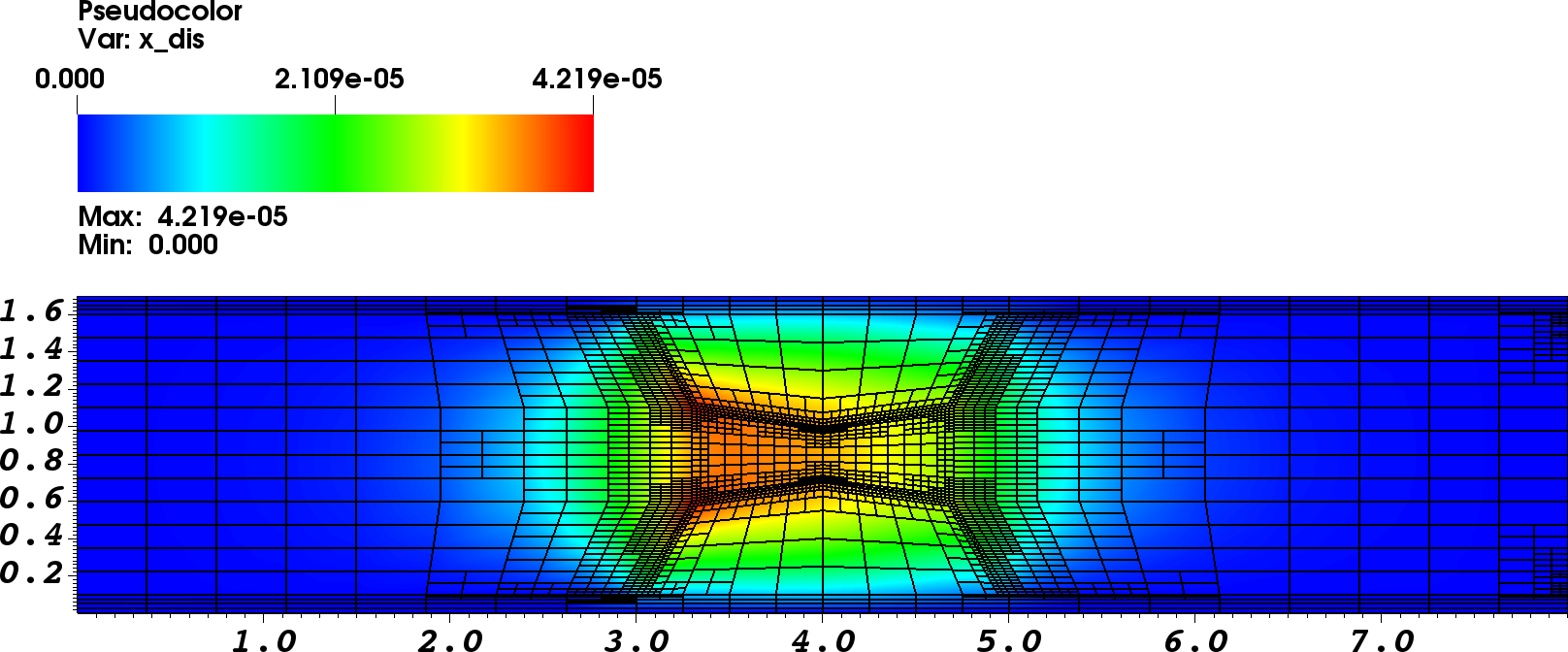}
    \caption{Example multigoal FSI. Primal solutions and adaptively refined mesh. From top to bottom: flow field in $x$ direction, pressure field, displacement field in $x$ direction.}
    \label{fig:fsi_multigoal_b}
\end{figure}

\begin{figure}[!ht]
    \centering
    \includegraphics[width=14cm]{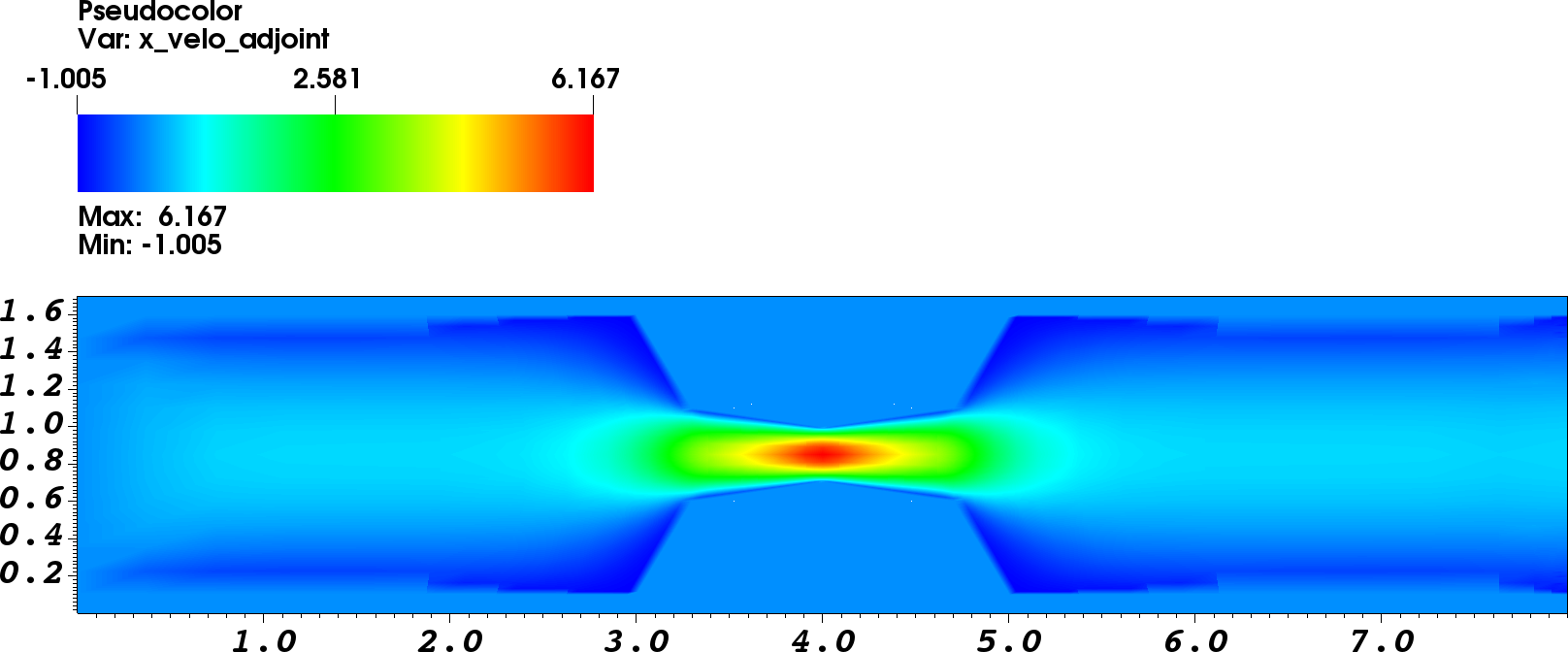}
    \includegraphics[width=14cm]{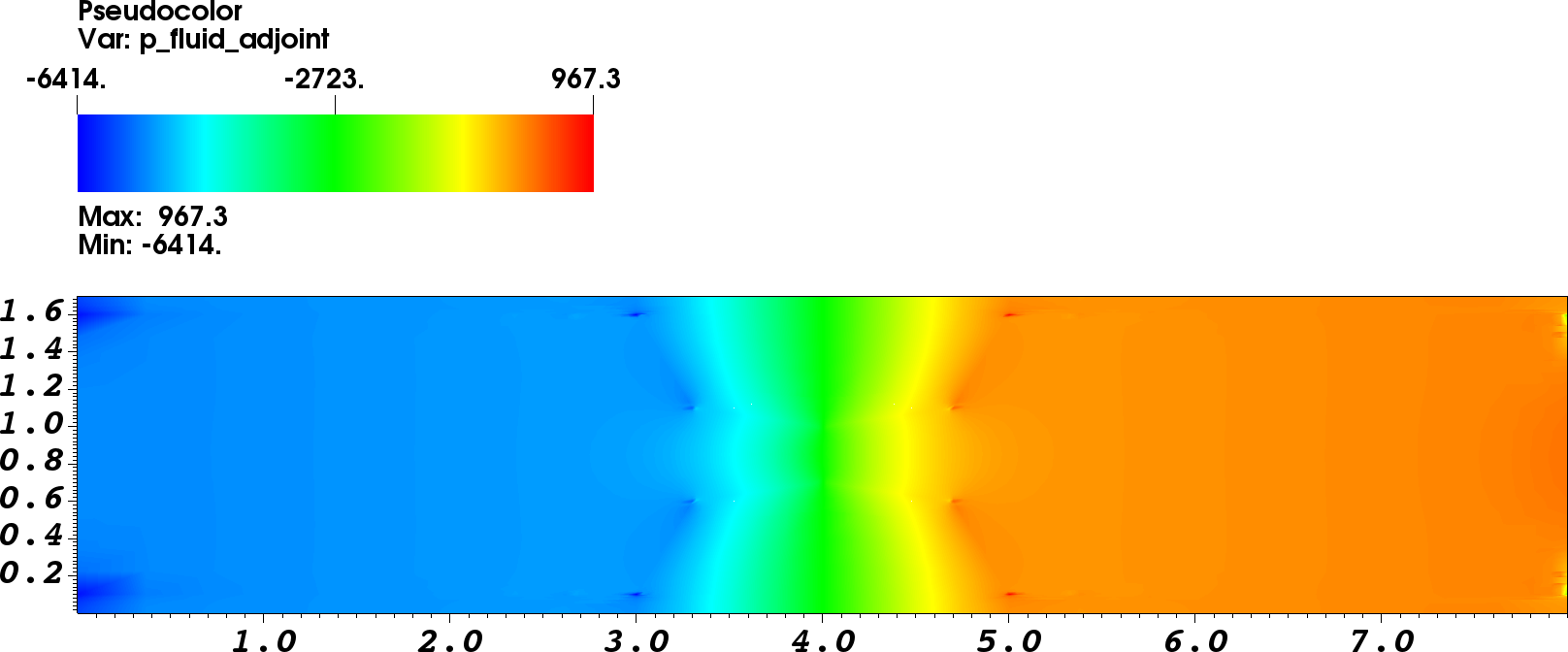}
    \includegraphics[width=14cm]{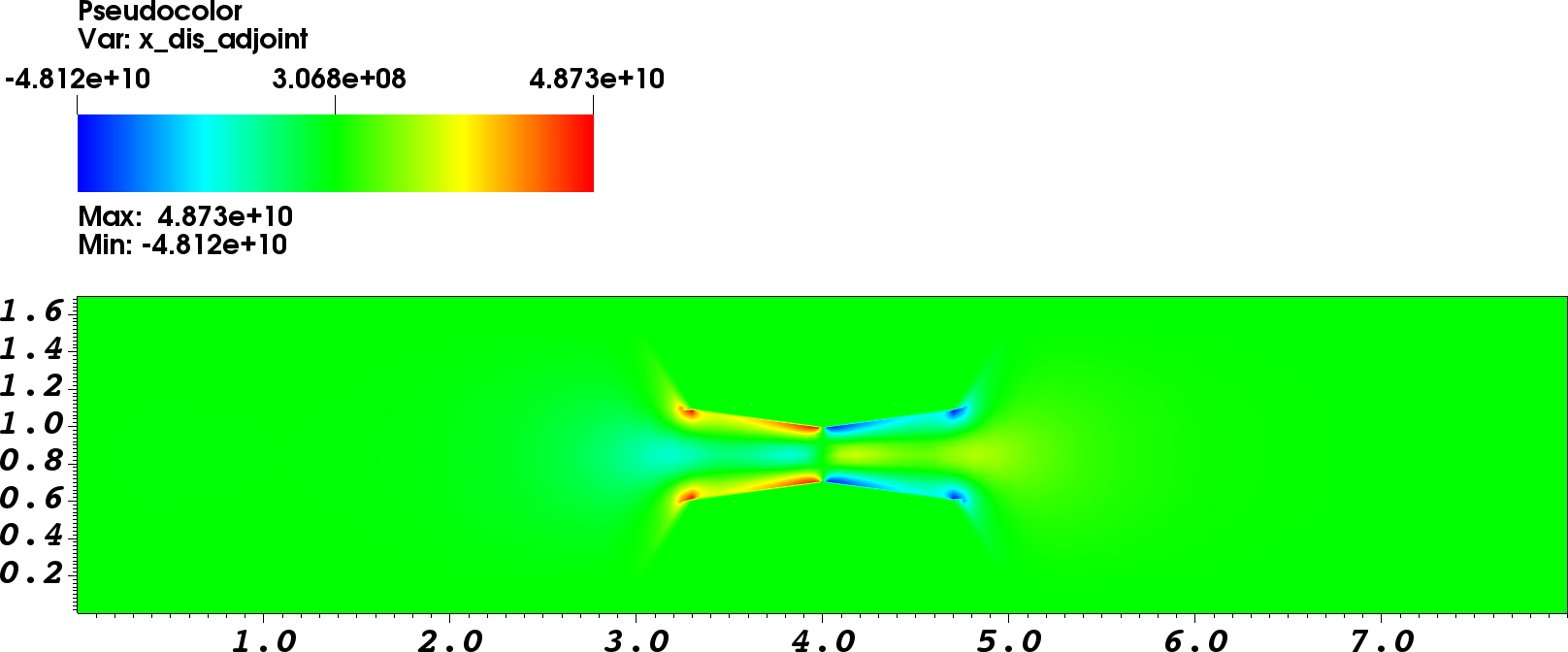}
    \caption{Example multigoal FSI. Adjoint solutions. From top to bottom: flow field in $x$ direction, pressure field, displacement field in $x$ direction.}
    \label{fig:fsi_multigoal_c}
\end{figure}

We choose two goal functionals $J_1(\U)$ and $J_2(\U)$ \franz{(with the notation $\U=(\vv,\u,p)$)}: stress in $x$ direction over the FSI interface, i.e., $J_1(\U) = \int_{\Gamma_i} J\hat\bsigma_f F^{-T}\cdot \mathbf{n} e_1\, ds$, with $e_1 = (1,0)^T$, and the point value $J_2(\u):= \bu_x = (6.5,0.1)$. We notice that the latter has a significant lower value than the first goal functional. The combined goal functional is $J_c(\U) := \omega_1 J_1(\U) + \omega_2 J_2(\u)$,
with $\omega_1 = \omega_2 = 1$.
Graphical illustrations of the adaptive meshes, primal solutions and adjoint solutions are provided in Figure~\ref{fig:fsi_multigoal_b} and Figure~\ref{fig:fsi_multigoal_c}.
The effectivity indices are quite satisfactory 
as shown in Table~\ref{tab:fsi_multigoal_d}.

\begin{table}[ht!]
    \centering
    \begin{tabular}{|r|l|l|l|l|l|}
    \hline
Dofs  &  Exact err   &    Est err   &      Est ind   &      Eff   &          Ind\\ \hline
$8561 $ &  $1.97e+00$    &   $1.37e+00$  &  $2.61e+00$  &   $6.93e-01$  &   $1.32e+00$\\
$15915$ &  $8.86e-01$    &   $6.52e-01$  &  $1.60e+00$  &   $7.36e-01$  &   $1.81e+00$\\
$31361$ &  $4.07e-01$    &   $3.63e-01$  &  $1.06e+00$  &   $8.93e-01$  &   $2.61e+00$\\
$55377$ &  $1.51e-01$    &   $2.22e-01$  &  $6.83e-01$  &   $1.47e+00$  &   $4.51e+00$\\ \hline
    \end{tabular}
    \caption{Example multigoal FSI. Degrees of freedom (DoFs) on the adaptively refined meshes. The exact error reduction (reference value numerically computed) $J_c(\U) - J_c(\U_h)$ with $\U=(\vv,\u,p)$, the estimated error reduction $\eta_h$, estimated indicator reduction, effectivity index, and indicator index. The indicator index and its estimator were introduced in~\cite{RiWi15_dwr}; here all local error indicators $\eta_K$ are taken as their absolute values in the sum to obtain the global estimator $\eta_h$.}
    \label{tab:fsi_multigoal_d}
\end{table}

\clearpage

\section{\revision{Discussion}}
\label{sec:discussion}


In computational biomechanics, errors can significantly impact the reliability of numerical simulations. These errors can be categorized into three main types: modeling errors, discretization errors, and numerical errors. Understanding and controlling these errors is crucial for developing clinically relevant and accurate biomechanical models. 
We provided here a methodology to control and reduce the discretization error and illustrated it. Let us however briefly discuss about the other sources of errors. 


\subsection{Modeling errors}
Modeling errors arise from simplifications and assumptions made during the formulation of a mathematical model. These errors can be introduced at various stages, including the definition of geometry, material properties, boundary conditions, and constitutive laws.
%
%
%
%
%
%
Modeling error can be defined as the discrepancy between the system output and the numerical prediction, assuming that other sources of error (such as discretization and numerical errors) are negligible. However, defining modeling error is inherently challenging due to the absence of a perfect model; all models are, to some extent, approximations of reality.
The sources of modeling errors can be broadly categorized as follows:
\begin{itemize}
    \item \emph{Dimensional Reduction}: Simplifications such as representing 3D structures as beams or shells, or assuming plane strain conditions.
    \item \emph{Geometric Simplification}: Approximations in the geometry, such as idealizing complex shapes or ignoring fine details.
    \item \emph{Boundary Conditions}: Simplifications in boundary conditions, for example, using Dirichlet conditions instead of more realistic Signorini conditions.
    \item \emph{Simplified PDEs}: Assumptions such as small strain elasticity or simplified fluid dynamics equations.
    \item \emph{Constitutive Laws}: Simplifications in the material behavior, such as using 
    isotropic hyperelastic laws instead of anisotropic ones.
    \item \emph{Parameter Calibration}: Errors arising from inaccurate or oversimplified parameter estimation.
    \item \emph{Neglecting Multiphysics and Multiscale Effects}: Ignoring interactions between different physical phenomena or scale-dependent behaviors.
\end{itemize}

\franz{Despite most of the effort in the computational mathematics community has been devoted to the discretization error, some works address errors due to approximations in the model. Let us mention for instance some well-known or recent works (the list is far from being exhaustive):
\begin{itemize}
    \item The general methodologies described in
    \cite{braack2003posteriori,oden2002estimation} based on dual weighted residuals,
    or more recently in \cite{Repin2024}, based upon a very general duality theory and functional identities.
    \item Still based on the dual weighted residual technique, methodologies to control the error and adapt the discretization in the context of parameter calibration, see 
\cite{becker2005,becker2022}.
\item Adaptive refinement for defeaturing problems~\cite{buffa2025}.
\item Estimation of modelling errors in the context of neutronics applications~\cite{taumhas2023}.
\end{itemize}}

\franz{In situations where }
one can assume a hierarchy of models, ranging from the most complex to the simplest. The modeling error associated with a coarse model $C$ can be quantified as $|J_C - J_F|$, where $J$ is a quantity of interest, and $F$ represents the fine ("perfect") model (which, in practice, is often the most complex available model).

In this context, we can sketch, following 
\cite{braack2003posteriori}
(see also \cite{oden2002estimation}) how the Dual Weighted Residual method can be adapted and useful.
%
Consider a fine model of the form: find $u_F \in V$ such that
\[
a(u_F, v) + a_\epsilon(u_F, v) = L(v), \qquad \forall v \in V,
\]
where $V$ is a vector space, $a(u_F, v)$ represents the main (coarse) part of the model, and $a_\epsilon(u_F, v)$ represents the complex part containing detailed physics. 
The source term is denoted by $L(v)$.
The quantity of interest is still given by a linear continuous functional $J(\cdot)$.
If the computation is performed using the coarse model: find $u_C \in V$ solution to
\[
a(u_C, v) = L(v), \qquad \forall v \in V,
\]
the modeling error is defined as:
\[
e_u := u_C - u_F \in V.
\]
This error satisfies the perturbed Galerkin orthogonality:
\[
a(e_u, v) = a(u_C, v) - a(u_F, v) = L(v) - L(v) - a_\epsilon(u_F, v) = -a_\epsilon(u_F, v).
\]
To analyze the error in the quantity of interest, we introduce the fine dual problem: find $z_F \in V$ solution to
\[
a(\phi, z_F) + a_\epsilon(\phi, z_F) = J(\phi), \qquad \forall \phi \in V.
\]
Using this dual problem with $\phi = e_u$ and the perturbed Galerkin orthogonality, we obtain:
\[
J(e_u) = a(e_u, z_F) + a_\epsilon(e_u, z_F) = -a_\epsilon(u_F, z_F) + a_\epsilon(u_C - u_F, z_F) = -a_\epsilon(u_C, z_F).
\]
In practice, a coarse dual problem is often used for increased efficiency: find $z_C \in V$ solution to
\[
a(\phi, z_C) = J(\phi), \qquad \forall \phi \in V.
\]
Starting from the error identity and splitting the right-hand side, we have:
\[
J(e_u) = -a_\epsilon(u_C, z_C) - a_\epsilon(u_C, z_F - z_C).
\]
A simple argument allows to show (see 
\cite{braack2003posteriori}) that the second term can be neglected in a first approximation:
\[
J(e_u) \simeq -a_\epsilon(u_C, z_C).
\]
This allows to use the complex model as a residual for the evaluation of the model error, with the primal and dual solutions computed solely with the coarse model.
As a result, this approach provides a systematic framework for quantifying and controlling modeling errors. 
In \cite{braack2003posteriori} more details can be found, particularly about how to treat nonlinear problems, and how to combine this method with Galerkin approximations.

\subsection{Discretization errors}

Discretization errors occur due to the approximation of continuous problems with discrete numerical methods, such as finite element methods (FEM). A posteriori error estimation acts mostly on this source of errors.
We recall that, for finite element or related (variational) approximations, the discretization error comes mostly from
\begin{enumerate}
\item \emph{the choice of the discrete variational formulation}, particularly if constraints such as incompressibility are treated in a mixed fashion (i.e., introducing extra variables), and if stabilizing or regularizing terms are present. Particularly these terms need to be taken into account into the a posteriori error estimator, and the choice of the stabilization terms or regularizing terms can have an influence on the discretization error. \franz{See, e.g.,~\cite{ainsworth2013} for a methodology to adapt the stabilization parameter in the context of Stokes equations}.

A difficulty can arise for some models (quasi-incompressible elasticity, some beams or plate models, etc) when phenomena such as \emph{locking} happens: the discretization error remains large while the mesh is really refined. Choosing the appropriate discretization technique for a specific model always remain of the utmost importance.

\item \emph{the polynomial order of finite elements}: traditional a priori error analysis show that higher order finite elements yield higher order methods, with much better accuracy. For instance, in a context of uniform refinement, quadratic finite elements perform better than linear finite elements, at the same price in terms of degrees of freedom. However, this is true for instance only if the continuous solution is regular enough. In the context of adaptive mesh refinement, linear finite elements can remain a competitive choice. 

Moreover, Lagrange finite elements perform poorly for order three and beyond, and in this case techniques such as Bernstein-B\'ezier or Isogeometric Analysis can be preferred, \franz{see, e.g.,~\cite{ainsworth2011,atroshchenko2018,hughes2005}}.

\item \emph{the mesh resolution:} as illustrated in this chapter, a coarse mesh can lead to significant discretization errors, while a fine uniform mesh may result in high computational costs. Goal oriented a posteriori error estimation allows to reduce the discretization error without global overrefinement of the mesh.
\end{enumerate}

%

\subsection{Numerical errors}

Numerical errors arise from the practical implementation of numerical algorithms, including linearization errors, numerical integration errors, iterative solver errors, and round-off errors. In general these errors are assumed to be of small magnitude in comparison to model and discretization errors in most of situations. However this point is not always thoroughly studied, and it can happen situations where the numerical errors can be significant, for instance for nonstationnary simulations with large time scales, with accumulation of round-off errors, or when iterative methods are stopped before convergence to save computational time. Many recent works in the adaptive finite element community have been devoted to this subject \franz{(see for instance~\cite{Papez2024} and references therein)}, but not in the context of computational biomechanics to the best of our knowledge.
%
%

\subsection{Balancing the errors}

The concept of \emph{balance of errors} has become the key paradigm in adaptive finite elements. It is also 
highly relevant in computational biomechanics.
Indeed, a classical convergence study, following the canonical validation/verification paradigm, allows, hopefully, to ensure in practice that discretization errors are negligible, but at the price of a finite element mesh with too much degrees of freedom and a high computational cost. With goal-oriented adaptivity driven by a posteriori error estimation, an optimized mesh can be obtained in a few iterations, that allows to ensure that the magnitude of the discretization error is below some threshold. This threshold can be fixed, for instance, in order to balance the modelling errors.



\section{Perspectives}
\label{sec:Conclusions}

\twick{
The goal-oriented error estimators based on the dual-weighted residual method from this work can be further extended to estimate model errors~\cite{oden2002estimation,braack2003posteriori}, 
balancing discretization and iteration errors, see, e.g.,~\cite{RaVih13b} (single goals) or \cite{EnLaWi18}(multiple goals), up to 
adaptive multiscale predictive modeling~\cite{Od18}. 
\franz{The dual-weighted residual method can also be extended in the context of parameter calibration: see for instance~\cite{becker2005,becker2022}.}
Another class of goal-oriented methods is concerned with theoretical investigations (see also the axioms of adaptivity~\cite{Carstensen20141195}) such as recent studies
towards optimal cost complexities~\cite{BeBrInMePr23,BringmannBrunnerPraetoriusStreitberger+2024}.
For time-dependent biomechanics problems, space-time modeling allows for Galerkin finite element discretizations both in space and time. Consequently, goal-oriented error estimates can be extended to the temporal domain that allow for adaptivity in space and time simultaneously. So far, space-time dual-weighted residual methods have been developed 
for heat and combustion problems~\cite{SchVe08,ThiWi24}, incompressible Navier-Stokes equations~\cite{BeRa12,RoThiKoeWi24}, and the p-Laplacian~\cite{ENDTMAYER2024286}. Using space-time modeling, but estimating the temporal error only, was developed for parabolic problems and for the incompressible Navier-Stokes equations in~\cite{MeiRi14,MeiRi15}, and for fluid-structure interaction in~\cite{FaiWi18}. An overview is given in~\cite{wick2024}.}

\appendix

\section{\franz{Incompressible hyperelasticity}}

\franz{In this appendix, we demonstrate how the methodology can be extended to treat (true) incompressible hyperelasticity, see as well \cite{bui:hal-04208610}.}

\subsection{\franz{Case of an incompressible material}}

\franz{To account for incompressibility, we can make the following changes in the setting of Section~\ref{sec:Methods}.} \franz{In the incompressible case, the virtual work of internal forces reads:} 
\begin{equation*} 
A(\bu,\bpp; \bv,\bq) := \int_{\Omega} \bPi(\bu,\bpp): \nabla_X \bv \, ~\dd\x
+ \int_{\Omega} (1-\mathrm{det}(\bC))\bq \, ~\dd\x, \qquad
\end{equation*}
where $\bC := \bF^T \cdot \bF$ denotes the right Cauchy-Green tensor,
and where $\bu$ and $\bv$ are admissible displacements and $\bpp$ and $\bq$ are admissible pressures. The condition $\mathrm{det}(\bC)=1$ is the incompressibility condition. 
The hyperelastic problem in weak form reads
\begin{equation}
\left\{
\begin{array}{l}
\text{Find a displacement }\bu, 
\text{with } 
\bu = \bu_D \text{ on } \Gamma_D
\text{ and a pressure }\bpp 
\text{ such that } \\\noalign{\smallskip}
A(\bu,\bpp; \bv,\bq) = L(\bv), \, \forall (\bv,\bq), 
\, \bv = \mathbf{0} \text{ on } \Gamma_D.
\end{array}
\right.
\label{eq:primal_weak2_incompress}
\end{equation}

Let $\mathcal{K}_h$ be a mesh of the domain $\Omega$.
Let us denote by $\bV_h \times \bQ_h$ the finite element pair that makes use of the lowest-order Taylor-Hood finite elements on $\mathcal{K}_h$ (continuous piecewise polynomials of order $2$ for the displacement and of order $1$ for the pressure). The finite element method to solve our hyperelastic problem reads 
%
\begin{equation}
\left\{
\begin{array}{l}
\text{Find a displacement }\bu_h \in \bV_h \text{, with } \bu_h = \bu_D^h \text{ on } \Gamma_D \text{ and a pressure }\bpp_h\in\bQ_h\text{ such that } \\\noalign{\smallskip}
A(\bu_h,\bpp_h; \bv_h,\bq_h) = L(\bv_h), \, \forall (\bv_h,\bq_h) \in \bV_h^0\times\bQ_h,
\end{array}
\right.
\label{eq:FEM_primal_weak_incompress}
\end{equation}
where $\bV_h^0$ is composed by the functions of $\bV_h$ vanishynig on $\Gamma_D$ and where $u_D^h$ is a finite element approximation of $u_D$, obtained for instance by Lagrange interpolation or by projection.

\subsection{Dual problem for computing the weights}
\label{sec_dual_problem}
\franz{For incompressible hyperelasticity, the nonlinear quantity of interest $Q$ reads} 
\begin{equation*}
|J(\bu,\bpp)-J(\bu_h,\bpp_h)|
\end{equation*}
and shall be minimized with the PDE as a constraint:
\begin{equation*}
    \min |J(\bu,\bpp)-J(\bu_h,\bpp_h)|
\quad\text{s.t. } A(\cdot,\cdot; \cdot,\cdot) = L(\cdot)
\end{equation*}
\franz{This 
implies to compute the following dual problem:
\begin{equation}
\left\{
\begin{array}{l}
\text{Find } ({\bz}_h,\bw_h) \in {\bV}_h^0  \times \bQ_h \text{ such that } \\\noalign{\smallskip}
(A')^*(\bu_h,\bpp_h| {\bz}_h,  {\bw}_h;{\bv}_h, {\bq}_h) = Q'(\bu_h,\bpp_h|{\bv}_h,\bq_h) \quad \forall ({\bv}_h,\bq_h) \in {\bV}_h^0  \times \bQ_h,
\end{array}
\right.
\label{eq:discrete_dual_problemIH}
\end{equation}
where $A'$ and $Q'$ denote the Fr\'echet derivative of $A$ and $Q$, respectively,
and $(A')^*$ is the adjoint form of $A'$.}

\subsection{The representation formula of Becker and Rannacher}

We introduce $r(\bu_h,\bpp_h;\bv,\bq)$ the residual of Problem 
\eqref{eq:primal_weak2_incompress}
as
\begin{equation}
r(\bu_h,\bpp_h;\bv,\bq) = L(\bv,\bq) - A(\bu_h,\bpp_h;\bv,\bq) \qquad \forall (\bv,\bq) \in \bV\times\bQ.
\end{equation}
This, roughly speaking, quantifies how well the hyperelasticity equations are approximated (it should tend to zero if the mesh is uniformly refined).
Thanks to the dual system \eqref{eq:discrete_dual_problemIH},
we obtain expression of the error on $Q$ 
as the best approximation term involving the residual and the (exact) dual solution
(see \cite[Proposition 2.3]{becker-rannacher-2001}): 
\begin{equation}
J(\bu,\bpp) - J(\bu_h,\bpp_h) = \min_{(\bv_h,\bq_h) \in \bV_h\times\bQ_h} r(\bu_h,\bpp_h; \bz - \bv_h, \bw - \bq_h) + R_m \label{eq:representation_DWRIH}
\end{equation}
where 
$R_m$ is the high-order remainder related to the error caused by the linearization of the nonlinear problem (the precise expression of which can be found in \cite{becker-rannacher-2001}). 
%

Proceeding as usual in \emph{a posteriori} error estimation, \emph{i.e.,}
after performing integration by parts on the residual $r$, we localize the different contributions to the goal-oriented error as follows:
\begin{equation}
\lvert J(\bu,\bpp) - J(\bu_h,\bpp_h) \rvert \leq \sum_{K \in \mathcal{K}_h} \eta_K{((\bu_K,\bpp_K),(\bz_K,\bw_K))} + H.O.T.
\end{equation}
In the above expression, $K$ denotes any cell of the mesh $\mathcal{K}_h$, and expressions such as $\bu_K$ denote the local restriction of the finite element variable $\bu_h$ of the cell $K$. Moreover, $H.O.T.$ denotes high order terms, that are not considered in the implementation. 

\subsection{Adaptive mesh refinement}

Using the error estimate on $J$, we implement a standard procedure for mesh refinement. 
As described in \cref{algo:adaptivity}, we start with an initial mesh called $mesh_i$, and by providing a guessed solution $\bu_i^{(0)}$, the nonlinear primal problem can be solved using Newton's method (see \cref{algo:newton_method}). 
Once accepting $\bu_i$ as the solution of the discrete primal problem, solving the discrete dual problem (see \cref{algo:dual_problem}) 
provides the dual solution $\bz_i$ on the current mesh $mesh_i$. The estimator $\eta_K$ is then computed by using the primal and dual solutions $\bu_i$ and $\bz_i$, respectively. From the estimator, different strategies can be used to mark the elements whose \emph{error} is high. In this paper, we use the D\"orfler marking strategy \cite{dorfler1996} (see Algorithm \ref{algo:Dorfler_making}).
%
%

\subsection{Expression of the estimator and algorithms}
We give below for each cell-wise contribution:
\begin{equation}
\eta_K = \left| \int_K \bR_{\bu} \cdot (\widehat\bz_h-i_h(\widehat\bz_h)) \mathrm{d}\Omega + \int_K \bR_{\bpp} \cdot (\widehat\bw_h-i_h(\widehat\bw_h)) \mathrm{d}\Omega +  \int_{\partial K} \bJ \cdot (\widehat\bz_h-i_h(\widehat\bz_h)) \mathrm{d}\gamma   \right|
\end{equation}
with, 
the interior residual
\begin{equation*}
\bR_{\bu} = \bB + 
\mathrm{div} \, \bPi(\bu_h)
\mbox{ and }\bR_{\bpp}=\det(\boldsymbol{C})-1
\end{equation*}
and the stress jump 
{
\begin{equation*}
\bJ = \begin{cases}
   - \frac{1}{2}   [[\bPi(\bu_h)]]
       & \text{ if } F \not \subset \Gamma, \\
      \bT -  \bPi(\u_h)  \cdot n_F                                 & \text{ if } F \subset \Gamma_N, \\
      \bzero                                                                     & \text{ if } F \subset \Gamma_D,
      \end{cases}
\end{equation*}
for each facet $F$, where $n_F$ is the exterior normal to the facet $F$ of $\Gamma_N$. }
\begin{algorithm}
\caption{Algorithm for mesh refinement}\label{algo:adaptivity}
\begin{algorithmic}
\State Select an initial triangulation $mesh_i$ of the domain $\Omega$
\State Guessed solution $(\bu_i^{(0)},\bpp_i^{(0)})$
\While{$ \sum_K \eta_K > \epsilon$}
\State $F(\bu_i,\bpp_i;\bv,\bq) \gets A(\bu_i,\bpp_i;\bv,\bq) - L(\bv,\bq)$
\State $\bu_i,\bpp_i \gets$ NewtonMethod $(F(\bu_i,\bpp_i;\bv,\bq),(\bu_i^{(0)},\bpp_i^{(0)}))$
\Comment{Problem \eqref{eq:primal_weak2}, see Algo \ref{algo:newton_method}} 
\State $\bz_i ,\bw_i\gets$ DualProblem $(\bu_i,\bpp_i, Q)$                      \Comment{Problem \eqref{eq:nonlinear_dual_problem}, see Algo \ref{algo:dual_problem}} 
\State $\eta_K \gets$ ComputeEstimator $(\bu_i,\bpp_i,\bz_i,\bw_i)$
\State markedElements $\gets$ DorflerMarking $(\eta_K, \alpha)$ \Comment{See Algo \cref{algo:Dorfler_making}}
\State $mesh_i \gets mesh_i$.refine(markedElements)
\State Compute $\sum_K \eta_K$
\EndWhile
\end{algorithmic}
\end{algorithm}
%
\begin{algorithm}
\caption{Solving a non-linear problem: NewtonMethod $(F(\bu_i^{(0)},\bpp_i^{(0)};\bv, \bpp),\bu_i^{(0)},\bpp_i^{(0)} )$}\label{algo:newton_method}
\begin{algorithmic}
\State $(\bu_k,\bpp_k) = (\bu_i^{(0)},\bpp_i^{(0)}) $
\While{$\lvert (\delta \bu,\delta\bpp) \rvert > \epsilon$}
 \State $ F'(\bu_k ,\bpp_k; \delta \bu ,\delta\bpp, \bv,\bq)  = -F(\bu_k, \bpp_k; \bv , \bq) $    \Comment{Solve for $(\delta \bu,\delta\bpp)$} 
 \State $ (\bu_{k+1},\bpp_{k+1}) \gets (\bu_k + \delta \bu,\bpp_k + \delta \bpp) $    \Comment{Update the solution}
 \State Compute $\lvert (\delta \bu,\delta\bpp) \rvert$
\EndWhile
\end{algorithmic}
\end{algorithm}

\begin{algorithm}
\caption{Solving the dual problem: DualProblem $(\bu_i,\bpp_i,Q)$}\label{algo:dual_problem}
\begin{algorithmic}
\State Compute $A'(\bu_i,\bpp_i| \bz_i,\bw_i;\bv,\bq) $
\State Compute $Q'(\bu_i,\bpp_i; \bv,\bq)$
\State $(\bz_i,\bw_i) \gets$ solve $(A'(\bu_i,\bpp_i| \bz_i,\bw;\bv,\bq) = Q'(\bu_i,\bpp_i; \bv,\bq) )$ \Comment{Solve the linear system}
\end{algorithmic}
\end{algorithm}
\begin{algorithm}
\caption{Mark elements after D\"orfler strategy by providing a element-wise estimator $\eta_K = [\eta_{K_1}, \eta_{K_2}, \dots \eta_{K_N}]$, and $0 < \alpha < 1$ a parameter which characterises the marking rate: the smaller the value of $\alpha$ is, the fewer the number of elements will be marked: DorflerMarking $(\eta_K, \alpha)$ }\label{algo:Dorfler_making}
\begin{algorithmic}
\State Sort the elements $K_i$ after descending order of the corresponding estimator $\eta_{K_i}$
\State Mark the first $M$ elements such that 
\[\mathrm{markedElements} \gets \min \left \{ M \in N \:\middle |\: \sum\limits_{i=1}^{M}\eta_{K_i}\geq\alpha\sum_{i=1}^N \eta_{K_i} \right \}. \]
\end{algorithmic}
\end{algorithm}

\clearpage

\bibliographystyle{abbrv}
\bibliography{biblio.bib}

@incollection{Papez2024,
  author    = {Jan Papež},
  title     = {Algebraic error in numerical PDEs and its estimation},
  booktitle = {Advances in Applied Mechanics},
  volume    = {58},
  year      = {2024},
  pages     = {377-427},
  doi       = {10.1016/bs.aams.2024.04.002},
  publisher = {Elsevier},
}

@Article{Margenberg2021,
  author        = {N. Margenberg and T. Richter},
  title         = {Parallel time-stepping for fluid-structure interactions},
  doi           = {10.1051/mmnp/2021005},
  eprint        = {1907.01252},
  pages         = {20},
  volume        = {16},
  journal       = {Mathematical Modelling of Natural Phenomena},
  year          = {2021},
}

@article{ainsworth2011,
 author = {Ainsworth, Mark and Andriamaro, Gaelle and Davydov, Oleg},
 title = {Bernstein-{B{\'e}zier} finite elements of arbitrary order and optimal assembly procedures},
 journal = {SIAM Journal on Scientific Computing},
 ajournal = {SIAM J. Sci. Comput.},
 issn = {1064-8275},
 volume = {33},
 number = {6},
 pages = {3087--3109},
 year = {2011},
 language = {English},
 doi = {10.1137/11082539X},
 zbMATH = {6008955},
 Zbl = {1237.65120}
}

@article{atroshchenko2018,
 author = {Atroshchenko, Elena and Tomar, Satyendra and Xu, Gang and Bordas, St{\'e}phane P. A.},
 title = {Weakening the tight coupling between geometry and simulation in isogeometric analysis: from sub- and super-geometric analysis to {Geometry}-{Independent} {Field} {approximaTion} ({GIFT})},
 journal = {International Journal for Numerical Methods in Engineering},
 ajournal = {Int. J. Numer. Methods Eng.},
 issn = {0029-5981},
 volume = {114},
 number = {10},
 pages = {1131--1159},
 year = {2018},
 language = {English},
 doi = {10.1002/nme.5778},
 zbMATH = {7878352},
 Zbl = {1548.74760}
}

@article{hughes2005,
 author = {Hughes, T. J. R. and Cottrell, J. A. and Bazilevs, Y.},
 title = {Isogeometric analysis: {CAD}, finite elements, {NURBS}, exact geometry and mesh refinement},
 journal = {Computer Methods in Applied Mechanics and Engineering},
 ajournal = {Comput. Methods Appl. Mech. Eng.},
 issn = {0045-7825},
 volume = {194},
 number = {39-41},
 pages = {4135--4195},
 year = {2005},
 language = {English},
 doi = {10.1016/j.cma.2004.10.008},
 zbMATH = {5142346},
 Zbl = {1151.74419}
}

@article{ainsworth2013,
 author = {Ainsworth, Mark and Allendes, Alejandro and Barrenechea, Gabriel R. and Rankin, Richard},
 title = {On the adaptive selection of the parameter in stabilized finite element approximations},
 journal = {SIAM Journal on Numerical Analysis},
 ajournal = {SIAM J. Numer. Anal.},
 issn = {0036-1429},
 volume = {51},
 number = {3},
 pages = {1585--1609},
 year = {2013},
 language = {English},
 doi = {10.1137/110837796},
 url = {strathprints.strath.ac.uk/44170/1/110837796.pdf},
 zbMATH = {6203689},
 Zbl = {1276.65076}
}

@article{taumhas2023,
 author = {Taumhas, Yonah Conjungo and Labeurthre, David and Madiot, Francois and Mula, Olga and Taddei, Tommaso},
 title = {Impact of physical model error on state estimation for neutronics applications},
 fjournal = {European Series in Applied and Industrial Mathematics (ESAIM): Proceedings and Surveys},
 journal = {ESAIM, Proc. Surv.},
 issn = {2267-3059},
 volume = {73},
 pages = {158--172},
 year = {2023},
 language = {English},
 doi = {10.1051/proc/202373158},
 zbMATH = {7829997}
}

@misc{buffa2025,
 author = {Buffa, Annalisa and Grappein, Denise and V{\'a}zquez, Rafael},
 title = {Adaptive refinement in defeaturing problems via an equilibrated flux a posteriori error estimator},
 year = {2025},
 howpublished = {Preprint, {arXiv}:2503.19784 [math.{NA}] (2025)},
 url = {https://arxiv.org/abs/2503.19784},
 nnote = {arXiv:2503.19784},
 arXiv = {arXiv:2503.19784}
}

@incollection{Repin2024,
  author    = {Sergey I. Repin},
  title     = {A posteriori error identities and estimates of modelling errors},
  booktitle = {Advances in Applied Mechanics},
  volume    = {58},
  year      = {2024},
  pages     = {245-293},
  doi       = {10.1016/bs.aams.2024.03.006},
  publisher = {Elsevier},
}

@article{fernandez2011,
 author = {Fern{\'a}ndez, Miguel {\'A}ngel},
 title = {Coupling schemes for incompressible fluid-structure interaction: implicit, semi-implicit and explicit},
 fjournal = {S\(\vec{\text{e}}\)MA Journal},
 journal = {S\(\vec{\text{e}}\)MA J.},
 issn = {2254-3902},
 volume = {55},
 pages = {59--108},
 year = {2011},
 language = {English},
 doi = {10.1007/BF03322593},
 zbMATH = {6054322},
 Zbl = {1242.76201}
}

@Book{Richter2017,
  author    = {Richter, Thomas},
  publisher = {Springer International Publishing},
  title     = {Fluid-structure Interactions},
  year      = {2017},
  address   = {Cham},
  isbn      = {9783319639703},
  number    = {118},
  series    = {Lecture Notes in Computational Science and Engineering},
  doi       = {10.1007/978-3-319-63970-3},
  issn      = {2197-7100},
  journal   = {Lecture Notes in Computational Science and Engineering},
  pagetotal = {1436},
  ppn_gvk   = {1724413511},
  subtitle  = {Models, Analysis and Finite Elements},
}

@Article{YangRichterJaegerNeussRadu2015,
  author  = {Y. Yang and T. Richter and W. Jaeger and M. Neuss-Radu},
  title   = {An ALE approach to mechano-chemical processes in fluid-structure interactions},
  doi     = {10.1002/fld.4345},
  number  = {4},
  pages   = {199--220},
  url     = {https://www.math.uni-magdeburg.de/~richter/pdf-files/2016-YangRichterJaegerNeussRadu-IJNMF.pdf},
  volume  = {84},
  ijournal = {I. J. Num. Meth. Fluids},
  journal = {International Journal of Numerical Methods in Fluids},
  year    = {2017},
}

@Article{YangJaegerNeussRaduRichter2015,
  author  = {Y. Yang and W. Jaeger and M. Neuss-Radu and T. Richter},
  title   = {Mathematical modeling and simulation of the evolution of plaques in blood vessels},
  doi     = {10.1007/s00285-015-0934-8},
  pages   = {973--996},
  volume  = {72},
  journal = {Journal of Mathematical Biology},
  year    = {2016},
}

@Book{Formaggia2009,
  editor    = {L. Formaggia and A. Quarteroni and A. Veneziani},
  publisher = {Springer Milan},
  title     = {Cardiovascular Mathematics},
  year      = {2009},
  isbn      = {9788847011526},
  doi       = {10.1007/978-88-470-1152-6},
}

@Article{Failer2021,
  author    = {Failer, Lukas and Minakowski, Piotr and Richter, Thomas},
  journal   = {Vietnam Journal of Mathematics},
  title     = {On the Impact of Fluid Structure Interaction in Blood Flow Simulations: Stenotic Coronary Artery Benchmark},
  year      = {2021},
  issn      = {2305-2228},
  month     = jan,
  number    = {1},
  pages     = {169--187},
  volume    = {49},
  doi       = {10.1007/s10013-020-00456-6},
  publisher = {Springer Science and Business Media LLC},
}

@Article{Faizal2019,
  author    = {Faizal, W.M. and Ghazali, N.N.N. and Badruddin, Irfan Anjum and Zainon, M.Z. and Yazid, Aznijar Ahmad and Ali, Mohamad Azlin Bin and Khor, C.Y. and Ibrahim, Norliza Binti and Razi, Roziana M.},
  journal   = {Computer Methods and Programs in Biomedicine},
  title     = {A review of fluid-structure interaction simulation for patients with sleep related breathing disorders with obstructive sleep},
  year      = {2019},
  issn      = {0169-2607},
  month     = oct,
  pages     = {105036},
  volume    = {180},
  doi       = {10.1016/j.cmpb.2019.105036},
  publisher = {Elsevier BV},
}

@Article{Arminio2024,
  author    = {Arminio, Mariachiara and Carbonaro, Dario and Morbiducci, Umberto and Gallo, Diego and Chiastra, Claudio},
  journal   = {Frontiers in Medical Technology},
  title     = {Fluid-structure interaction simulation of mechanical aortic valves: a narrative review exploring its role in total product life cycle},
  year      = {2024},
  issn      = {2673-3129},
  month     = jul,
  volume    = {6},
  doi       = {10.3389/fmedt.2024.1399729},
  publisher = {Frontiers Media SA},
}

@Article{Causin2005,
  author    = {Causin, P. and Gerbeau, J.F. and Nobile, F.},
  journal   = {Computer Methods in Applied Mechanics and Engineering},
  title     = {Added-mass effect in the design of partitioned algorithms for fluid–structure problems},
  year      = {2005},
  issn      = {0045-7825},
  mmonth     = oct,
  number    = {42–44},
  pages     = {4506--4527},
  volume    = {194},
  doi       = {10.1016/j.cma.2004.12.005},
  publisher = {Elsevier BV},
}

@Article{Heil2008,
  author    = {Heil, Matthias and Hazel, Andrew L. and Boyle, Jonathan},
  journal   = {Computational Mechanics},
  title     = {Solvers for large-displacement fluid–structure interaction problems: segregated versus monolithic approaches},
  year      = {2008},
  issn      = {1432-0924},
  month     = mar,
  number    = {1},
  pages     = {91--101},
  volume    = {43},
  doi       = {10.1007/s00466-008-0270-6},
  publisher = {Springer Science and Business Media LLC},
}

@Article{Bessonov2015,
  author    = {Bessonov, N. and Sequeira, A. and Simakov, S. and Vassilevskii, Yu. and Volpert, V.},
  journal   = {Mathematical Modelling of Natural Phenomena},
  title     = {Methods of Blood Flow Modelling},
  year      = {2015},
  issn      = {1760-6101},
  month     = dec,
  number    = {1},
  pages     = {1--25},
  volume    = {11},
  doi       = {10.1051/mmnp/201611101},
  editor    = {Volpert, V.},
  publisher = {EDP Sciences},
}

@Article{Holzapfel2000,
  author    = {Holzapfel, Gerhard A. and Gasser, Thomas C. and Ogden, Ray W.},
  title = {A New Constitutive Framework for Arterial Wall Mechanics and a Comparative Study of Material Models},
  journal   = {Journal of Elasticity},
  year      = {2000},
  issn      = {0374-3535},
  number    = {1/3},
  pages     = {1--48},
  volume    = {61},
  doi       = {10.1023/a:1010835316564},
  publisher = {Springer Science and Business Media LLC},
}

@Article{Richter2013,
  author    = {Richter, Thomas},
  journal   = {Journal of Computational Physics},
  title     = {{A Fully Eulerian formulation for fluid–structure-interaction problems}},
  year      = {2013},
  issn      = {0021-9991},
  mmonth     = jan,
  pages     = {227--240},
  volume    = {233},
  doi       = {10.1016/j.jcp.2012.08.047},
  publisher = {Elsevier BV},
}

@Article{Wick2013,
  author    = {Wick, Thomas},
  journal   = {Computer Methods in Applied Mechanics and Engineering},
  title     = {{Fully Eulerian fluid–structure interaction for time-dependent problems}},
  year      = {2013},
  issn      = {0045-7825},
  month     = mar,
  pages     = {14--26},
  volume    = {255},
  doi       = {10.1016/j.cma.2012.11.009},
  publisher = {Elsevier BV},
}

@Article{Dunne2006,
  author    = {Dunne, Th.},
  journal   = {International Journal for Numerical Methods in Fluids},
  title     = {{An Eulerian approach to fluid–structure interaction and goal‐oriented mesh adaptation}},
  year      = {2006},
  issn      = {1097-0363},
  mmonth     = mar,
  number    = {9–10},
  pages     = {1017--1039},
  volume    = {51},
  doi       = {10.1002/fld.1205},
  publisher = {Wiley},
}

@Article{Wahl2021,
  author    = {von Wahl, Henry and Richter, Thomas and Frei, Stefan and Hagemeier, Thomas},
  journal   = {Physics of Fluids},
  title     = {Falling balls in a viscous fluid with contact: Comparing numerical simulations with experimental data},
  year      = {2021},
  issn      = {1089-7666},
  month     = mar,
  number    = {3},
  volume    = {33},
  doi       = {10.1063/5.0037971},
  publisher = {AIP Publishing},
}

@Article{Cottet2008,
  author    = {Cottet, Georges-Henri and Maitre, Emmanuel and Milcent, Thomas},
  journal   = {ESAIM: Mathematical Modelling and Numerical Analysis},
  title     = {Eulerian formulation and level set models for incompressible fluid-structure interaction},
  year      = {2008},
  issn      = {1290-3841},
  month     = apr,
  number    = {3},
  pages     = {471--492},
  volume    = {42},
  doi       = {10.1051/m2an:2008013},
  publisher = {EDP Sciences},
}

@Article{Wick2011,
  author    = {Wick, Thomas},
  journal   = {Computers \& Structures},
  title     = {Fluid-structure interactions using different mesh motion techniques},
  year      = {2011},
  issn      = {0045-7949},
  month     = jul,
  number    = {13–14},
  pages     = {1456--1467},
  volume    = {89},
  doi       = {10.1016/j.compstruc.2011.02.019},
  publisher = {Elsevier BV},
}

@InBook{Richter2015,
  author    = {Richter, Thomas and Wick, Thomas},
  pages     = {377--400},
  publisher = {Springer International Publishing},
  title     = {On Time Discretizations of Fluid-Structure Interactions},
  year      = {2015},
  isbn      = {9783319233215},
  booktitle = {Multiple Shooting and Time Domain Decomposition Methods},
  doi       = {10.1007/978-3-319-23321-5_15},
  issn      = {2191-3048},
}

@Article{Wall2008,
  author    = {Wall, Wolfgang A. and Rabczuk, Timon},
  journal   = {International Journal for Numerical Methods in Fluids},
  title     = {Fluid–structure interaction in lower airways of CT‐based lung geometries},
  year      = {2008},
  issn      = {1097-0363},
  month     = apr,
  number    = {5},
  pages     = {653--675},
  volume    = {57},
  doi       = {10.1002/fld.1763},
  publisher = {Wiley},
}

@Article{Cardillo2021,
  author    = {Cardillo, Giulia and Camporeale, Carlo},
  journal   = {Journal of Fluids and Structures},
  title     = {Modeling fluid–structure interactions between cerebro-spinal fluid and the spinal cord},
  year      = {2021},
  issn      = {0889-9746},
  mmonth     = apr,
  pages     = {103251},
  volume    = {102},
  doi       = {10.1016/j.jfluidstructs.2021.103251},
  publisher = {Elsevier BV},
}

@Article{Astorino2009,
  author    = {Astorino, Matteo and Gerbeau, Jean-Frédéric and Pantz, Olivier and Traoré, Karim-Frédéric},
  journal   = {Computer Methods in Applied Mechanics and Engineering},
  title     = {Fluid–structure interaction and multi-body contact: Application to aortic valves},
  year      = {2009},
  issn      = {0045-7825},
  month     = sep,
  number    = {45–46},
  pages     = {3603--3612},
  volume    = {198},
  doi       = {10.1016/j.cma.2008.09.012},
  publisher = {Elsevier BV},
}

@PHDTHESIS{Wick2012,
   author = {T. Wick},
   title = {Adaptive {F}inite {E}lement {S}imulation of 
{F}luid-{S}tructure {I}nteraction 
with {A}pplication to {H}eart-{V}alve {D}ynamics},
   school = {University of Heidelberg},
   year = {2011},
   url = {http://www.ub.uni-heidelberg.de/archiv/12992}
}

@Article{Failer2020,
  author    = {Failer, L. and Richter, T.},
  journal   = {Journal of Scientific Computing},
  title     = {A Parallel {N}ewton Multigrid Framework for Monolithic Fluid-Structure Interactions},
  year      = {2020},
  issn      = {1573-7691},
  mmonth     = jan,
  number    = {2},
  volume    = {82},
  doi       = {10.1007/s10915-019-01113-y},
  publisher = {Springer Science and Business Media LLC},
}

@Article{Richter2015a,
  author    = {Richter, Thomas},
  journal   = {International Journal for Numerical Methods in Engineering},
  title     = {A monolithic geometric multigrid solver for fluid-structure interactions in {ALE} formulation},
  year      = {2015},
  issn      = {0029-5981},
  mmonth     = may,
  number    = {5},
  pages     = {372--390},
  volume    = {104},
  doi       = {10.1002/nme.4943},
  publisher = {Wiley},
}

@Article{Boffi2024,
  author    = {Boffi, Daniele and Credali, Fabio and Gastaldi, Lucia and Scacchi, Simone},
  journal   = {Mathematics and Computers in Simulation},
  title     = {{A parallel solver for fluid–structure interaction problems with Lagrange multiplier}},
  year      = {2024},
  issn      = {0378-4754},
  mmonth     = jun,
  pages     = {406--424},
  volume    = {220},
  doi       = {10.1016/j.matcom.2024.01.027},
  publisher = {Elsevier BV},
}

@Article{Zee2011,
  author    = {van der Zee, K.G. and van Brummelen, E.H. and Akkerman, I. and de Borst, R.},
  journal   = {Computer Methods in Applied Mechanics and Engineering},
  title     = {Goal-oriented error estimation and adaptivity for fluid–structure interaction using exact linearized adjoints},
  year      = {2011},
  issn      = {0045-7825},
  month     = sep,
  number    = {37–40},
  pages     = {2738--2757},
  volume    = {200},
  doi       = {10.1016/j.cma.2010.12.010},
  publisher = {Elsevier BV},
}

@Article{Richter2012,
  author    = {Richter, Th.},
  journal   = {Computer Methods in Applied Mechanics and Engineering},
  title     = {Goal-oriented error estimation for fluid–structure interaction problems},
  year      = {2012},
  issn      = {0045-7825},
  mmonth     = jun,
  pages     = {28--42},
  volume    = {223–224},
  doi       = {10.1016/j.cma.2012.02.014},
  publisher = {Elsevier BV},
}

@Article{Wick2012a,
  author    = {Wick, Thomas},
  journal   = {Archive of Mechanical Engineering},
  title     = {Goal-Oriented Mesh Adaptivity for Fluid-Structure Interaction with Application to Heart-Valve Settings},
  year      = {2012},
  issn      = {0004-0738},
  month     = jan,
  number    = {1},
  pages     = {73--99},
  volume    = {59},
  doi       = {10.2478/v10180-012-0005-2},
  publisher = {Polish Academy of Sciences Chancellery},
}

@PhdThesis{Failerdiss,
  author   = {Failer, Lukas},
  school   = {Technische Universität München},
  title    = {Optimal Control of Time-Dependent Nonlinear Fluid-Structure Interaction},
  year     = {2017},
  language = {en},
  pages    = {151},
  url      = {https://mediatum.ub.tum.de/1355578},
}

@Article{Od18,
 Author = {Oden, J. Tinsley},
 Title = {Adaptive multiscale predictive modelling},
 FJournal = {Acta Numerica},
 Journal = {Acta Numerica},
 ISSN = {0962-4929},
 Volume = {27},
 Pages = {353--450},
 Year = {2018},
 Language = {English},
 DOI = {10.1017/S096249291800003X},
 zbMATH = {7099158},
 Zbl = {1430.60006}
}

@misc{Duprez2019,
author = "Michel Duprez and Stéphane Bordas and Marek Bucki and Huu Phuoc Bui and Franz Chouly and Vanessa Lleras and Claudio Lobos and Alexei Lozinski and Pierre-Yves Rohan and Satyendra Tomar",
title = "{Quantifying discretization errors for soft-tissue simulation in computer assisted surgery: a preliminary study}",
year = "2019",
mmonth = "5",
url = "https://figshare.com/articles/software/Quantifying_discretization_errors_for_soft-tissue_simulation_in_computer_assisted_surgery_a_preliminary_study/8128178",
doi = "10.6084/m9.figshare.8128178.v1",
note = {Figshare. DOI:\url{10.6084/m9.figshare.8128178.v1}} 
}

@article{bui:hal-04208610,
  TITLE = {{Enhancing biomechanical simulations based on a posteriori error estimates: the potential of Dual Weighted Residual-driven adaptive mesh refinement}},
  AUTHOR = {Bui, Huu Phuoc and Duprez, Michel and Rohan, Pierre-Yves and Lejeune, Arnaud and Bordas, St{\'e}phane Pierre Alain and Bucki, Marek and Chouly, Franz},
  journal = {International Journal for Numerical Methods in Biomedical Engineering},
url = {https://doi.org/10.1002/cnm.3897},
pages = {e3897--n/a},
year = {2025},
note = {e3897 cnm.3897},
}

@article {suttmeier2004,
    AUTHOR = {Suttmeier, F. T.},
     TITLE = {Reliable approximation of weight factors entering
              residual-based error bounds for {FE}-discretisations},
   JOURNAL = {Computing},
  FJOURNAL = {Computing. Archives for Scientific Computing},
    VOLUME = {73},
      YEAR = {2004},
    NUMBER = {3},
     PAGES = {199--205},
      ISSN = {0010-485X,1436-5057},
   MRCLASS = {65N30 (35J25 65N15)},
  MRNUMBER = {2106248},
       DOI = {10.1007/s00607-004-0082-2},
       URL = {https://doi.org/10.1007/s00607-004-0082-2},
}

@article{wick2024,
  title={A Posteriori Single- and Multi-Goal Error Control and Adaptivity for Partial Differential Equations},
  author={Endtmayer, Bernhard and Langer, Ulrich and Richter, Thomas and Schafelner, Andreas and Wick, Thomas},
  journal={Advances in Applied Mechanics (AAMS)},
  volume={59},
  pages={19-108},
  year={2024},
  publisher={Elsevier}
}

@Article{wick2017pum,
 Author = {Endtmayer, Bernhard and Wick, Thomas},
 Title = {A partition-of-unity dual-weighted residual approach for multi-objective goal functional error estimation applied to elliptic problems},
 Journal = {Computational Methods in Applied Mathematics},
 aJournal = {Comput. Methods Appl. Math.},
 ISSN = {1609-4840},
 Volume = {17},
 Number = {4},
 Pages = {575--599},
 Year = {2017},
 Language = {English},
 DOI = {10.1515/cmam-2017-0001},
 zbMATH = {7174725},
 Zbl = {1434.65252}
}

@Article{endtmayer2021,
 Author = {Endtmayer, Bernhard and Langer, Ulrich and Wick, Thomas},
 Title = {Reliability and efficiency of {DWR}-type a posteriori error estimates with smart sensitivity weight recovering},
 Journal = {Computational Methods in Applied Mathematics},
 aJournal = {Comput. Methods Appl. Math.},
 ISSN = {1609-4840},
 Volume = {21},
 Number = {2},
 Pages = {351--371},
 Year = {2021},
 Language = {English},
 DOI = {10.1515/cmam-2020-0036},
 zbMATH = {7433985},
 Zbl = {1476.65297}
}

@Book{ainsworth-oden-2000,
  Title                    = {A posteriori error estimation in finite element analysis},
  Author                   = {Ainsworth, Mark and Oden, J. Tinsley},
  Publisher                = {Wiley-Interscience, New York},
  Year                     = {2000},
  Series                   = {Pure and Applied Mathematics},

  ddoi                      = {10.1002/9781118032824},
  ISBN                     = {0-471-29411-X},
  Mrclass                  = {65-02 (65N15)},
  Mrnumber                 = {1885308},
  Mrreviewer               = {Ricardo G. Dur{\'a}n},
  Pages                    = {xx+240},
  uurl                      = {http://dx.doi.org/10.1002/9781118032824}
}

@Article{becker-rannacher-2001,
  Title                    = {An optimal control approach to a posteriori error estimation
 in finite element methods},
  Author                   = {Becker, Roland and Rannacher, Rolf},
  aJournal                  = {Acta Numer.},
  Year                     = {2001},
  Pages                    = {1--102},
  Volume                   = {10},

  DDoi                      = {10.1017/S0962492901000010},
  journal                 = {Acta Numerica},
  ISSN                     = {0962-4929},
  Mrclass                  = {65N15 (49M25 65N30)},
  Mrnumber                 = {2009692},
  Mrreviewer               = {Wen Bin Liu},
  uurl                      = {http://dx.doi.org/10.1017/S0962492901000010}
}

@Article{becker-rannacher-1996,
  Title                    = {A feed-back approach to error control in finite element
 methods: basic analysis and examples},
  Author                   = {Becker, R. and Rannacher, R.},
  aJournal                  = {East-West J. Numer. Math.},
  Year                     = {1996},
  Number                   = {4},
  Pages                    = {237--264},
  Volume                   = {4},

  journal                 = {East-West Journal of Numerical Mathematics},
  ISSN                     = {0928-0200},
  Mrclass                  = {65N30 (65M60)},
  Mrnumber                 = {1430239},
	uurl = {https://pdfs.semanticscholar.org/a85c/d5dbe857581bac6e6f0317194650448106f9.pdf}
}

@article{bijar2016atlas,
  title={Atlas-Based Automatic Generation of Subject-Specific Finite Element Tongue Meshes},
  author={Bijar, Ahmad and Rohan, Pierre-Yves and Perrier, Pascal and Payan, Yohan},
  journal={Annals of Biomedical Engineering},
  volume={44},
  number={1},
  pages={16--34},
  year={2016},
  publisher={Springer},
	uurl={https://doi.org/10.1007/s10439-015-1497-y}
}

@Article{bucki-2011,
  Title                    = {Jacobian-based Repair Method for Finite Element Meshes after Registration},
  Author                   = {M. Bucki and C. Lobos and Y. Payan and N. Hitschfeld},
  Journal                  = {Engineering with Computers},
  Year                     = {2011},
  Number                   = {3},
  Pages                    = {285--297},
  Volume                   = {27},

  ddoi                      = {10.1007/s00366-010-0198-2},
  uurl                      = {http://dx.doi.org/10.1007/s00366-010-0198-2}
}

@book{cowin2001cardiovascular,
  title={Cardiovascular soft tissue mechanics},
  author={Cowin, Stephen Corteen and Humphrey, Jay Dowell},
  year={2001},
  publisher={Springer}
}

@Article{dorfler1996,
  Title                    = {A convergent adaptive algorithm for {P}oisson's equation},
  Author                   = {D{\"o}rfler, Willy},
  aJournal                  = {SIAM J. Numer. Anal.},
  Year                     = {1996},
  Number                   = {3},
  Pages                    = {1106--1124},
  Volume                   = {33},

  Coden                    = {SJNAAM},
  ddoi                      = {10.1137/0733054},
  journal                 = {SIAM Journal on Numerical Analysis},
  ISSN                     = {0036-1429},
  Mrclass                  = {65N50 (65N55)},
  Mrnumber                 = {1393904},
  Mrreviewer               = {S. F. McCormick},
  uurl                      = {http://dx.doi.org/10.1137/0733054}
}

@Book{ern_guermond2004,
  Title                    = {Theory and practice of finite elements},
  Author                   = {Ern, Alexandre and Guermond, Jean-Luc},
  Publisher                = {Springer-Verlag, New York},
  Year                     = {2004},
  Series                   = {Applied Mathematical Sciences},
  Volume                   = {159},

  ddoi                      = {10.1007/978-1-4757-4355-5},
  ISBN                     = {0-387-20574-8},
  Mrclass                  = {65-02 (65M60 65N30 74S05 76M10 78M10)},
  Mrnumber                 = {2050138},
  Mrreviewer               = {R. S. Anderssen},
  Pages                    = {xiv+524},
  uurl                      = {http://dx.doi.org/10.1007/978-1-4757-4355-5}
}

@article{meunier2008mechanical,
  title={Mechanical experimental characterisation and numerical modelling of an unfilled silicone rubber},
  author={Meunier, Luc and Chagnon, Gr{\'e}gory and Favier, Denis and Org{\'e}as, Laurent and Vacher, Pierre},
  journal={Polymer Testing},
  volume={27},
  number={6},
  pages={765--777},
  year={2008},
  publisher={Elsevier}
}

@InCollection{nochetto-2009,
  Title                    = {Theory of adaptive finite element methods: an introduction},
  Author                   = {Nochetto, Ricardo H. and Siebert, Kunibert G. and Veeser,
 Andreas},
  Booktitle                = {Multiscale, nonlinear and adaptive approximation},
  Publisher                = {Springer, Berlin},
  Year                     = {2009},
  Pages                    = {409--542},

  DDoi                      = {10.1007/978-3-642-03413-8_12},
  Mrclass                  = {65N30 (65N15 65N50)},
  Mrnumber                 = {2648380},
  Mrreviewer               = {Nicolae Pop},
  uurl                      = {http://dx.doi.org/10.1007/978-3-642-03413-8_12}
}

@Book{payan-ohayon-2017,
  Title                    = {Biomechanics of living organs: hyperelastic constitutive laws for finite element modeling},
  Author                   = {Payan, Yohan and Ohayon, Jacques},
	series = {{Academic Press Series in Biomedical Engineering}},
  Publisher                = {Elsevier},
  Year                     = {2017},
	uurl = {https://doi.org/10.1016/c2015-0-00832-2}
}

@Article{rognes-logg-2013,
  Title                    = {Automated goal-oriented error control {I}: {S}tationary
 variational problems},
  Author                   = {Rognes, Marie E. and Logg, Anders},
  AJournal                  = {SIAM J. Sci. Comput.},
  Year                     = {2013},
  Number                   = {3},
  Pages                    = {C173--C193},
  Volume                   = {35},

  DDoi                      = {10.1137/10081962X},
  Journal                 = {SIAM Journal on Scientific Computing},
  ISSN                     = {1064-8275},
  Mrclass                  = {65N30 (65N50)},
  Mrnumber                 = {3048221},
  Mrreviewer               = {JaEun Ku},
  uurl                      = {http://dx.doi.org/10.1137/10081962X}
}

@Book{verfurth-2013,
  Title                    = {A posteriori error estimation techniques for finite element
 methods},
  Author                   = {Verf{\"u}rth, R{\"u}diger},
  Publisher                = {Oxford University Press, Oxford},
  Year                     = {2013},
  SSeries                   = {Numerical Mathematics and Scientific Computation},

  ddoi                      = {10.1093/acprof:oso/9780199679423.001.0001},
  ISBN                     = {978-0-19-967942-3},
  Mrclass                  = {65N30 (35J25 35K20 35Q30 35Q74 65N15)},
  Mrnumber                 = {3059294},
  Mrreviewer               = {Manfred Dobrowolski},
  Pages                    = {xx+393},
  uurl                      = {http://dx.doi.org/10.1093/acprof:oso/9780199679423.001.0001}
}

@article{Carstensen20141195,
title = "Axioms of adaptivity ",
journal = "Computers and Mathematics with Applications ",
volume = "67",
number = "6",
pages = "1195 - 1253",
year = "2014",
note = "",
issn = "0898-1221",
doi = "http://doi.org/10.1016/j.camwa.2013.12.003",
url = "http://www.sciencedirect.com/science/article/pii/S0898122113006822",
author = "C. Carstensen and M. Feischl and M. Page and D. Praetorius"
}

@article{BringmannBrunnerPraetoriusStreitberger+2024,
	author = {Philipp Bringmann and Maximilian Brunner and Dirk Praetorius and Julian Streitberger},
	doi = {10.1515/jnma-2023-0150},
	journal = {Journal of Numerical Mathematics},
	lastchecked = {2025-03-14},
	title = {{Optimal complexity of goal-oriented adaptive FEM for nonsymmetric linear elliptic PDEs}},
	title2 = {},
	url = {https://doi.org/10.1515/jnma-2023-0150},
	year = {2024}
}

@article{BeBrInMePr23,
 author = {Becker, Roland and Brunner, Maximilian and Innerberger, Michael and Melenk, Jens Markus and Praetorius, Dirk},
 title = {Cost-optimal adaptive iterative linearized {FEM} for semilinear elliptic {PDEs}},
 journal = {European Series in Applied and Industrial Mathematics (ESAIM): Mathematical Modelling and Numerical Analysis},
 qjournal = {ESAIM, Math. Model. Numer. Anal.},
 issn = {0764-583X},
 volume = {57},
 number = {4},
 pages = {2193--2225},
 year = {2023},
 language = {English},
 doi = {10.1051/m2an/2023036},
 zbMATH = {7739210},
 Zbl = {1523.65090}
}

@article {giles-suli-2002,
    AUTHOR = {Giles, Michael B. and S\"uli, Endre},
     TITLE = {Adjoint methods for {PDE}s: a posteriori error analysis and
              postprocessing by duality},
  AJOURNAL = {Acta Numer.},
  JOURNAL = {Acta Numerica},
    VOLUME = {11},
      YEAR = {2002},
     PAGES = {145--236},
      ISSN = {0962-4929},
   MRCLASS = {65N15 (35G15 35J25 65N30)},
  MRNUMBER = {2009374},
MRREVIEWER = {Lucas J\'odar},
       dddoi = {10.1017/S096249290200003X},
       uurl = {http://dx.doi.org/10.1017/S096249290200003X},
}

@article {maday-patera-2000,
    AUTHOR = {Maday, Yvon and Patera, Anthony T.},
     TITLE = {Numerical analysis of a posteriori finite element bounds for
              linear functional outputs},
  AJOURNAL = {Math. Models Methods Appl. Sci.},
  JOURNAL = {Mathematical Models and Methods in Applied Sciences},
    VOLUME = {10},
      YEAR = {2000},
    NUMBER = {5},
     PAGES = {785--799},
      ISSN = {0218-2025},
   MRCLASS = {65M60 (65M12)},
  MRNUMBER = {1763202},
MRREVIEWER = {Erwin Schechter},
       dddoi = {10.1142/S0218202500000409},
       uurl = {http://dx.doi.org/10.1142/S0218202500000409},
}

@article {nochetto-veeser-verani-2009,
    AUTHOR = {Nochetto, Ricardo H. and Veeser, Andreas and Verani, Marco},
     TITLE = {A safeguarded dual weighted residual method},
  AJOURNAL = {IMA J. Numer. Anal.},
  JOURNAL = {IMA Journal of Numerical Analysis},
    VOLUME = {29},
      YEAR = {2009},
    NUMBER = {1},
     PAGES = {126--140},
      ISSN = {0272-4979},
   MRCLASS = {65N15},
  MRNUMBER = {2470943},
       dddoi = {10.1093/imanum/drm026},
       uurl = {http://dx.doi.org/10.1093/imanum/drm026},
}

@article {paraschivoiu-peraire-patera-1997,
    AUTHOR = {Paraschivoiu, Marius and Peraire, Jaime and Patera, Anthony
              T.},
     TITLE = {A posteriori finite element bounds for linear-functional
              outputs of elliptic partial differential equations},
      NOTE = {Symposium on Advances in Computational Mechanics, Vol. 2
              (Austin, TX, 1997)},
  AJOURNAL = {Comput. Methods Appl. Mech. Engrg.},
  JOURNAL = {Computer Methods in Applied Mechanics and Engineering},
    VOLUME = {150},
      YEAR = {1997},
    NUMBER = {1-4},
     PAGES = {289--312},
      ISSN = {0045-7825},
   MRCLASS = {65N30 (76M10)},
  MRNUMBER = {1487947},
MRREVIEWER = {Jacques Rappaz},
       dddoi = {10.1016/S0045-7825(97)00086-8},
       uurl = {http://dx.doi.org/10.1016/S0045-7825(97)00086-8},
}

@article {prudhomme-oden-1999,
    AUTHOR = {Prudhomme, S. and Oden, J. T.},
     TITLE = {On goal-oriented error estimation for elliptic problems:
              application to the control of pointwise errors},
      NOTE = {New advances in computational methods (Cachan, 1997)},
  AJOURNAL = {Comput. Methods Appl. Mech. Engrg.},
  JOURNAL = {Computer Methods in Applied Mechanics and Engineering},
    VOLUME = {176},
      YEAR = {1999},
    NUMBER = {1-4},
     PAGES = {313--331},
      ISSN = {0045-7825},
   MRCLASS = {65N30 (74S05)},
  MRNUMBER = {1665351},
       dddoi = {10.1016/S0045-7825(98)00343-0},
       uurl = {http://dx.doi.org/10.1016/S0045-7825(98)00343-0},
}

@article {MeiRi14,
    AUTHOR = {Meidner, Dominik and Richter, Thomas},
     TITLE = {Goal-oriented error estimation for the fractional step theta
              scheme},
   AJOURNAL = {Comput. Methods Appl. Math.},
  JOURNAL = {Computational Methods in Applied Mathematics},
    VOLUME = {14},
      YEAR = {2014},
    NUMBER = {2},
     PAGES = {203--230},
      ISSN = {1609-4840},
   MRCLASS = {65M60 (65M15 65M50)},
  MRNUMBER = {3187923},
MRREVIEWER = {Dami{\'a}n P. Ginestar},
       DOI = {10.1515/cmam-2014-0002},
       URL = {http://dx.doi.org/10.1515/cmam-2014-0002},
}

@article {MeiRi15,
    AUTHOR = {Meidner, Dominik and Richter, Thomas},
     TITLE = {A posteriori error estimation for the fractional step theta
              discretization of the incompressible {N}avier-{S}tokes
              equations},
   AJOURNAL = {Comput. Methods Appl. Mech. Engrg.},
  JOURNAL = {Computer Methods in Applied Mechanics and Engineering},
    VOLUME = {288},
      YEAR = {2015},
     PAGES = {45--59},
      ISSN = {0045-7825},
   MRCLASS = {65M60 (35Q30 65M15 76D05)},
  MRNUMBER = {3327016},
       DOI = {10.1016/j.cma.2014.11.031},
       URL = {http://dx.doi.org/10.1016/j.cma.2014.11.031},
}

@article {larsson-2002,
    AUTHOR = {Larsson, Fredrik and Hansbo, Peter and Runesson, Kenneth},
     TITLE = {Strategies for computing goal-oriented a posteriori error
              measures in non-linear elasticity},
   AJOURNAL = {Internat. J. Numer. Methods Engrg.},
  JOURNAL = {International Journal for Numerical Methods in Engineering},
    VOLUME = {55},
      YEAR = {2002},
    NUMBER = {8},
     PAGES = {879--894},
      ISSN = {0029-5981},
   MRCLASS = {74G15 (74B20)},
  MRNUMBER = {1928975},
       dddoi = {10.1002/nme.513},
       uurl = {http://dx.doi.org/10.1002/nme.513},
}

@article {whiteley-tavener-2014,
    AUTHOR = {Whiteley, J. P. and Tavener, S. J.},
     TITLE = {Error estimation and adaptivity for incompressible
              hyperelasticity},
  FJOURNAL = {Internat. J. Numer. Methods Engrg.},
  JOURNAL = {International Journal for Numerical Methods in Engineering},
    VOLUME = {99},
      YEAR = {2014},
    NUMBER = {5},
     PAGES = {313--332},
      ISSN = {0029-5981},
   MRCLASS = {65N30 (65N15 74B20)},
  MRNUMBER = {3233187},
       dddoi = {10.1002/nme.4677},
       uurl = {http://dx.doi.org/10.1002/nme.4677},
}

@article{wick2016goal,
  title={{Goal functional evaluations for phase-field fracture using PU-based DWR mesh adaptivity}},
  author={Wick, Thomas},
  journal={Computational Mechanics},
  volume={57},
  number={6},
  pages={1017--1035},
  year={2016},
  publisher={Springer},
	uurl = {https://doi.org/10.1007/s00466-016-1275-1},
}

@article{Wi14_fsi_eale_heart,
   author = {T. Wick},
   title = {Flapping and Contact {FSI} Computations with the 
  Fluid-Solid Interface-Tracking/Interface-Capturing Technique
  and Mesh Adaptivity},
  journal = {Computational Mechanics},   
   year = {2014},
   volume = {53},
   number = {1},
   pages = {29-43},

}

@article{floch_vulnerable_2009,
	title = {Vulnerable {atherosclerotic} {plaque} {elasticity} {reconstruction} {based} on a {segmentation}-{driven} {optimization} {procedure} {using} {strain} {measurements}: {theoretical} {framework}},
	volume = {28},
	issn = {0278-0062},
	shorttitle = {Vulnerable {Atherosclerotic} {Plaque} {Elasticity} {Reconstruction} {Based} on a {Segmentation}-{Driven} {Optimization} {Procedure} {Using} {Strain} {Measurements}},
	dddoi = {10.1109/TMI.2009.2012852},
	uurl = {https://doi.org/10.1109/TMI.2009.2012852},
	number = {7},
	journal = {IEEE Transactions on Medical Imaging},
	author = {Le Floc'h, S. and Ohayon, J. and Tracqui, P. and Finet, G. and Gharib, A. M. and Maurice, R. L. and Cloutier, G. and Pettigrew, R. I.},
	year = {2009},
	pages = {1126--1137},
}

@article{gratsch_posteriori_2005,
	title = {A posteriori error estimation techniques in practical finite element analysis},
	volume = {83},
	issn = {0045-7949},
	uurl = {http://www.sciencedirect.com/science/article/pii/S0045794904003165},
	dddoi = {10.1016/j.compstruc.2004.08.011},
	number = {4},
	journal = {Computers \& Structures},
	author = {Gr\"{a}tsch, Thomas and Bathe, Klaus-Jürgen},
	year = {2005},
	pages = {235--265},
}

@article{braack2003posteriori,
    AUTHOR = {Braack, Malte and Ern, Alexandre},
     TITLE = {A posteriori control of modeling errors and discretization
              errors},
   AJOURNAL = {Multiscale Model. Simul.},
  JOURNAL = {SIAM Journal on Multiscale Modeling \& Simulation},
    VOLUME = {1},
      YEAR = {2003},
    NUMBER = {2},
     PAGES = {221--238},
      ISSN = {1540-3459},
   MRCLASS = {65N15 (65N30 65N50)},
  MRNUMBER = {1990196},
       dddoi = {10.1137/S1540345902410482},
       uurl = {https://doi.org/10.1137/S1540345902410482},
}

@article{oden2002estimation,
  title={Estimation of modeling error in computational mechanics},
  author={Oden, J Tinsley and Prudhomme, Serge},
  journal={Journal of Computational Physics},
  volume={182},
  number={2},
  pages={496--515},
  year={2002},
  publisher={Elsevier},
	uurl={https://doi.org/10.1006/jcph.2002.7183}
}

@article{endtmayer2018two,
 author = {Endtmayer, B. and Langer, U. and Wick, T.},
 title = {Two-side a posteriori error estimates for the dual-weighted residual method},
 journal = {SIAM Journal on Scientific Computing},
 qjournal = {SIAM J. Sci. Comput.},
 issn = {1064-8275},
 volume = {42},
 number = {1},
 pages = {a371--a394},
 year = {2020},
 language = {English},
 doi = {10.1137/18M1227275},
 zbMATH = {7170181},
 Zbl = {1440.65200}
}

@article {Costa2001,
	author = {Costa, K. D. and Holmes, J. W. and Mcculloch, A. D.},
	title = {Modelling cardiac mechanical properties in three dimensions},
	volume = {359},
	number = {1783},
	pages = {1233--1250},
	year = {2001},
	doi = {10.1098/rsta.2001.0828},
	publisher = {The Royal Society},
	issn = {1364-503X},
	URL = {http://rsta.royalsocietypublishing.org/content/359/1783/1233},
	eprint = {http://rsta.royalsocietypublishing.org/content/359/1783/1233.full.pdf},
	journal = {Philosophical Transactions of the Royal Society of London A: Mathematical, Physical and Engineering Sciences}
}

@incollection{OHAYON2017,
title = "Chapter 2 - {H}yperelastic Models for Contractile Tissues: Application to Cardiovascular Mechanics",
editor = "Yohan Payan and Jacques Ohayon",
booktitle = "Biomechanics of Living Organs",
publisher = "Academic Press",
address = "Oxford",
pages = "31 - 58",
year = "2017",
isbn = "978-0-12-804009-6",
doi = "https://doi.org/10.1016/B978-0-12-804009-6.00002-X",
url = "http://www.sciencedirect.com/science/article/pii/B978012804009600002X",
author = "Jacques Ohayon and Davide Ambrosi and Jean-Louis Martiel"
}

@article{Bovendeerd1992,
title = "Dependence of local left ventricular wall mechanics on myocardial fiber orientation: {A} model study",
journal = "Journal of Biomechanics",
volume = "25",
number = "10",
pages = "1129 - 1140",
year = "1992",
issn = "0021-9290",
doi = "https://doi.org/10.1016/0021-9290(92)90069-D",
url = "http://www.sciencedirect.com/science/article/pii/002192909290069D",
author = "P.H.M. Bovendeerd and T. Arts and J.M. Huyghe and D.H. van Campen and R.S. Reneman"
}

@article{Smith2004,
  author={Nicolas P Smith},
  title={A computational study of the interaction between coronary blood flow and myocardial mechanics},
  journal={Physiological Measurement},
  volume={25},
  number={4},
  pages={863},
  url={http://stacks.iop.org/0967-3334/25/i=4/a=007},
  year={2004},
}

@article{Mulvany1979,
  title={The active tension-length curve of vascular smooth muscle related to its cellular components},
  author={Maria Mulvany and David M. Warshaw},
  journal={The Journal of General Physiology},
  year={1979},
  volume={74},
  pages={85 - 104}
}

@article{Rachev1997,
title = "Theoretical study of the effect of stress-dependent remodeling on arterial geometry under hypertensive conditions",
journal = "Journal of Biomechanics",
volume = "30",
number = "8",
pages = "819 - 827",
year = "1997",
issn = "0021-9290",
doi = "https://doi.org/10.1016/S0021-9290(97)00032-8",
url = "http://www.sciencedirect.com/science/article/pii/S0021929097000328",
author = "Alexander Rachev"
}

@article{Mielczarek2012,
author = {Bo{\.z}ena Mielczarek and Justyna Uzia{\l{l}}ko-Mydlikowska},
title ={Application of computer simulation modeling in the health care sector: a survey},
journal = {SIMULATION},
volume = {88},
number = {2},
pages = {197-216},
year = {2012},
doi = {10.1177/0037549710387802},

URL = { 
        https://doi.org/10.1177/0037549710387802
    
},
eprint = { 
        https://doi.org/10.1177/0037549710387802
    
}
}

@article{zienkiewicz1992superconvergent,
  title={The superconvergent patch recovery ({SPR}) and adaptive finite element refinement},
  author={Zienkiewicz, OC and Zhu, JZ},
  journal={Computer Methods in Applied Mechanics and Engineering},
  volume={101},
  number={1},
  pages={207--224},
  year={1992},
  publisher={Elsevier}
}

@ARTICLE{Bui2016_TBME, 
author={H. P. Bui and S. Tomar and H. Courtecuisse and S. Cotin and S. P. A. Bordas}, 
journal={IEEE Transactions on Biomedical Engineering}, 
title={Real-Time Error Control for Surgical Simulation}, 
year={2018}, 
volume={65}, 
number={3}, 
pages={596-607}, 
doi={10.1109/TBME.2017.2695587}, 
ISSN={0018-9294}, 
mmonth={March},}

@article {Bui2017_error_brain_ijnmbe,
author = {Bui, Huu Phuoc and Tomar, Satyendra and Courtecuisse, Hadrien and Audette, Michel and Cotin, Stéphane and Bordas, Stéphane P.A.},
title = {Controlling the error on target motion through real-time mesh adaptation: Applications to deep brain stimulation},
journal = {International Journal for Numerical Methods in Biomedical Engineering},
issn = {2040-7947},
url = {http://dx.doi.org/10.1002/cnm.2958},
doi = {10.1002/cnm.2958},
pages = {e2958--n/a},
year = {2017},
note = {e2958 cnm.2958},
}

@ARTICLE{Heuveline2003,
    author = {V. Heuveline and R. Rannacher},
    title = {Duality-Based Adaptivity in the Hp-Finite Element Method},
    journal = {Journal of Numerical Mathematics},
    year = {2003},
    volume = {11},
    pages = { 95-113}
}

@article {duprez2020,
    AUTHOR = {Duprez, Michel and Bordas, St\'{e}phane Pierre Alain and Bucki,
              Marek and Bui, Huu Phuoc and Chouly, Franz and Lleras, Vanessa
              and Lobos, Claudio and Lozinski, Alexei and Rohan, Pierre-Yves
              and Tomar, Satyendra},
     TITLE = {Quantifying discretization errors for soft tissue simulation
              in computer assisted surgery: a preliminary study},
AJOURNAL = {Appl. Math. Model.},
  JOURNAL = {Applied Mathematical Modelling},
    VOLUME = {77},
      YEAR = {2020},
    NUMBER = {part 1},
     PAGES = {709--723},
      ISSN = {0307-904X},
   MRCLASS = {92C50 (65N30 74L15)},
  MRNUMBER = {3994022},
       DOI = {10.1016/j.apm.2019.07.055},
       URL = {https://doi.org/10.1016/j.apm.2019.07.055},
}

@article {granzow2018,
    AUTHOR = {Granzow, Brian N. and Oberai, Assad A. and Shephard, Mark S.},
     TITLE = {Adjoint-based error estimation and mesh adaptation for
              stabilized finite deformation elasticity},
   AJOURNAL = {Comput. Methods Appl. Mech. Engrg.},
  JOURNAL = {Computer Methods in Applied Mechanics and Engineering},
    VOLUME = {337},
      YEAR = {2018},
     PAGES = {263--280},
      ISSN = {0045-7825},
   MRCLASS = {65N30 (65N12 65N15 74B20)},
  MRNUMBER = {3801780},
       DOI = {10.1016/j.cma.2018.03.035},
       URL = {https://doi.org/10.1016/j.cma.2018.03.035},
}

@misc{Duprez2023,
author = "Michel Duprez and Arnaud Lejeune and Franz Chouly and Stéphane Bordas and Huu Phuoc Bui",
title = "{DWR-hyperelastic-soft-tissue}",
year = "2023",
mmonth = "4",
url = "https://figshare.com/articles/software/DWR-hyperelastic-soft-tissue/22548634",
doi = "10.6084/m9.figshare.22548634.v1",
jjournal="figshare",
note = {Figshare. DOI:\url{10.6084/m9.figshare.22548634.v1}}
}

@article {gonzalezestrada2014,
    AUTHOR = {Gonz\'{a}lez-Estrada, O. A. and Nadal, E. and R\'{o}denas, J. J. and
              Kerfriden, P. and Bordas, S. P. A. and Fuenmayor, F. J.},
     TITLE = {Mesh adaptivity driven by goal-oriented locally equilibrated
              superconvergent patch recovery},
   AJOURNAL = {Comput. Mech.},
  JOURNAL = {Computational Mechanics},
    VOLUME = {53},
      YEAR = {2014},
    NUMBER = {5},
     PAGES = {957--976},
      ISSN = {0178-7675},
   MRCLASS = {74B05 (65N50 74S05)},
  MRNUMBER = {3188432},
       DOI = {10.1007/s00466-013-0942-8},
       URL = {https://doi.org/10.1007/s00466-013-0942-8},
}

@article {becker2005,
    AUTHOR = {Becker, Roland and Vexler, Boris},
     TITLE = {Mesh refinement and numerical sensitivity analysis for
              parameter calibration of partial differential equations},
   AJOURNAL = {J. Comput. Phys.},
  JOURNAL = {Journal of Computational Physics},
    VOLUME = {206},
      YEAR = {2005},
    NUMBER = {1},
     PAGES = {95--110},
      ISSN = {0021-9991},
   MRCLASS = {65N50},
  MRNUMBER = {2135836},
       DOI = {10.1016/j.jcp.2004.12.018},
       URL = {https://doi.org/10.1016/j.jcp.2004.12.018},
}

@ARTICLE{SchVe08,
   author = {Michael Schmich and Boris Vexler},
   title = {Adaptivity with Dynamic Meshes for Space-Time 
Finite Element Discretizations of Parabolic Equations},
   journal = {SIAM Journal on Scientific Computing},
   year = {2008},
   volume = {30},
   number = {1},
   pages = {369 - 393}
}

@article{RiWi15_dwr,
title = "Variational localizations of the dual weighted residual estimator ",
journal = "Journal of Computational and Applied Mathematics ",
volume = "279",
number = "0",
pages = "192 - 208",
year = "2015",
note = "",
author = "T. Richter and T. Wick"
}

@article{ThiWi24,
  author  = {Jan Philipp Thiele and Thomas Wick},
  title   = {Numerical Modeling and Open-Source Implementation of Variational {Partition-of-Unity} Localizations of Space-Time {Dual-Weighted Residual} Estimators for Parabolic Problems},
  journal = {Journal of Scientific Computing},
  year    = {2024},
  volume  = {99},
  number  = {25},
  OPTpages   = {237-264}
}

@article{RoThiKoeWi24,
  url = {https://doi.org/10.1515/cmam-2022-0200},
  title = {Tensor-{p}roduct {s}pace-{t}ime {g}oal-{o}riented {e}rror {c}ontrol and {a}daptivity {w}ith {P}artition-of-{U}nity {D}ual-{W}eighted {R}esiduals for {n}onstationary {f}low {p}roblems},
  author = {Julian Roth and Jan Philipp Thiele and Uwe Köcher and Thomas Wick},
  pages = {185--214},
  volume = {24},
  number = {1},
  ajournal = {Comput. Methods Appl. Math.},
  journal = {Computational Methods in Applied Mathematics},
  doi = {doi:10.1515/cmam-2022-0200},
  year = {2024},
  lastchecked = {2024-01-15}
}

@article{ENDTMAYER2024286,
title = {Goal-oriented adaptive space-time finite element methods for regularized parabolic p-Laplace problems},
journal = {Computers \& Mathematics with Applications},
volume = {167},
pages = {286-297},
year = {2024},
author = {B. Endtmayer and U. Langer and A. Schafelner}
}

@incollection {HaHou03,
    AUTHOR = {Hartmann, R. and Houston, P.},
     TITLE = {Goal-oriented a posteriori error estimation for multiple
              target functionals},
 BOOKTITLE = {Hyperbolic problems: theory, numerics, applications},
     PAGES = {579--588},
 PUBLISHER = {Springer, Berlin},
      YEAR = {2003},
   MRCLASS = {35A15 (49M15 65M15 65M60)},
  MRNUMBER = {2053207},
}

@article {Ha08,
    AUTHOR = {Hartmann, R.},
     TITLE = {Multitarget error estimation and adaptivity in aerodynamic
              flow simulations},
   AJOURNAL = {SIAM J. Sci. Comput.},
  JOURNAL = {SIAM Journal on Scientific Computing},
    VOLUME = {31},
      YEAR = {2008},
    NUMBER = {1},
     PAGES = {708--731},
      ISSN = {1064-8275},
   MRCLASS = {76M10 (65N15 76G25)},
  MRNUMBER = {2460796},
       DOI = {10.1137/070710962},
       URL = {https://doi.org/10.1137/070710962},
}

@article {EnLaWi18,
    AUTHOR = {Endtmayer, B. and Langer, U. and Wick, T.},
     TITLE = {Multigoal-Oriented Error Estimates for Non-linear Problems},
   AJOURNAL = {J. Numer. Math.},
  JOURNAL = {Journal of Numerical Mathematics},
    VOLUME = {27},
      publisher={De Gruyter},
      YEAR = {2019},
    NUMBER = {4},
     PAGES = {215--236},
      ISSN = {1609-4840},
   MRCLASS = {65N30 (35Q74 49M15 65M60)},
  MRNUMBER = {3709050},
       URL = {https://doi.org/10.1515/jnma-2018-0038},
}

@article{RaVih13b,
	Author = {Rolf Rannacher and Jevgeni Vihharev},
	Journal = {Journal of Numerical Mathematics},
	Number = {1},
	Pages = {23-61},
	Title = {Adaptive finite element analysis of nonlinear problems: balancing of discretization 
and iteration errors},
	Volume = {21},
	Year = {2013}
}

@ARTICLE{BeRa12,
   author = {Michael Besier and Rolf Rannacher},
   title = {{Goal-oriented space-time adaptivity in the finite element Galerkin method for the computation of nonstationary incompressible flow}},
   qjournal = {Int. J. Num. Meth. Fluids},
   journal = {International Journal of Numerical Methods in Fluids},
   year = {2012},
   volume = {70},
   number = {},
   pages = {1139-1166}
}

@article{FaiWi18,
title = "Adaptive time-step control for nonlinear fluid-structure interaction",
journal = "Journal of Computational Physics",
volume = "366",
pages = "448 - 477",
year = "2018",
author = "Lukas Failer and Thomas Wick"
}

@article{FrKnStWeWi25,
author = {Stefan Frei and Tobias Knoke and Marc C. Steinbach and Anne-Kathrin Wenske and Thomas Wick},
title = {{Numerical simulations of fully Eulerian fluid-structure contact interaction using a ghost-penalty cut finite element approach}},
journal = {Advances in Computational Science and Engineering},
year = {2025},
volume = {3},
pages = {74-94},
issn = {},
doi = {10.3934/acse.2025005},
url = {https://www.aimsciences.org/article/id/67e645751d15e67fbdaccb6d}
}

@book{BoGaNe14,
  title={Fluid-Structure Interaction and Biomedical Applications},
  author={Bodn{\'a}r, T. and Galdi, G.P. and Ne{\v{c}}asov{\'a}, {\v{S}}.},
  isbn={9783034808224},
  series={Advances in Mathematical Fluid Mechanics},
  url={https://books.google.de/books?id=gknPBAAAQBAJ},
  year={2014},
  publisher={Springer Basel}
}

@article{RazzaqDamanikHronOuazziTurek:2012,
author = {M.~Razzaq and H.~Damanik and J.~Hron and A.~Ouazzi and S.~Turek},
title = {{FEM multigrid techniques for fluid-structure interaction with application to hemodynamics}},
journal = {Applied Numerical Mathematics},
volume = {62},
number = {9},
pages = {1156--1170},
year = {2012}
}

@BOOK{BaTaTe13,
  title = {Computational Fluid-Structure Interaction: {M}ethods and {A}pplications},
  publisher = {Wiley},
  year = {2013},
  author = {Y. Bazilevs and K. Takizawa and T.E. Tezduyar},
  volume = {},
  series = {},
  address = {}
}

@BOOK{FrHoRiWiYa17,
   author = {S. Frei and B. Holm and T. Richter and T. Wick and H. Yang},
   title = {Fluid-structure interactions: Fluid-Structure Interaction: Modeling, Adaptive Discretisations and Solvers},
   year = {2017},
   publisher = {de Gruyter}
}

@article{JoLaWi19_fsi,
author = {Jodlbauer, D. and Langer, U. and Wick, T.},
title = {Parallel block-preconditioned monolithic solvers for fluid-structure interaction problems},
journal = {International Journal for Numerical Methods in Engineering},
ajournal = {Int. J. Num. Meth. Eng.},
volume = {117},
number = {6},
pages = {623-643},
year = {2019}
}

@ARTICLE{BaQuaQua08,
   author = {S. Badia and Q. Quaini and A. Quarteroni},
   title = {Splitting methods based on algebraic factorization
  for fluid-structure interaction},
   ajournal = {SIAM J. Sci. Comput.},
   journal = {SIAM Journal on Scientific Computing},
   year = {2008},
   volume = {30},
   number = {4},
   pages = {1778-1805}
}

@ARTICLE{BaQuaQua08b,
   author = {S. Badia and Q. Quaini and A. Quarteroni},
   title = {Modular vs. non-modular preconditioners for fluid-structure interaction
systems with large added-mass effect},
   ajournal = {Comput. Methods. Appl. Mech. Engrg.},
   journal = {Computer Methods in Applied Mechanics and Engineering},
   year = {2008},
   volume = {197},
   number = {},
   pages = {4216-4232}
}

@article{gee2011truly,
  title={Truly monolithic algebraic multigrid for fluid--structure interaction},
  author={Gee, MW and K{\"u}ttler, U and Wall, WA},
  journal={International Journal for Numerical Methods in Engineering},
  qjournal={Int. J. Numer. Meth. Engrg.},
  volume={85},
  number={8},
  pages={987--1016},
  year={2011},
  publisher={John Wiley \& Sons, Ltd.}
}

@article{Croelal10,
title = "Fluid-structure interaction simulation of aortic blood flow",
journal = "Computers and Fluids",
volume = "43",
number = "1",
pages = "46 - 57",
year = "2011",
issn = "0045-7930",
doi = "10.1016/j.compfluid.2010.11.032",
url = "http://www.sciencedirect.com/science/article/pii/S0045793010003452",
author = "Paolo Crosetto and Philippe Reymond and Simone Deparis and Dimitrios Kontaxakis and Nikolaos Stergiopulos and Alfio Quarteroni"
}

@article{CrDeFouQua11,
author = {P. Crosetto and S. Deparis and G. Fourestey and A. Quarteroni},
title = {Parallel algorithms for fluid-structure interaction problems in 
haemodynamics},
journal = {SIAM Journal on Scientific Computing},
volume = {33},
number = {4},
pages = {1598-1622},
year = {2011}
}

@article{AhEndtSteiWi22,
title = {Multigoal-oriented error estimation and mesh adaptivity for fluid–structure interaction},
journal = {Journal of Computational and Applied Mathematics},
volume = {412},
pages = {114315},
year = {2022},
issn = {0377-0427},
doi = {https://doi.org/10.1016/j.cam.2022.114315},
url = {https://www.sciencedirect.com/science/article/pii/S0377042722001340},
author = {K. Ahuja and B. Endtmayer and M.C. Steinbach and T. Wick}
}

@ARTICLE{StTeBe03,
   author = {K. Stein and T. Tezduyar and R. Benney},
   title = {Mesh moving techniques for 
fluid-structure interactions with large displacements},
   journal = {Journal of Applied Mechanics},
   year = {2003},
   volume = {70},
   number = {},
   pages = {58-63}
}

@inbook{DoHuePoRo04,
   author = {J. Donea and A. Huerta and J.-Ph. Ponthot and A. Rodriguez-Ferran},
   title = {Arbitrary {L}agrangian-{E}ulerian {m}ethods},
   pages = {1-25},
   publisher = {John Wiley and Sons},
   year = {2004},
   series = {Encyclopedia of Computational Mechanics},
   address = {}
}

@ARTICLE{CouShk05,
   author = {Daniel Coutand and Steve Shkoller},
   title = {Motion of an Elastic Solid inside an Incompressible Viscous Fluid},
   ajournal = {Arch. Rational Mech. Anal.},
journal = {Archive for Rational Mechanics and Analysis},
   year = {2005},
   volume = {},
   number = {},
   pages = {25-102}
}

@ARTICLE{CouShk06,
   author = {Daniel Coutand and Steve Shkoller},
   title = {The interaction between quasilinear elastodynamics and the 
{N}avier-{S}tokes equations},
   ajournal = {Arch. Rational Mech. Anal.},
   journal = {Archive for Rational Mechanics and Analysis},
   year = {2006},
   volume = {179},
   number = {},
   pages = {303-352}
}

@Article{dealII96,
  author  = {Pasquale C. Africa and Daniel Arndt and Wolfgang Bangerth and Bruno Blais and
             Marc Fehling and Rene Gassm{\"o}ller and Timo Heister and Luca Heltai and
             Sebastian Kinnewig and Martin Kronbichler and Matthias Maier and Peter Munch and
             Magdalena Schreter-Fleischhacker and Jan P. Thiele and Bruno Turcksin and
             David Wells and Vladimir Yushutin},
  title   = {The deal.II library, Version 9.6},
  journal = {Journal of Numerical Mathematics},
  year    = 2024,
  volume  = 32,
  number  = 4,
  pages   = {369--380},
  doi     = {10.1515/jnma-2024-0137}
}

@Article{deal2020,
  title   = {The {deal.II} finite element library: Design, features, and insights},
  author  = {Daniel Arndt and Wolfgang Bangerth and Denis Davydov and
             Timo Heister and Luca Heltai and Martin Kronbichler and
             Matthias Maier and Jean-Paul Pelteret and Bruno Turcksin and
             David Wells},
  journal = {Computers {\&} Mathematics with Applications},
  year    = {2021},
  DOI     = {10.1016/j.camwa.2020.02.022},
  pages   = {407-422},
  volume  = {81},
  issn    = {0898-1221},
  url     = {https://arxiv.org/abs/1910.13247}
}

@Article{RichterWick2013
,
  author  = {T. Richter and T. Wick},
  title   = {Optimal Control and Parameter Estimation for Stationary Fluid-Structure Interaction Problems},
  doi     = {10.1137/120893239},
  number  = {5},
  pages   = {B1085--B1104},
  volume  = {35},
  journal = {SIAM J. Sci. Comput.},
  year    = {2013},
}

@Article{FailerRichter2019opt,
  author = 	 {L. Failer and T. Richter},
  title = 	 {A Newton multigrid framework for optimal control of fluid-structure interactions},
  journal = 	 {Optimization and Engineering},
  year = 	 {2020},
  doi = {10.1007/s11081-020-09498-8},
}

@INPROCEEDINGS{Wi21_git_YIC,
   author = {T. Wick},
   title = {{Adjoint-based methods for optimization and goal-oriented error control applied to fluid-structure interaction: implementation of a partition-of-unity dual-weighted residual estimator for stationary forward FSI problems in deal.II}},
  booktitle = {Book of Extended Abstracts of the 6th ECCOMAS Young Investigators Conference
7th-9th July 2021, Valencia, Spain},
  year = {2021},
  editor = {},
  pages = {},
  publisher = {ECCOMAS},
  doi = {https://doi.org/10.4995/YIC2021.2021.12332},
  url = {https://doi.org/10.4995/YIC2021.2021.12332}
}

@article{ashraf2022fluid,
  title={Fluid--structure interaction modelling of the upper airway with and without obstructive sleep apnea: a review},
  author={Ashraf, Walid and Jacobson, Natasha and Popplewell, Neil and Moussavi, Zahra},
  journal={Medical \& Biological Engineering \& Computing},
  volume={60},
  number={7},
  pages={1827--1849},
  year={2022},
  publisher={Springer}
}

@InProceedings{fc2006,
author="Chouly, Franz
and Van Hirtum, Annemie
and Lagr{\'e}e, Pierre-Yves
and Paoli, Jean-Roch
and Pelorson, Xavier
and Payan, Yohan",
editor="Harders, Matthias
and Sz{\'e}kely, G{\'a}bor",
title="Simulation of the Retroglossal Fluid-Structure Interaction During Obstructive Sleep Apnea",
booktitle="Biomedical Simulation",
year="2006",
publisher="Springer Berlin Heidelberg",
address="Berlin, Heidelberg",
pages="48--57",
isbn="978-3-540-36010-0"
}

@Article{engmann2013,
author ="Engmann, Jan and Burbidge, Adam S.",
title  ="Fluid mechanics of eating{,} swallowing and digestion – overview and perspectives",
journal  ="Food Funct.",
year  ="2013",
volume  ="4",
issue  ="3",
pages  ="443-447",
publisher  ="The Royal Society of Chemistry",
doi  ="10.1039/C2FO30184A",
url  ="http://dx.doi.org/10.1039/C2FO30184A"}

@article{sundstrom2025,
title = {Fluid-Structure Interaction Analysis of Aerodynamic and Elasticity Forces During Vocal Fold Vibration},
journal = {Journal of Voice},
volume = {39},
number = {2},
pages = {293-303},
year = {2025},
issn = {0892-1997},
doi = {https://doi.org/10.1016/j.jvoice.2022.08.030},
url = {https://www.sciencedirect.com/science/article/pii/S0892199722002715},
author = {Elias Sundström and Liran Oren and Charles {Farbos de Luzan} and Ephraim Gutmark and Sid Khosla}
}

@article{mittal2013,
   author = "Mittal, Rajat and Erath, Byron D. and Plesniak, Michael W.",
   title = "Fluid Dynamics of Human Phonation and Speech", 
   journal= "Annual Review of Fluid Mechanics",
   year = "2013",
   volume = "45",
   number = "Volume 45, 2013",
   pages = "437-467",
   doi = "https://doi.org/10.1146/annurev-fluid-011212-140636",
   url = "https://www.annualreviews.org/content/journals/10.1146/annurev-fluid-011212-140636",
   publisher = "Annual Reviews",
   issn = "1545-4479",
   type = "Journal Article",
  }

@article{heil2011,
   author = "Heil, Matthias and Hazel, Andrew L.",
   title = "Fluid-Structure Interaction in Internal Physiological Flows", 
   journal= "Annual Review of Fluid Mechanics",
   year = "2011",
   volume = "43",
   number = "Volume 43, 2011",
   pages = "141-162",
   doi = "https://doi.org/10.1146/annurev-fluid-122109-160703",
   url = "https://www.annualreviews.org/content/journals/10.1146/annurev-fluid-122109-160703",
   publisher = "Annual Reviews",
   issn = "1545-4479",
   type = "Journal Article",
  }

@article{becker2022,
 author = {Becker, Roland and Innerberger, Michael and Praetorius, Dirk},
 title = {Adaptive {FEM} for parameter-errors in elliptic linear-quadratic parameter estimation problems},
 fjournal = {SIAM Journal on Numerical Analysis},
 journal = {SIAM J. Numer. Anal.},
 issn = {0036-1429},
 volume = {60},
 number = {3},
 pages = {1450--1471},
 year = {2022},
 language = {English},
 doi = {10.1137/21M1458077},
 zbMATH = {7554726},
 Zbl = {1514.65161}
}

@Article{Hsu2014,
author={Hsu, Ming-Chen
and Kamensky, David
and Bazilevs, Yuri
and Sacks, Michael S.
and Hughes, Thomas J. R.},
title={Fluid--structure interaction analysis of bioprosthetic heart valves: significance of arterial wall deformation},
journal={Computational Mechanics},
year={2014},
month={Oct},
day={01},
volume={54},
number={4},
pages={1055-1071},
issn={1432-0924},
doi={10.1007/s00466-014-1059-4},
url={https://doi.org/10.1007/s00466-014-1059-4}
}

@article{FENG2024116724,
	author = {Liuyang Feng and Hao Gao and Xiaoyu Luo},
	doi = {10.1016/j.cma.2023.116724},
	issn = {0045-7825},
	journal = {Computer Methods in Applied Mechanics and Engineering},
	pages = {116724},
	title = {Whole-heart modelling with valves in a fluid–structure interaction framework},
	url = {https://www.sciencedirect.com/science/article/pii/S0045782523008472},
	volume = {420},
	year = {2024}
}

@article{VANLOON2006806,
	author = {R. {van Loon} and P.D. Anderson and F.N. {van de Vosse}},
	doi = {10.1016/j.jcp.2006.01.032},
	issn = {0021-9991},
	journal = {Journal of Computational Physics},
	number = {2},
	pages = {806–823},
	title = {A fluid–structure interaction method with solid-rigid contact for heart valve dynamics},
	url = {https://www.sciencedirect.com/science/article/pii/S0021999106000416},
	volume = {217},
	year = {2006}
}

@article{Fuma2023,
	author = {Fumagalli, Ivan and Polidori, Rebecca and Renzi, Francesca and Fusini, Laura and Quarteroni, Alfio and Pontone, Gianluca and Vergara, Christian},
	journal = {International Journal for Numerical Methods in Biomedical Engineering},
	number = {6},
	pages = {e3704},
	title = {Fluid-structure interaction analysis of transcatheter aortic valve implantation},
	volume = {39},
	year = {2023}
}

@article{Hoffman2021,
	author = {{Hiromi Spühler}, Jeannette and Hoffman, Johan},
	journal = {International Journal for Numerical Methods in Engineering},
	number = {19},
	pages = {5258–5278},
	title = {An interface-tracking unified continuum model for fluid-structure interaction with topology change and full-friction contact with application to aortic valves},
	volume = {122},
	year = {2021}
}

\end{document}